\documentclass[12pt,leqno]{article}
\usepackage[mathscr]{eucal}
\usepackage{amsmath,amssymb,latexsym,theorem,bbm,amsfonts}
\usepackage{verbatim}
\usepackage{color}
\usepackage{appendix}
\usepackage{array,epsfig,theorem,bbm,graphics,booktabs,url}

\newcommand{\CC}{\mathbb{C}}

\newcommand{\NN}{\mathbb{N}}

\newcommand{\RR}{\mathbb{R}}

\newcommand{\ZZ}{\mathbb{Z}}

\newcommand{\tcB}{\widetilde{\cB}}

\newcommand{\tB}{\widetilde{B}}

\newcommand{\tF}{\widetilde{F}}

\newcommand{\tth}{\widetilde{h}}

\newcommand{\tcY}{\widetilde{\cY}}

\newcommand{\bzero}{{\boldsymbol{0}}}

\newcommand{\cA}{{\mathcal A}}
\newcommand{\cB}{{\mathcal B}}

\newcommand{\cD}{{\mathcal D}}

\newcommand{\cF}{{\mathcal F}}

\newcommand{\cN}{{\mathcal N}}

\newcommand{\cX}{{\mathcal X}}

\newcommand{\cY}{{\mathcal Y}}
\newcommand{\cZ}{{\mathcal Z}}
\newcommand{\cW}{{\mathcal W}}

\newcommand{\dd}{\mathrm{d}}
\newcommand{\ee}{\mathrm{e}}
\newcommand{\ff}{\mathrm{f}}
\newcommand{\ii}{\mathrm{i}}

\newcommand{\EE}{\operatorname{\mathbb{E}}}

\newcommand{\PP}{\operatorname{\mathbb{P}}}
\newcommand{\OO}{\operatorname{O}}

\newcommand{\var}{\operatorname{Var}}
\newcommand{\cov}{\operatorname{Cov}}
\newcommand{\sign}{\operatorname{sign}}

\renewcommand{\Re}{\operatorname{Re}}

\newcommand{\tS}{\widetilde{S}}

\newcommand{\hS}{\widehat{S}}
\newcommand{\tT}{\widetilde{T}}

\newcommand{\vare}{\varepsilon}

\renewcommand{\mid}{\,|\,}

\renewcommand{\leq}{\leqslant}
\renewcommand{\geq}{\geqslant}

\newcommand{\stoch}{\stackrel{\PP}{\longrightarrow}}
\newcommand{\distr}{\stackrel{\cD}{\longrightarrow}}
\newcommand{\distrf}{\stackrel{\cD_\ff}{\longrightarrow}}
\newcommand{\distre}{\stackrel{\cD}{=}}

\newcommand{\as}{\stackrel{{\mathrm{a.s.}}}{\longrightarrow}}

\newcommand{\bbone}{\mathbbm{1}}

\newcommand{\proofend}{\hfill\mbox{$\Box$}}

\numberwithin{equation}{section}

\theoremstyle{change} \theorembodyfont{\em}
\newtheorem{Lem}{Lemma.}[section]
\newtheorem{Thm}[Lem]{Theorem.}
\newtheorem{Pro}[Lem]{Proposition.}
\newtheorem{Cor}[Lem]{Corollary.}

\theorembodyfont{\rm}
\newtheorem{Rem}[Lem]{Remark.}

\begin{document}

\begin{center}
 {\bfseries\Large Iterated limits for aggregation of randomized}\\[1mm]
 {\bfseries\Large INAR(1) processes with Poisson innovations}

\vskip0.5cm

 {\sc\large
  M\'aty\'as $\text{Barczy}^{*,\diamond}$,
  Fanni $\text{Ned\'enyi}^{**}$,
  Gyula $\text{Pap}^{**}$}

\end{center}

\vskip0.2cm

\noindent
 * Faculty of Informatics, University of Debrecen,
   Pf.~12, H--4010 Debrecen, Hungary.

\noindent
 ** Bolyai Institute, University of Szeged,
     Aradi v\'ertan\'uk tere 1, H--6720 Szeged, Hungary.

\noindent e--mails: barczy.matyas@inf.unideb.hu (M. Barczy),
                    nfanni@math.u-szeged.hu (F. Ned\'enyi),
                    papgy@math.u-szeged.hu (G. Pap).

\noindent $\diamond$ Corresponding author.

\vskip0.2cm


\renewcommand{\thefootnote}{}
\vspace*{0.2cm}
\footnote{The research has been supported by the DAAD-M\"OB Research Grant
 No.~55757 partially financed by the German Federal Ministry of Education and
 Research (BMBF).
}

\vspace*{-10mm}

\begin{abstract}
We discuss joint temporal and contemporaneous aggregation of \ $N$ \ independent
 copies of strictly stationary INteger-valued AutoRegressive processes of
 order 1 (INAR(1)) with random coefficient \ $\alpha \in (0, 1)$ \ and with idiosyncratic Poisson
 innovations.
Assuming that \ $\alpha$ \ has a density function of the form
 \ $\psi(x) (1 - x)^\beta$, \ $x \in (0, 1)$, \ with
 \ $\lim_{x\uparrow 1} \psi(x) = \psi_1 \in (0, \infty)$, \ different limits of
 appropriately centered and scaled aggregated partial sums are shown to exist for
 \ $\beta \in (-1, 0)$, \ $\beta = 0$, \ $\beta \in (0, 1)$ \ or
 \ $\beta \in (1, \infty)$, \ when taking first the limit as \ $N \to \infty$ \ and then
 the time scale \ $n \to \infty$, \ or vice versa.
In fact, we give a partial solution to an open problem of Pilipauskait\.e and Surgailis
 \cite{PilSur} by replacing the random-coefficient AR(1) process with a certain randomized
 INAR(1) process.
\end{abstract}

\section{Introduction}

The aggregation problem is concerned with the relationship between individual (micro)
 behavior and aggregate (macro) statistics.
There exist different types of aggregation.
The scheme of contemporaneous (also called cross-sectional) aggregation of
 random-coefficient AR(1) models was firstly proposed by Robinson \cite{Rob} and Granger
 \cite{Gra} in order to obtain the long memory phenomena in aggregated time series.
See also Gon\c{c}alves and Gouri\'eroux \cite{GonGou}, Zaffaroni \cite{Zaf2004},
 Oppenheim and Viano \cite{OppVia2004}, Celov et al. \cite{CelLeiPhi} and Beran et al.
 \cite{BerSchGho} on the aggregation of more general time-series models with finite variance.
Puplinskait\.e and Surgailis \cite{PupSur1,PupSur2} discussed aggregation of
 random-coefficient AR(1) processes with infinite variance and innovations in the domain
 of attraction of a stable law.
Related problems for some network traffic models were studied in Willinger et al.\
 \cite{WilTaqSheWil}, Taqqu et al.\ \cite{TaqWilShe}, Gaigalas and Kaj \cite{GaiKaj} and
 Dombry and Kaj \cite{DomKaj}, where independent and centered ON/OFF processes are
 aggregated, in Mikosch et al.\ \cite{MikResRooSte}, where aggregation of
 \ $\mathrm{M}/\mathrm{G}/\infty$ \ queues with heavy-tailed activity periods are
 investigated, in Pipiras et al.\ \cite{PipTaqLev}, where integrated renewal or
 renewal-reward processes are considered, or in Igl\'oi and Terdik \cite{IglTer},
 where the limit behavior of the aggregate of certain random-coefficient Ornstein--Uhlenbeck
 processes is examined.
On page 512 in Jirak \cite{Jir} one can find a lot of references for papers dealing with the
 aggregation of continuous time stochastic processes.

The present paper extends some of the results in Pilipauskait\.e and Surgailis
 \cite{PilSur}, which discusses the limit behavior of sums
 \begin{equation}\label{SAR}
  S_t^{(N,n)} := \sum_{j=1}^N \sum_{k=1}^{\lfloor nt \rfloor} X_k^{(j)} ,
  \qquad t \in [0, \infty) ,
 \qquad N, n \in \{1, 2, \ldots\} ,
 \end{equation}
 where \ $(X_k^{(j)})_{k\in\{0,1,\ldots\}}$, \ $j \in \{1, 2, \ldots\}$, \  are
 independent copies of a stationary random-coefficient AR(1) process
 \begin{equation}\label{RCAR}
  X_k = a X_{k-1} + \vare_k , \qquad k \in \{1, 2, \ldots\} ,
 \end{equation}
 with standardized independent and identically distributed (i.i.d.) innovations
 \ $(\vare_k)_{k\in\{1,2,\ldots\}}$ \ having \ $\EE(\vare_1) = 0$ \ and
 \ $\var(\vare_1) = 1$, \ and a random coefficient \ $a$ \ with values in \ $[0, 1)$,
 \ being independent of \ $(\vare_k)_{k \in \{1,2,\ldots\}}$ \ and admitting a probability
 density function of the form
 \begin{equation}\label{varphi}
  \psi(x) (1 - x)^\beta , \qquad x \in [0, 1) ,
 \end{equation}
 where \ $\beta \in (-1, \infty)$ \ and \ $\psi$ \ is an integrable
 function on \ $[0, 1)$ \ having
 a limit \ $\lim_{x\uparrow 1} \psi(x) = \psi_1 > 0$.
\ Here the distribution of \ $X_0$ \ is chosen as the unique stationary distribution of
 the model \eqref{RCAR}.
Its existence was shown in Puplinskait\.e and Surgailis \cite[Proposition 1]{PupSur1}.
We point out that they considered so-called idiosyncratic innovations, i.e., the
 innovations \ $(\vare^{(j)}_k)_{k\in\ZZ_+}$, \ $j \in \NN$, \ belonging  to
 \ $(X^{(j)}_k)_{k\in\ZZ_+}$, \ $j \in \NN$, \ are independent.
In \cite{PilSur} they derived scaling limits of the finite dimensional distributions of
 \ $(A_{N,n}^{-1} S_t^{(N,n)})_{t\in[0,\infty)}$, \ where \ $A_{N,n}$ \ are some scaling
 factors and first \ $N \to \infty$ \ and then \ $n \to \infty$, \ or vice versa, or
 both \ $N$ \ and \ $n$ \ increase to infinity, possibly with different rates.
Very recently,  Pilipauskait\.e and Surgailis \cite{PilSur2} extended their results
 in \cite{PilSur} from the case of idiosyncratic innovations to the case of common
 innovations, i.e., when \ $(\vare^{(j)}_k)_{k\in\ZZ_+}=(\vare^{(1)}_k)_{k\in\ZZ_+}$, \ $j \in \NN$.

The aim of the present paper is to extend the results of Pilipauskait\.e and Surgailis
 \cite[Theorem 2.1]{PilSur} concerning iterated scaling limits to the case of
 certain randomized first-order Integer-valued AutoRegressive (INAR(1)) processes.
The theory and application of integer-valued time series models are rapidly developing and
 important topics, see, e.g., Steutel and van Harn \cite{SteHar} and Wei{\ss}
 \cite{Wei}.
The INAR(1) process is among the most fertile integer-valued time series models, and
 it was first introduced by McKenzie \cite{McK} and Al-Osh and Alzaid \cite{AloAlz}.
An INAR(1) time series model is a stochastic process \ $(X_k)_{k\in\{0,1,\ldots\}}$
 \ satisfying the recursive equation
 \begin{equation}\label{INAR1}
  X_k = \sum_{j=1}^{X_{k-1}} \xi_{k,j} + \vare_k, \qquad k \in \{1, 2, \ldots\} ,
 \end{equation}
 where \ $(\vare_k)_{k\in\{1,2,\ldots\}}$ \ are i.i.d.\ non-negative integer-valued
 random variables, \ $(\xi_{k,j})_{k,j\in\{1,2,\ldots\}}$ \ are i.i.d.\ Bernoulli random
 variables with mean \ $\alpha \in [0, 1]$, \ and \ $X_0$ \ is a non-negative
 integer-valued random variable such that \ $X_0$,
 \ $(\xi_{k,j})_{k,j\in\{1,2,\ldots\}}$ \ and \ $(\vare_k)_{k\in\{1,2,\ldots\}}$ \ are
 independent.
By using the binomial thinning operator \ $\alpha\,\circ$ \ due to Steutel and van Harn
 \cite{SteHar},
the INAR(1) model in \eqref{INAR1} can be written as
 \begin{equation}\label{SteHarINAR1}
  X_k = \alpha \circ X_{k-1} + \vare_k , \qquad k \in \{1, 2, \ldots\} ,
 \end{equation}
which form captures the resemblance with the AR model.
We note that an INAR(1) process can also be considered as a special branching process
 with immigration having Bernoulli offspring distribution.

Leonenko et al.\ \cite{LeoSavZhi} introduced the aggregation
 \ $\sum_{j=1}^\infty X^{(j)}$ \ of a sequence of independent stationary INAR(1)
 processes \ $X^{(j)}$, \ $j \in \NN$, \ where
 \ $X^{(j)}_k = \alpha^{(j)}\circ X^{(j)}_{k-1} + \vare^{(j)}_k$, \ $k \in \ZZ$,
 \ $j \in \NN$.
\ Under appropriate conditions on \ $\alpha^{(j)}$, \ $j \in \NN$, \ and on the
 distributions of \ $\vare^{(j)}$, \ $j \in \NN$, \ they showed that the process
 \ $\sum_{j=1}^\infty X^{(j)}$ \ is well-defined in \ $L^2$-sense and it has long memory.

We will consider a certain randomized INAR(1) process \ $(X_k)_{k\in\ZZ_+}$ \ with randomized thinning parameter
 \ $\alpha$, \ given formally by the recursive equation
 \begin{equation}\label{randomINAR1}
  X_k = \alpha \circ X_{k-1} + \vare_k , \qquad k \in \{1, 2, \ldots\} ,
 \end{equation}
 where \ $\alpha$ \ is a random variable with values in \ $(0, 1)$ \  and \ $X_0$ \ is
 some appropriate random variable.
\ This means that, conditionally on \ $\alpha$, \ the process \ $(X_k)_{k\in\ZZ_+}$ \ is an INAR(1)
 process with thinning parameter \ $\alpha$.
\ Conditionally on \ $\alpha$, \ the i.i.d.\ innovations
 \ $(\vare_k)_{k\in\{1,2,\ldots\}}$ \ are supposed to have a Poisson distribution with
 parameter \ $\lambda \in (0, \infty)$, \ and the conditional distribution of the initial
 value \ $X_0$ \ given \ $\alpha$ \ is supposed to be the unique stationary distribution,
 namely, a Poisson distribution with parameter \ $\lambda/(1-\alpha)$.
\ For a rigorous construction of this process, see Section \ref{RINAR}.
Here we only note that \ $(X_k)_{k\in\ZZ_+}$ \ is a strictly stationary sequence, but it
 is not even a Markov chain (so it is not an INAR(1) process) if \ $\alpha$ \ is not
 degenerate, see Appendix \ref{Markov}.
Let us also remark that the choice of Poisson-distributed innovations serves a technical purpose.
It allows us to calculate and use the explicit stationary distribution and the joint generator function given in \eqref{help1_alter}.
The authors are planning to try releasing this assumption and giving more general results in future research.

Note that there is another way of randomizing the INAR(1) model \eqref{SteHarINAR1},
 the so-called random-coefficient INAR(1) process (RCINAR(1)), proposed by Zheng et
 al.\ \cite{ZheBasDat2007} and Leonenko et al.\ \cite{LeoSavZhi}.
It differs from \eqref{randomINAR1}, namely, it is a process
 formally given by the recursive equation
 \[
  X_k = \alpha_k \circ X_{k-1} + \vare_k , \qquad k \in \{1, 2, \ldots\} ,
 \]
 where \ $(\alpha_k)_{k\in\{1, 2, \ldots\}}$ \ is an i.i.d.\ sequence of random
 variables with values in \ $[0, 1]$.
\ An RCINAR(1) process can be considered as a special kind of branching processes with
 immigration in a random environment, see Key \cite{Key}, where a rigorous construction
 is given on the state space of the so-called genealogical trees.

In the paper first we examine a strictly stationary INAR(1) process
 \eqref{SteHarINAR1} with deterministic thinning and Poisson innovation, and in
 Section \ref{prel} an explicit formula is given for the joint generator function of
 \ $(X_0, X_1, \ldots, X_k)$, \ $k \in \{0, 1, \ldots\}$.
\ In Section \ref{INAR} we consider independent copies of this stationary INAR(1)
 process supposing idiosyncratic Poisson innovations.
Applying the natural centering by the expectation, in Propositions
 \ref{simple_aggregation}, \ref{simple_aggregation2} and in Theorem
 \ref{double_aggregation}, we derive scaling limits for the
 contemporaneously, the temporally and the joint temporally and
 contemporaneously aggregated processes, respectively.
In Section \ref{RINAR} first we give a construction of the stationary randomized
 INAR(1) process \eqref{randomINAR1}.
Considering independent copies of this randomized INAR(1) process, we
 discuss the limit behavior of the temporal and contemporaneous aggregation of
 these processes, both with centering by the expectation and by the conditional expectation,
 see Propositions \ref{simple_aggregation_random}--\ref{simple_aggregation_random4}.
Then, assuming that the distribution of \ $\alpha$ \ has the form \eqref{varphi}, we
 prove iterated limit theorems for the joint temporally and contemporaneously
 aggregated processes in case of both centralizations, see Theorems
 \ref{joint_aggregation_random}--\ref{joint_aggregation_random_5.5}.
As a consequence of our results, we formulate limit theorems
 with centering by the empirical mean as well, see Corollary \ref{Cor1}.
Note that we have separate results for the different ranges of \ $\beta$ \ (namely,
 \ $\beta \in (-1, 0)$, \ $\beta = 0$, \ $\beta \in (0, 1)$ \ and
 \ $\beta \in (1, \infty)$), \ the different orders of the iterations, and the different centralizations.
 The case \ $\beta = 1$ \ is not covered in this paper, nor in
 Pilipauskait\.e and Surgailis \cite{PilSur} for the random coefficient AR(1) processes.
We discuss this case for both models in Ned\'enyi and Pap \cite{NedPap}.
Section \ref{Proofs} contains the proofs.
In the appendices we discuss the non-Markov property of the randomized INAR(1) model,
 some approximations of the exponential function and some
 of its integrals, and an integral representation of the fractional Brownian motion
 due to Pilipauskait\.e and Surgailis \cite{PilSur}.
We consider three kinds of centralizations (by the conditional and the unconditional expectations
 and by the empirical mean).
In Pilipauskait\.e and Surgailis \cite{PilSur} centralization does not
 appear since they aggregate centered processes.
In Jirak \cite{Jir} the role of centralizations by the conditional and the unconditional expectations
 is discussed, where an asymptotic theory of aggregated linear processes is developed, and the limit
 distribution of a large class of linear and nonlinear functionals of such processes are
 determined.

All in all, we have similar limit theorems for randomized INAR(1) processes that
 Pilipauskait\.e and Surgailis \cite[Theorem 2.1]{PilSur} have for random coefficients
 AR(1) processes.
On page 1014, Pilipauskait\.e and Surgailis \cite{PilSur} formulated an open problem that
 concerns possible existence and description of limit distribution of the double sum
 \eqref{SAR} for general i.i.d.\ processes \ $(X^{(j)}_t)_{t\in\RR_+}$, \ $j \in \NN$.
\ We solve this open problem for some randomized INAR(1) processes.
Since INAR(1) processes are special branching processes with immigration, based on our
 results, later on, one may proceed with general branching processes with immigration.
The techniques of our proofs differ from those of Pilipauskait\.e and
 Surgailis \cite{PilSur} in many cases, for a somewhat detailed comparison,
 see the beginning of Section \ref{Proofs}.

\section{Generator function of finite-dimensional distributions of
         Galton--Watson branching processes with immigration}
\label{prel}

Let \ $\ZZ_+$, \ $\NN$, \ $\RR$, \ $\RR_+$, \ and \ $\CC$ \ denote the set of
 non-negative integers, positive integers, real numbers, non-negative real numbers,
 and complex numbers, respectively.
The Borel \ $\sigma$-field on \ $\RR$ \ is denoted by \ $\cB(\RR)$.
\ Every random variable in this section will be defined on a fixed probability
 space \ $(\Omega, \cA, \PP)$.

For each \ $k, j \in \ZZ_+$, \ the number of individuals in the \ $k^\mathrm{th}$
 \ generation will be denoted by \ $X_k$, \ the number of offsprings produced by
 the \ $j^\mathrm{th}$ \ individual belonging to the \ $(k-1)^\mathrm{th}$
 \ generation will be denoted by \ $\xi_{k,j}$, \ and the number of immigrants in
 the \ $k^\mathrm{th}$ \ generation will be denoted by \ $\vare_k$.
\ Then we have
 \[
  X_k = \sum_{j=1}^{X_{k-1}} \xi_{k,j} + \vare_k , \qquad k \in \NN ,
 \]
 where we define \ $\sum_{j=1}^0 := 0$.
\ Here \ $\bigl\{ X_0, \, \xi_{k,j}, \, \vare_k : k, j \in \NN \bigr\}$ \ are
 supposed to be independent nonnegative integer-valued random variables.
Moreover, \ $\{ \xi_{k,j} : k, j \in \NN \}$ \ and \ $\{ \vare_k : k \in \NN \}$
 \ are supposed to consist of identically distributed random variables,
 respectively.

Let us introduce the generator functions
 \[
   F_k(z) := \EE(z^{X_k}) , \quad k \in \ZZ_+ , \qquad
   G(z) := \EE(z^{\xi_{1,1}}) , \qquad
   H(z) := \EE(z^{\vare_1})
 \]
 for \ $z \in D := \{ z \in \CC : |z| \leq 1\}$.
\ First we observe that for each \ $k \in \NN$, \ the conditional generator
 function \ $\EE(z_k^{X_k} \mid X_{k-1})$ \ of \ $X_k$ \ given \ $X_{k-1}$
 \ takes the form
 \begin{equation}\label{square-1}
   \EE(z_k^{X_k} \mid X_{k-1})
   = \EE\Bigl(z_k^{\sum_{j=1}^{X_{k-1}} \xi_{k,j} + \vare_k}
               \,\Big|\, X_{k-1}\Bigr)
   = \EE(z_k^{\vare_k}) \prod_{j=1}^{X_{k-1}} \EE(z_k^{\xi_{k,j}})
    = H(z_k) \, G(z_k)^{X_{k-1}}
 \end{equation}
 for \ $z_k \in D$, \ where we define \ $\prod_{j=1}^0 := 1$.
\ The aim of the following discussion is to calculate the joint generator functions
 of the finite dimensional distributions of \ $(X_k)_{k\in\ZZ_+}$.
\ Using \eqref{square-1}, we also have the recursion
 \begin{align*}
  F_k(z)
  = \EE(\EE(z^{X_k} \mid X_{k-1}))
   = \EE(H(z) \, G(z)^{X_{k-1}})
  = H(z) \EE\bigl(G(z)^{X_{k-1}}\bigr)
   = H(z) \, F_{k-1}(G(z))
 \end{align*}
 for \ $z \in D$ \ and \ $k  \in \NN$.
\ Put \ $G_{(0)}(z) := z$ \ and \ $G_{(1)}(z) := G(z)$ \ for \ $z \in D$,
 \ and introduce the iterates \ $G_{(k+1)}(z) := G_{(k)}(G(z))$, \ $z \in D$,
 \ $k \in \NN$.
\ The above recursion yields
 \[
   F_k(z)
   = H(z) \, H(G(z)) \cdots H(G_{(k-1)}(z)) \, F_0(G_{(k)}(z))
   = F_0(G_{(k)}(z)) \prod_{j=0}^{k-1} H(G_{(j)}(z))
 \]
 for \ $z \in D$ \ and \ $k \in \NN$.
\ Supposing that \ $\EE(\xi_{1,1}) = G'(1-) < 1$, \ $0<\PP(\xi_{1,1} = 0) < 1$,
 \ $0<\PP(\xi_{1,1} = 1)$ \ and \ $0<\PP(\vare_1 = 0) < 1$, \ the Markov chain
 \ $(X_k)_{k\in\ZZ_+}$ \ is irreducible and aperiodic.
Further, it is ergodic (positive recurrent) if and only if
 \ $\sum_{\ell=1}^\infty \log(\ell) \PP(\vare_1 = \ell) < \infty$,
  \ and in this case the unique stationary distribution has the generator function
 \[
  \tF(z) = \prod_{j=0}^\infty H(G_{(j)}(z)) , \qquad z \in D ,
 \]
 see, e.g., Seneta \cite[Chapter 5]{Sen2} and Foster and Williamson \cite[Theorem, part (iii)]{FosWil}.

Consider the special case with Bernoulli offspring and Poisson immigration
 distributions, namely,
 \begin{align}\label{help_stac_INAR1_1}
  \begin{split}
  & \PP(\xi_{1,1} = 1) = \alpha = 1 - \PP(\xi_{1,1} = 0) , \\
  & \PP(\vare_1 = \ell)
     = \frac{\lambda^\ell}{\ell!} \ee^{-\lambda} , \qquad
       \ell \in \ZZ_+ ,
  \end{split}
 \end{align}
 with \ $\alpha \in (0, 1)$ \ and \ $\lambda \in (0, \infty)$.
\ With the special choices \eqref{help_stac_INAR1_1}, the Galton--Watson process
 with immigration \ $(X_k)_{k\in\ZZ_+}$ \ is an INAR(1) process with Poisson
 innovations.
Then
 \[
   G(z) = 1 - \alpha + \alpha z , \qquad
   H(z) = \sum_{\ell=0}^\infty \frac{z^\ell\lambda^\ell}{\ell!}
          \ee^{-\lambda}
        = \ee^{\lambda(z-1)} , \qquad z \in \CC ,
 \]
 hence
 \[
   G_{(j)}(z) = 1 - \alpha^j + \alpha^j z , \qquad z \in \CC , \quad \ j \in \NN .
 \]
Indeed, by induction, for all \ $j \in \ZZ_+$, \
 \[
    G_{(j+1)}(z) = G(G_{(j)}(z))
                 = \alpha G_{(j)}(z) + 1 - \alpha
                 = \alpha(1-\alpha^j + \alpha^j z) + 1 - \alpha
                 = 1 - \alpha^{j+1} + \alpha^{j+1} z .
 \]
Since \ $\EE(\xi_{1,1}) = G'(1-) = \alpha \in (0, 1)$, \ $\PP(\xi_{1,1} = 0) = 1 - \alpha \in(0,1)$,
 \ $\PP(\xi_{1,1} = 1) = \alpha >0$, \ $\PP(\vare_1=0) = \ee^{-\lambda}\in(0,1)$, \ and
 \begin{align*}
  \sum_{\ell=1}^\infty \log(\ell) \frac{\lambda^\ell}{\ell!}\ee^{-\lambda}
  \leq \sum_{\ell=1}^\infty \ell \frac{\lambda^\ell}{\ell!}\ee^{-\lambda}
  =\EE(\varepsilon_1)
  =\lambda < \infty ,
 \end{align*}
 the Markov chain \ $(X_k)_{k\in\ZZ_+}$ \ has a unique stationary distribution
 admitting a generator function of the form
 \[
  \tF(z)
  = \prod_{j=0}^\infty \ee^{\alpha^j\lambda(z-1)}
  = \ee^{(1-\alpha)^{-1}\lambda(z-1)} , \qquad z \in \CC ,
 \]
 thus it is a Poisson distribution with expectation \ $(1-\alpha)^{-1} \lambda$.

Suppose now that the initial distribution is a Poisson distribution with
 expectation \ $(1-\alpha)^{-1} \lambda$, \ hence the Markov chain
 \ $(X_k)_{k\in\ZZ_+}$ \ is strictly stationary and
 \begin{equation}\label{square0}
  F_0(z_0) = \EE(z_0^{X_0}) = \ee^{(1-\alpha)^{-1}\lambda(z_0-1)} ,
  \qquad z_0 \in \CC .
 \end{equation}
By induction, one can derive the following result, formulae for the joint
 generator function of \ $(X_0, X_1, \ldots, X_k)$, \ $k \in \ZZ_+$.

\begin{Pro}\label{Pro_egyuttes_mom_gen}
Under \eqref{help_stac_INAR1_1} and supposing that the distribution of
 \ $X_0$ \ is Poisson distribution with expectation \ $(1-\alpha)^{-1} \lambda$,
 \ the joint generator function of \ $(X_0, X_1, \ldots, X_k)$, \ $k \in \ZZ_+$,
 \ takes the form
 \begin{equation}\label{help1_alter}
  \begin{aligned}
   F_{0,\dots, k}(z_0, \dots, z_k)
   &:= \EE(z_0^{X_0} z_1^{X_1} \cdots z_k^{X_k}) \\
   &= \exp\left\{\frac{\lambda}{1-\alpha}
                 \sum_{0\leq i\leq j\leq k}
                  \alpha^{j-i} (z_i - 1) z_{i+1} \cdots z_{j-1} (z_j - 1)\right\}
  \end{aligned}
 \end{equation}
 for all \ $k \in \NN$ \ and \ $z_0, \ldots, z_k \in \CC$, \ where, for \ $i = j$,
 \ the term in the sum above is \ $z_i - 1$.
\ Alternatively, one can write up the joint generator function as
 \begin{align}\label{help1}
  F_{0,\ldots,k}(z_0, \ldots, z_k)
  = \exp\left\{\lambda
               \sum_{0\leq i\leq j\leq k}
                (1-\alpha)^{K_{i,j,k}} \alpha^{j-i}
                (z_i z_{i+1} \cdots z_j - 1)\right\} ,
 \end{align}
 where
 \[
   K_{i,j,k}
   := \begin{cases}
       -1 & \text{if \ $i = 0$ \ and \ $j = k$,} \\
       0 & \text{if \ $i = 0$ \ and \ $0 \leq j \leq k - 1$,} \\
       0 & \text{if \ $1 \leq i \leq k$ \ and \ $j = k$,} \\
       1 & \text{if \ $1 \leq i \leq j \leq k - 1$.} \\
      \end{cases}
 \]
\end{Pro}

\begin{Rem}
Under the conditions of Proposition \ref{Pro_egyuttes_mom_gen},
 the distribution of \ $(X_0, X_1)$ \ can be represented using independent Poisson distributed random variables.
Namely, if \ $U$, $V$ \ and \ $W$ \ are independent Poisson distributed random variables with
 parameters \ $\lambda (1-\alpha)^{-1} \alpha$, \ $\lambda$ \ and \ $\lambda$, \ respectively, then
 \ $(X_0, X_1) \distre (U+V, U+W)$. \
Indeed, for all \ $z_0, z_1 \in \CC$,
 \begin{align*}
  \EE(z_0^{U+V} z_1^{U+W})
  &= \EE( (z_0 z_1)^U z_0^V z_1^W )
   = \EE( (z_0 z_1)^U ) \EE(z_0^V) \EE( z_1^W ) \\
  &= \ee^{\lambda(1-\alpha)^{-1}\alpha(z_0z_1-1)} \ee^{\lambda(z_0 - 1)}
     \ee^{\lambda(z_1 - 1)}  ,
 \end{align*}
 as desired.
Further, note that formula \eqref{help1} shows that \ $(X_0, \ldots, X_k)$ \
 has a \ $(k+1)$-variate Poisson distribution, see, e.g., Johnson et al. \cite[(37.85)]{Johnson}.%
\proofend
\end{Rem}

\section{Iterated aggregation of INAR(1) processes with Poisson innovations}
\label{INAR}

Let \ $(X_k)_{k\in\ZZ_+}$ \ be an INAR(1) process with offspring and immigration
 distributions given in \eqref{help_stac_INAR1_1} and with initial distribution
 given in \eqref{square0}, hence the process is strictly stationary.
Let \ $X^{(j)} = (X^{(j)}_k)_{k\in\ZZ_+}$, \ $j \in \NN$, \ be a sequence of
 independent copies of the stationary INAR(1) process \ $(X_k)_{k\in\ZZ_+}$.

First we consider a simple aggregation procedure.
For each \ $N \in \NN$, \ consider the stochastic process
 \ $S^{(N)} = (S^{(N)}_k)_{k\in\ZZ_+}$ \ given by
 \begin{equation}\label{SN}
   S^{(N)}_k
   := \sum_{j=1}^N (X^{(j)}_k - \EE(X^{(j)}_k)) , \qquad k \in \ZZ_+ ,
 \end{equation}
 where \ $\EE(X^{(j)}_k) = \lambda (1-\alpha)^{-1}$, \ $k \in \ZZ_+$,
 \ $j \in \NN$, \ since the stationary distribution is Poisson with expectation
 \ $(1-\alpha)^{-1} \lambda$.
\ We will use \ $\distrf$ \ or \ $\cD_\ff\text{-}\hspace*{-1mm}\lim$ \ for the weak
 convergence of the finite dimensional distributions, and \ $\distr$ \ for the weak
 convergence of stochastic processes with sample
 paths in \ $D(\RR_+,\RR)$, \ where \ $D(\RR_+,\RR)$ \ denotes the space of real-valued c\`adl\`ag
 functions defined on \ $\RR_+$.
\ The almost sure convergence is denoted by \ $\as$.

\begin{Pro}\label{simple_aggregation}
We have
 \[
   N^{-\frac{1}{2}} S^{(N)} \distrf \cX \qquad \text{as \ $N \to \infty$,}
 \]
 where \ $\cX = (\cX_k)_{k\in\ZZ_+}$ \ is a stationary Gaussian process with zero
 mean and covariances
 \begin{align}\label{cov}
  \EE(\cX_0 \cX_k)
  = \cov(X_0, X_k)
  = \frac{\lambda\alpha^k}{1-\alpha} , \qquad k \in \ZZ_+ .
 \end{align}
\end{Pro}


\begin{Pro}\label{simple_aggregation2}
We have
 \[
   \biggl(n^{-\frac{1}{2}}
          \sum_{k=1}^{\lfloor nt \rfloor}
           S^{(1)}_k\biggr)_{t\in\RR_+}
   = \biggl(n^{-\frac{1}{2}}
            \sum_{k=1}^{\lfloor nt \rfloor}
             (X_k^{(1)} - \EE(X_k^{(1)}))\biggr)_{t\in\RR_+}
   \distr \frac{\sqrt{\lambda(1+\alpha)}}{1-\alpha} B
 \]
 as \ $n \to \infty$, \ where \ $B = (B_t)_{t\in\RR_+}$ \ is a standard Brownian
 motion.
\end{Pro}

Note that Propositions \ref{simple_aggregation} and \ref{simple_aggregation2} are
 about the scaling of the space-aggregated process \ $S^{(N)}$ \ and the
 time-aggregated process
 \ $\bigl(\sum_{k=1}^{\lfloor nt \rfloor} S^{(1)}_k\bigr)_{t\in\RR_+}$,
 \ respectively.

For each \ $N, n \in \NN$, \ consider the stochastic process
 \ $S^{(N,n)} = (S^{(N,n)}_t)_{t\in\RR_+}$ \ given by
 \begin{equation}\label{SNnt}
  S^{(N,n)}_t
  := \sum_{j=1}^N \sum_{k=1}^{\lfloor nt \rfloor} (X^{(j)}_k - \EE(X^{(j)}_k)) ,
  \qquad t \in \RR_+ .
 \end{equation}


\begin{Thm}\label{double_aggregation}
We have
 \[
   \cD_\ff\text{-}\hspace*{-1mm}\lim_{N\to\infty} \,
   \cD_\ff\text{-}\hspace*{-1mm}\lim_{n\to\infty} \,
    (nN)^{-\frac{1}{2}} S^{(N,n)}
   =
   \cD_\ff\text{-}\hspace*{-1mm}\lim_{n\to\infty} \,
   \cD_\ff\text{-}\hspace*{-1mm}\lim_{N\to\infty} \,
    (nN)^{-\frac{1}{2}} S^{(N,n)}
   = \frac{\sqrt{\lambda(1+\alpha)}}{1-\alpha} B ,
 \]
 where \ $B = (B_t)_{t\in\RR_+}$ \ is a standard Brownian motion.
\end{Thm}

\section{Iterated aggregation of randomized INAR(1) processes with
          Poisson innovations}
\label{RINAR}

Let \ $\lambda \in (0, \infty)$, \ and let \ $\PP_\alpha$ \ be a probability
 measure on \ $(0, 1)$.
\ Then there exist a probability space \ $(\Omega, \cA, \PP)$, \ a random variable
 \ $\alpha$ \ with distribution \ $\PP_\alpha$ \ and random variables
 \ $\{X_0, \, \xi_{k,j}, \, \vare_k : k, j \in \NN\}$, \ conditionally independent
 given \ $\alpha$ \ on \ $(\Omega, \cA, \PP)$ \ such that
 \begin{gather}\label{help_stac_RINAR1_1}
  \PP(\xi_{k,j} = 1 \mid \alpha)
  = \alpha = 1 - \PP(\xi_{k,j} = 0 \mid \alpha) , \qquad
  k, j \in \NN , \\
  \PP(\vare_k = \ell \mid \alpha)
  = \frac{\lambda^\ell}{\ell!} \ee^{-\lambda} , \qquad \ell \in \ZZ_+ ,
  \qquad k \in \NN , \label{help_stac_RINAR1_2} \\
  \PP(X_0 = \ell \mid \alpha)
  = \frac{\lambda^\ell}{\ell!(1-\alpha)^\ell} \ee^{-(1-\alpha)^{-1}\lambda} ,
  \qquad \ell \in \ZZ_+ . \label{help_stac_RINAR1_3}
 \end{gather}
(Note that the conditional distribution of \ $\vare_k$ \ does not depend on
 \ $\alpha$.)
\ Indeed, for each \ $n \in \NN$, \ by Ionescu Tulcea's theorem (see, e.g.,
 Shiryaev \cite[II. \S\,9, Theorem 2]{Shi}), there exist a probability space
 \ $(\Omega_n, \cA_n, \PP_n)$ \ and random variables \ $\alpha^{(n)}$,
 \ $X_0^{(n)}$, \ $\vare_k^{(n)}$ \ and \ $\xi_{k,j}^{(n)}$ \ for
 \ $k, j \in \{1, \ldots, n\}$ \ on \ $(\Omega_n, \cA_n, \PP_n)$ \ such that
 \begin{align*}
  &\PP_n(\text{$\alpha^{(n)} \in B$, \ $X_0^{(n)} = x_0$,
               \ $\vare_k^{(n)} = \ell_k$, \ $\xi_{k,j}^{(n)} = x_{k,j}$ \ for all
               \ $k, j \in \{1, \ldots, n\}$}) \\
  &= \int_B
      p_n\left(a, x_0, (\ell_k)_{k=1}^n, (x_{k,j})_{k,j=1}^n\right)
      \, \PP_\alpha(\dd a)
 \end{align*}
 for all \ $B \in \cB(\RR)$, \ $x_0 \in \ZZ_+$, \ $(\ell_k)_{k=1}^n \in \ZZ_+^n$,
 \ $(x_{k,j})_{k,j=1}^n \in \{0, 1\}^{n\times n}$, \ with
 \begin{align*}
  p_n\left(a, x_0, (\ell_k)_{k=1}^n, (x_{k,j})_{k,j=1}^n\right)
  := \frac{\lambda^{x_0}}{x_0!(1-a)^{x_0}} \ee^{-(1-a)^{-1}\lambda}
     \prod_{k=1}^n \frac{\lambda^{\ell_k}}{\ell_k!} \ee^{-\lambda}
     \prod_{k,j=1}^n a^{x_{k,j}} (1 - a)^{1-x_{k,j}} ,
 \end{align*}
 since the mapping
 \ $(0, 1) \ni a
    \mapsto p_n \left(a, x_0, (\ell_k)_{k=1}^n, (x_{k,j})_{k,j=1}^n\right)$
 \ is Borel measurable for all \ $x_0 \in \ZZ_+$, \ $(\ell_k)_{k=1}^n \in \ZZ_+^n$,
 \ $(x_{k,j})_{k,j=1}^n \in \{0, 1\}^{n\times n}$, \ and
 \begin{align*}
  \sum\Bigl\{p_n\left(a, x_0, (\ell_k)_{k=1}^n, (x_{k,j})_{k,j=1}^n\right)
             : x_0 \in \ZZ_+, \, (\ell_k)_{k=1}^n \in \ZZ_+^n,
                \, (x_{k,j})_{k,j=1}^n \in \{0, 1\}^{n\times n}\Bigr\}
  = 1
 \end{align*}
 for all \ $a \in (0, 1)$.
\ Then the Kolmogorov consistency theorem implies the existence of a probability
 space \ $(\Omega, \cA, \PP)$ \ and random variables \ $\alpha$, \ $X_0$,
 \ $\vare_k$ \ and \ $\xi_{k,j}$ \ for \ $k, j \in \NN$ \ on \ $(\Omega, \cA, \PP)$
 \ with the desired properties \eqref{help_stac_RINAR1_1},
 \eqref{help_stac_RINAR1_2} and \eqref{help_stac_RINAR1_3}, since for all
 \ $n \in \NN$, \ we have
 \begin{align*}
  &\sum\bigl\{p_{n+1}\bigl(a, x_0, (\ell_k)_{k=1}^{n+1},
                           (x_{k,j})_{k,j=1}^{n+1}\bigr) \\
  &\phantom{\sum\bigl\{}
              : \ell_{n+1} \in \ZZ_+ ,
                \, (x_{n+1,j})_{j=1}^n , (x_{k,n+1})_{k=1}^n \in \{0, 1\}^n ,
                \, x_{n+1,n+1} \in \{0, 1\}\bigr\} \\
  &= p_n\bigl(a, x_0, (\ell_k)_{k=1}^n, (x_{k,j})_{k,j=1}^n\bigr) .
 \end{align*}
Define a process \ $(X_k)_{k\in\ZZ_+}$ \ by
 \[
  X_k = \sum_{j=1}^{X_{k-1}} \xi_{k,j} + \vare_k , \qquad k \in \NN .
 \]
By Section \ref{prel}, conditionally on \ $\alpha$, \ the process
 \ $(X_k)_{k\in\ZZ_+}$ \ is a strictly stationary INAR(1) process
 with thinning parameter \ $\alpha$ \ and with Poisson innovations.
Moreover, by the law of total probability, it is also (unconditionally) strictly
 stationary but it is not a Markov chain (so it is not an INAR(1) process) if
 \ $\alpha$ \ is not degenerate, see Appendix \ref{Markov}.
The process \ $(X_k)_{k\in\ZZ_+}$ \ can be called a randomized INAR(1)
 process with Poisson innovations, and the distribution of \ $\alpha$ \ is the
 so-called mixing distribution of the model.
The conditional generator function of \ $X_0$ \ given \ $\alpha \in (0, 1)$ \ has
 the form
 \[
  F_0( z_0\mid \alpha)
  := \EE( z_0^{X_0} \mid \alpha)
  = \ee^{(1-\alpha)^{-1}\lambda({z_0}-1)} ,
  \qquad z_0 \in \CC ,
 \]
 and the conditional expectation of \ $X_0$ \ given \ $\alpha$ \ is
 \ $\EE(X_0 \mid \alpha) = (1-\alpha)^{-1}\lambda$.
\ Here and in the sequel conditional expectations like
 \ $\EE( z_0^{X_0} \mid \alpha)$ \ or \ $\EE(X_0 \mid \alpha)$ \ are meant in the
 generalized sense, see, e.g., in Stroock \cite[\S\,5.1.1]{Str}.
The joint conditional generator function of \ $X_0, X_1, \ldots, X_k$ \ given
 \ $\alpha$ \ will be denoted by \ $F_{0,\ldots,k}(z_0, \ldots, z_k\mid\alpha)$,
 \ $z_0, \ldots, z_k \in \CC$.

Let \ $\alpha^{(j)}$, \ $j \in \NN$, \ be a sequence of independent copies of the
 random variable \ $\alpha$, \ and let \ $(X^{(j)}_k)_{k\in\ZZ_+}$, \ $j \in \NN$,
 \ be a sequence of independent copies of the process \ $(X_k)_{k\in\ZZ_+}$
 \ with idiosyncratic innovations (i.e., the innovations
 \ $(\vare^{(j)}_k)_{k\in\ZZ_+}$, $j\in\NN$, \ belonging to
 \ $(X^{(j)}_k)_{k\in\ZZ_+}$, \ $j \in \NN$, \ are independent) such that
 \ $(X^{(j)}_k)_{k\in\ZZ_+}$ \ conditionally on \ $\alpha^{(j)}$ \ is a strictly
 stationary INAR(1) process with thinning parameter \ $\alpha^{(j)}$ \ and
 with Poisson innovations for all \ $j \in \NN$.

First we consider a simple aggregation procedure.
For each \ $N \in \NN$, \ consider the stochastic process
 \ $\tS^{(N)} = (\tS^{(N)}_k)_{k\in\ZZ_+}$ \ given by
 \[
   \tS^{(N)}_k
   := \sum_{j=1}^N (X^{(j)}_k - \EE(X^{(j)}_k \mid \alpha^{(j)}))
   = \sum_{j=1}^N \Bigl(X^{(j)}_k - \frac{\lambda}{1-\alpha^{(j)}}\Bigr) ,
   \qquad k \in \ZZ_+ .
 \]

\begin{Pro}\label{simple_aggregation_random}
If \ $\EE\bigl(\frac{1}{1-\alpha}\bigr) < \infty$, \ then
 \[
   N^{-\frac{1}{2}} \tS^{(N)}
   \distrf \tcY \qquad \text{as \ $N \to \infty$,}
 \]
 where \ $(\tcY_k)_{k\in\ZZ_+}$ \ is a stationary Gaussian process with zero mean
 and covariances
 \begin{equation}\label{covariance}
  \EE(\tcY_0 \tcY_k)
  = \cov\left(X_0 - \frac{\lambda}{1-\alpha},
              X_k - \frac{\lambda}{1-\alpha}\right)
  = \lambda \EE\Bigl(\frac{\alpha^k}{1-\alpha}\Bigr) , \qquad k \in \ZZ_+ .
 \end{equation}
\end{Pro}


\begin{Pro}\label{simple_aggregation_random3}
We have
 \[
   \biggl(n^{-\frac{1}{2}}
          \sum_{k=1}^{\lfloor nt \rfloor}
           \tS^{(1)}_k\biggr)_{t\in\RR_+}
   = \biggl(n^{-\frac{1}{2}}
            \sum_{k=1}^{\lfloor nt \rfloor}
            (X^{(1)}_k - \EE(X^{(1)}_k \mid \alpha^{(1)}))\biggr)_{t\in\RR_+}
   \distrf \frac{\sqrt{\lambda(1+\alpha)}}{1-\alpha} B
 \]
 as \ $n \to \infty$, \ where \ $B = (B_t)_{t\in\RR_+}$ \ is a standard Brownian
 motion, independent of \ $\alpha$.
\end{Pro}

In the next two propositions, which are counterparts of Propositions \ref{simple_aggregation}
 and \ref{simple_aggregation2}, we point out that the usual centralization leads to
 limit theorems similar to Propositions \ref{simple_aggregation_random} and
 \ref{simple_aggregation_random3}, but with an occasionally
  different scaling and with a different limit process.
We use again the notation \ $S^{(N)} = (S^{(N)}_k)_{k\in\ZZ_+}$ \ given in \eqref{SN}
 for the simple aggregation (with the usual centralization) of the randomized process.


\begin{Pro}\label{simple_aggregation_random2}
If \ $\EE\bigl(\frac{1}{(1-\alpha)^2}\bigr) < \infty$, \ then
 \[
   N^{-\frac{1}{2}} S^{(N)}
   \distrf \cY
   \qquad \text{as \ $N \to \infty$,}
 \]
 where \ $\cY = (\cY_k)_{k\in\ZZ_+}$ \ is a stationary Gaussian process with zero mean
 and covariances
 \[
   \EE(\cY_0 \cY_k)
   = \cov(X_0 , X_k)
   = \lambda \EE\Bigl(\frac{\alpha^k}{1-\alpha}\Bigr)
     + \lambda^2 \var \left(\frac{1}{1-\alpha}\right) , \qquad k \in \ZZ_+.
 \]
\end{Pro}


\begin{Pro}\label{simple_aggregation_random4}
If \ $\EE\bigl(\frac{1}{1-\alpha}\bigr) < \infty$, \ then
 \begin{align*}
   \Biggl(n^{-1}
          \sum_{k=1}^{\lfloor nt \rfloor}
           S^{(1)}_k\Biggr)_{t\in\RR_+}
   = \Biggl(n^{-1}
            \sum_{k=1}^{\lfloor nt \rfloor}
            (X^{(1)}_k - \EE(X^{(1)}_k))\Biggr)_{t\in\RR_+}
   \distrf
   \biggl(\biggl(\frac{\lambda}{1-\alpha}
                 - \EE\biggl(\frac{\lambda}{1-\alpha}\biggr)\biggr)t\biggr)_{t\in\RR_+}
 \end{align*}
 as \ $n \to \infty$.
\end{Pro}

In Proposition \ref{simple_aggregation_random4} the limit process is simply a line with
 a random slope.

In the forthcoming Theorems
 \ref{joint_aggregation_random}--\ref{joint_aggregation_random_5.5}, we assume that the
 distribution of the random variable \ $\alpha$, \ i.e., the mixing distribution, has a
 probability density of the form
 \begin{equation}\label{alpha}
  \psi(x) (1 - x)^\beta , \qquad x \in (0, 1) ,
 \end{equation}
 where \ $\psi$ \ is a function on \ $(0, 1)$ \ having a limit
 \ $\lim_{x\uparrow 1} \psi(x) = \psi_1 \in (0, \infty)$.
\ Note that necessarily \ $\beta \in (-1, \infty)$ \ (otherwise
 \ $\int_0^1 \psi(x) (1 - x)^\beta \, \dd x = \infty$), \
 the function \ $(0,1)\ni x\mapsto \psi(x)$ \ is integrable on \ $(0,1)$,
 \ and the function \ $(0, 1) \ni x \mapsto \psi(x) (1 - x)^\beta$ \ is regularly varying at the point 1
 (i.e., \ $(0,\infty)\ni x\mapsto \psi(1-\frac{1}{x}) \,x^{-\beta}$ \ is regularly varying at infinity).
Further, in case of
 \ $\psi(x) = \frac{\Gamma(a+\beta+2)}{\Gamma(a+1)\Gamma(\beta+1)} x^a$,
 \ $x \in (0, 1)$,  \ with some \ $a \in (-1, \infty)$, \ the random variable
 \ $\alpha$ \ is Beta distributed with parameters \ $a + 1$ \ and \ $\beta + 1$.
\ The special case of Beta mixing distribution is an important one from the
 historical point of view, since the Nobel prize winner Clive W. J. Granger used
 Beta distribution as a mixing distribution for random coefficient AR(1) processes,
 see Granger \cite{Gra}.

\begin{Rem}\label{condexp}
Under the condition \eqref{alpha}, for each \ $\ell \in \NN$, \ the expectation
 \ $\EE\bigl(\frac{1}{(1-\alpha)^\ell}\bigr)$ \ is finite if and only if
 \ $\beta > \ell - 1$.
\ Indeed, if \ $\beta > \ell - 1$, \ then, by choosing \ $\vare \in (0, 1)$ \ with
 \ $\sup_{a\in(1-\vare,1)} \psi(a) \leq 2 \psi_1$, \ we have
 \ $\EE\bigl(\frac{1}{(1-\alpha)^\ell}\bigr) = I_1(\vare) + I_2(\vare)$, \ where
 \begin{align*}
  I_1(\vare)
  &:= \int_0^{1-\vare} \psi(a) (1-a)^{\beta-\ell} \, \dd a
  \leq \vare^{\beta-\ell} \int_0^{1-\vare} \psi(a) \, \dd a < \infty , \\
  I_2(\vare)
  &:= \int_{1-\vare}^1 \psi(a) (1-a)^{\beta-\ell} \, \dd a
  \leq 2 \psi_1 \int_{1-\vare}^1 (1-a)^{\beta-\ell} \, \dd a
  = \frac{2\psi_1\vare^{\beta-\ell+1}}{\beta-\ell+1} < \infty .
 \end{align*}
Conversely, if \ $\beta \leq \ell - 1$, \ then, by choosing \ $\vare \in (0, 1)$
 \ with \ $\sup_{a\in(1-\vare,1)} \psi(a) \geq \psi_1/2$, \ we have
 \begin{align*}
  \EE\left(\frac{1}{(1-\alpha)^\ell}\right)
  \geq \int_{1-\vare}^1 \psi(a) (1-a)^{\beta-\ell} \, \dd a
  \geq \frac{\psi_1}{2} \int_{1-\vare}^1 (1-a)^{\beta-\ell} \, \dd a
  = \infty .
 \end{align*}
This means that in case of \ $\beta \in (-1, 0]$, \ the processes
 \ $S^{(N,n)} = (S_t^{(N,n)})_{t\in\RR_+}$, \ $N, n \in \NN$, \ given in
 \eqref{SNnt} are not defined for the randomized INAR(1) process introduced in
 this section with mixing distribution given in \eqref{alpha}.
Moreover, the Propositions \ref{simple_aggregation_random},
 \ref{simple_aggregation_random3}, \ref{simple_aggregation_random2} and
 \ref{simple_aggregation_random4} are valid in case of \ $\beta > 0$,
 \ $\beta > -1$, \ $\beta > 1$ \ and \ $\beta > 0$, \ respectively.
\proofend
\end{Rem}

For each \ $N, n \in \NN$, \ consider the stochastic process
 \ $\widetilde S^{(N,n)} = (\widetilde S_t^{(N,n)})_{t\in\RR_+}$ \ given by
 \[
   \widetilde S_t^{(N,n)}
   := \sum_{j=1}^N \sum_{k=1}^{\lfloor nt \rfloor}
       (X^{(j)}_k - \EE(X^{(j)}_k \mid \alpha^{(j)})),
   \qquad t\in\RR_+.
 \]

\begin{Rem}
If \ $\beta>0$, \ then the covariances of the strictly stationary process
 \ $(X_k - \EE(X_k\mid \alpha))_{k\in\ZZ_+} = (X_k - \frac{\lambda}{1-\alpha})_{k\in\ZZ_+}$ \
 exist and take the form
 \[
   \cov\big(X_0 - \EE(X_0\mid \alpha), X_k - \EE(X_k\mid \alpha)\big)
       = \EE\left(\frac{\lambda \alpha^k}{1-\alpha}\right),
       \qquad k\in\ZZ_+,
 \]
 see \eqref{help15}.
Further,
 \begin{align*}
  \sum_{k=0}^\infty \Big\vert    \cov(X_0 - \EE(X_0\mid \alpha), X_k - \EE(X_k\mid \alpha)) \Big\vert
    & = \sum_{k=0}^\infty \EE\left(\frac{\lambda \alpha^k}{1-\alpha}\right)
      = \lambda \EE\left(\frac{1}{1-\alpha} \sum_{k=0}^\infty \alpha^k \right)\\
    & = \lambda \EE\left(\frac{1}{(1-\alpha)^2} \right),
 \end{align*}
 which is finite if and only if \ $\beta>1$, \ see Remark \ref{condexp}.
This means that the strictly stationary process
 \ $(X_k - \EE(X_k\mid \alpha))_{k\in\ZZ_+}$ \ has short memory (i.e., it has summable covariances) if \ $\beta>1$, \, and long memory if \ $\beta\in(0,1]$ \ (i.e., it has non-summable covariances).
\proofend
\end{Rem}


For \ $\beta \in (0, 2)$, \ let
 \ $(\cB_{1-\frac{\beta}{2}}(t))_{t\in\RR_+}$ \ denote a fractional Brownian
 motion with parameter \ $1 - \beta/2$, \ that is a Gaussian process with zero mean
 and covariance function
 \begin{align}\label{fractional_BM}
  \cov(\cB_{1-\frac{\beta}{2}}(t_1),\cB_{1-\frac{\beta}{2}}(t_2))
  =\frac{t_1^{2-\beta}+t_2^{2-\beta}-|t_2-t_1|^{2-\beta}}{2},
  \qquad t_1, t_2 \in \RR_+ .
 \end{align}

In Appendix \ref{App2} we recall an integral representation of the fractional
 Brownian motion \ $(\cB_{1-\frac{\beta}{2}}(t))_{t\in\RR_+}$ \ due to
 Pilipauskait\.{e} and Surgailis \cite{PilSur} in order to connect our forthcoming
 results with the ones in Pilipauskait\.{e} and Surgailis \cite{PilSur} and in
 Puplinskait\.{e} and Surgailis \cite{PupSur1}, \cite{PupSur2}.

The next three results are limit theorems for appropriately scaled versions of
 \ $\widetilde S^{(N,n)}$, \ first taking the limit \ $N \to \infty$ \ and then
 \ $n \to \infty$ \ in the case \ $\beta \in (-1,1)$, \
 which are counterparts of (2.7), (2.8) and (2.9) of Theorem 2.1 in Pilipauskait\.e and Surgailis
 \cite{PilSur}, respectively.

\begin{Thm}\label{joint_aggregation_random}
If \ $\beta \in (0, 1)$, \ then
 \[
   \cD_\ff\text{-}\hspace*{-1mm}\lim_{n\to\infty} \,
   \cD_\ff\text{-}\hspace*{-1mm}\lim_{N\to\infty} \,
    n^{-1+\frac{\beta}{2}} N^{-\frac{1}{2}} \,
   \widetilde S^{(N,n)}
   = \sqrt{\frac{2\lambda\psi_1\Gamma(\beta) }{(2-\beta)(1-\beta)}}
     \cB_{1-\frac{\beta}{2}} .
 \]
\end{Thm}

\begin{Thm}\label{joint_aggregation_random_2}
If \ $\beta \in (-1, 0)$, \ then
 \[
   \cD_\ff\text{-}\hspace*{-1mm}\lim_{n\to\infty} \,
   \cD_\ff\text{-}\hspace*{-1mm}\lim_{N\to\infty} \,
   n^{-1} N^{-\frac{1}{2(1+\beta)}} \,
   \tS^{(N,n)}
   = (V_{2(1+\beta)} t)_{t\in\RR_+} ,
 \]
 where \ $V_{2(1+\beta)}$ \ is a symmetric \ $2(1+\beta)$-stable random variable
 (not depending on \ $t$) \ with characteristic function
 \[
   \EE(\ee^{\ii \theta V_{2(1+\beta)}})
   =\ee^{-K_\beta|\theta|^{2(1+\beta)}}, \qquad \theta \in \RR ,
 \]
 where
 \[
   K_\beta
   := \psi_1\left(\frac{\lambda}{2}\right)^{1+\beta}
      \frac{\Gamma(-\beta)}{1+\beta} .
 \]
\end{Thm}

\begin{Thm}\label{joint_aggregation_random_3}
If \ $\beta = 0$, \ then
 \[
   \cD_\ff\text{-}\hspace*{-1mm}\lim_{n\to\infty} \,
   \cD_\ff\text{-}\hspace*{-1mm}\lim_{N\to\infty} \,
    n^{-1} (N\log N)^{-\frac{1}{2}} \,
   \tS^{(N,n)}
   = (W_{\lambda\psi_1} t)_{t\in\RR_+} ,
 \]
 where \ $W_{\lambda\psi_1}$ \ is a normally distributed random variable with mean
 zero and with variance \ $\lambda\psi_1$.
\end{Thm}

The next result is a limit theorem for an appropriately scaled version of
 \ $\tS^{(N,n)}$, \ first taking the limit \ $n \to \infty$ \ and then
 \ $N \to \infty$ \ in the case \ $\beta \in (-1, 1)$,
 \ which is a counterpart of (2.10) of Theorem 2.1 in Pilipauskait\.e and
 Surgailis \cite{PilSur}.

\begin{Thm}\label{joint_aggregation_random_4}
If \ $\beta \in (-1, 1)$, \ then
 \[
   \cD_\ff\text{-}\hspace*{-1mm}\lim_{N\to\infty} \,
   \cD_\ff\text{-}\hspace*{-1mm}\lim_{n\to\infty} \,
   N^{-\frac{1}{1+\beta}}n^{-\frac{1}{2}} \,
   \tS^{(N,n)}
   = \cY_{1+\beta} ,
 \]
 where
 \ $\cY_{1+\beta} = \big(\cY_{1+\beta}(t)
    := \sqrt{Y_{(1+\beta)/2}} \, B_t\big)_{t \in \RR_+}$
 \ is a \ $(1+\beta)$-stable L\'evy process.
Here \ $Y_{(1+\beta)/2}$ \ is a positive \ $\frac{1+\beta}{2}$-stable random
 variable with Laplace transform
 \ $\EE(\ee^{-\theta Y_{(1+\beta)/2}})
    = \ee^{-k_\beta \theta^{\frac{1+\beta}{2}}}$, \ $\theta \in \RR_+$,
 \ and with characteristic function
 \[
  \EE(\ee^{\ii\theta Y_{(1+\beta)/2}})
      = \exp\left\{ -k_\beta \vert\theta\vert^{\frac{1+\beta}{2}}
                     \ee^{-\ii \sign(\theta) \frac{\pi(1+\beta)}{4}}\right\},
    \qquad \theta\in\RR,
 \]
 where
 \[
   k_\beta
   := \frac{(2\lambda)^{\frac{1+\beta}{2}} \psi_1}{1+\beta}
      \Gamma\left(\frac{1-\beta}{2}\right) ,
 \]
 and \ $(B_t)_{t\in\RR_+}$ \ is an independent standard Wiener process.
\end{Thm}

Next we show an iterated scaling limit theorem where the order of
 the iteration can be arbitrary in the case \ $\beta \in (1, \infty)$,
 \ which is a counterpart of Theorem 2.3 in Pilipauskait\.e and Surgailis
 \cite{PilSur}.

\begin{Thm}\label{joint_aggregation_random_5}
If \ $\beta \in (1, \infty)$, \ then
 \[
   \cD_\ff\text{-}\hspace*{-1mm}\lim_{n\to\infty} \,
   \cD_\ff\text{-}\hspace*{-1mm}\lim_{N\to\infty} \,
   (nN)^{-\frac{1}{2}} \tS^{(N,n)}
   = \cD_\ff\text{-}\hspace*{-1mm}\lim_{N\to\infty} \,
     \cD_\ff\text{-}\hspace*{-1mm}\lim_{n\to\infty} \,
     (nN)^{-\frac{1}{2}} \tS^{(N,n)}
   = \sigma B,
 \]
 where \ $\sigma^2 := \lambda \EE((1+\alpha)(1-\alpha)^{-2})$
 \ and \ $(B_t)_{t\in\RR_+}$ \ is a standard Wiener process.
\end{Thm}

By Remark \ref{condexp}, if \ $\beta > 1$, \ then
 \ $\EE\bigl(\frac{1}{(1-\alpha)^2}\bigr) < \infty$, \ and hence
 \ $\sigma^2 < \infty$, \ where \ $\sigma^2$ \ is given in
 Theorem \ref{joint_aggregation_random_5}.

In the next theorems we consider the usual centralization with \ $\EE(X^{(j)}_k)$
 \ in the cases \ $\beta \in (0, 1)$ \ and \ $\beta > 1$.
\ These are the counterparts of Theorems \ref{joint_aggregation_random},
 \ref{joint_aggregation_random_4} and \ref{joint_aggregation_random_5}.
Recall that, due to Remark \ref{condexp}, the expectation
 \ $\EE(X_0) = \EE\bigl(\frac{\lambda}{1-\alpha}\bigr)$ \ is finite if and only if
 \ $\beta > 0$, \ so Theorems \ref{joint_aggregation_random_2} and
 \ref{joint_aggregation_random_3} can not have counterparts in this sense.

\begin{Thm}\label{joint_aggregation_random_4.5}
If \ $\beta \in (0, 1)$, \ then
 \begin{align*}
  \cD_\ff\text{-}\hspace*{-1mm}\lim_{n\to\infty} \,
   \cD_\ff\text{-}\hspace*{-1mm}\lim_{N\to\infty} \,
   n^{-1} N^{-\frac{1}{1+\beta}}
   S^{(N,n)}
  = \cD_\ff\text{-}\hspace*{-1mm}\lim_{N\to\infty} \,
   \cD_\ff\text{-}\hspace*{-1mm}\lim_{n\to\infty} \,
   n^{-1} N^{-\frac{1}{1+\beta}}
   S^{(N,n)}
   = \bigl(Z_{1+\beta} \,t\bigr)_{t\in\RR_+} ,
 \end{align*}
 where \ $Z_{1+\beta}$ \ is a \ $(1+\beta)$-stable random variable with
 characteristic function
 \ $\EE(\ee^{\ii \theta Z_{1+\beta}}) = \ee^{-|\theta|^{1+\beta}\omega_\beta(\theta)}$,
 \ $\theta\in\RR$, \ where
 \[
   \omega_\beta(\theta)
   := \frac{\psi_1\Gamma(1-\beta)\lambda^{1+\beta}}
           {-\beta(1+\beta)} \ee^{-\ii\pi\sign(\theta)(1+\beta)/2} ,
   \qquad \theta \in \RR .
 \]
\end{Thm}

\begin{Thm}\label{joint_aggregation_random_5.5}
If \ $\beta \in (1, \infty)$, \ then
 \begin{align*}
  \cD_\ff\text{-}\hspace*{-1mm}\lim_{n\to\infty} \,
   \cD_\ff\text{-}\hspace*{-1mm}\lim_{N\to\infty} \,
   n^{-1} N^{-\frac{1}{2}} S^{(N,n)}
  = \cD_\ff\text{-}\hspace*{-1mm}\lim_{N\to\infty} \,
     \cD_\ff\text{-}\hspace*{-1mm}\lim_{n\to\infty} \,
     n^{-1} N^{-\frac{1}{2}} S^{(N,n)}
   = (W_{\lambda^2\var((1-\alpha)^{-1})} \,t)_{t\in\RR_+} ,
 \end{align*}
 where \ $W_{\lambda^2\var((1-\alpha)^{-1})}$ \ is a normally distributed
 random variable with mean zero and with variance
 \ $\lambda^2\var((1-\alpha)^{-1})$.
\end{Thm}

In case of Theorems \ref{joint_aggregation_random_2},
 \ref{joint_aggregation_random_3}, \ref{joint_aggregation_random_4.5} and
 \ref{joint_aggregation_random_5.5} the limit processes are lines with random
 slopes.

We point out that the processes of doubly indexed partial sums, \ $S^{(N,n)}$ \ and
 \ $\tS^{(N,n)}$ \ contain the expected or conditional expected values of the processes \ $X^{(j)}$, \ $j\in\NN$.
\ Therefore, in a statistical testing, they could not be used directly.
So we consider a similar process
 \[
 \hS^{(N,n)}_t:=\sum_{j=1}^N \sum_{k=1}^{\lfloor nt\rfloor} \left[X_k^{(j)}-\frac{\sum_{\ell=1}^nX_\ell^{(j)}}{n}\right], \qquad t\in\RR_+,
 \]
 which does not require the knowledge of the expectation or conditional expectation of the processes  \ $X^{(j)}$, \ $j\in\NN$.
Note that the summands in \ $\hS^{(N,n)}_t$ \ have \ 0 \ conditional means with respect to \ $\alpha$, \
 so we do not need any additional centering.
Moreover, \ $\hS^{(N,n)}$ \ is related to the two previously examined processes in the following way: in case of \ $\beta\in(0,\infty)$ \ (which ensures the existence of \ $\EE(X_k^{(j)})$, \ $k\in\ZZ_+$), \ we have
 \begin{equation*}
  \hS^{(N,n)}_t
  =
 \sum_{j=1}^N \sum_{k=1}^{\lfloor nt\rfloor}\left[
 X^{(j)}_k-\EE(X^{(j)}_k)
 -
 \frac{\sum_{\ell=1}^n(X^{(j)}_\ell-\EE(X^{(j)}_\ell))}{n}
 \right]=S^{(N,n)}_t-\frac{\lfloor nt\rfloor}{n}S^{(N,n)}_1,
 \end{equation*}
 and in case of \ $\beta\in(-1,\infty)$,
 \begin{equation*}
 \hS^{(N,n)}_t
 =
 \sum_{j=1}^N \sum_{k=1}^{\lfloor nt\rfloor}\left[X^{(j)}_k-\EE(X^{(j)}_k\mid \alpha^{(j)})
 -
 \frac{\sum_{\ell=1}^n(X^{(j)}_\ell-\EE(X^{(j)}_\ell\mid \alpha^{(j)}))}{n}
 \right]=\tS^{(N,n)}_t-\frac{\lfloor nt\rfloor}{n}\tS^{(N,n)}_1
 \end{equation*}
 for every \ $t\in\RR_+$.
\ Therefore, by Theorem \ref{joint_aggregation_random}, Theorem \ref{joint_aggregation_random_4},
 and Theorem \ref{joint_aggregation_random_5}, using Slutsky's lemma, the following limit theorems hold.

\begin{Cor}\label{Cor1}
If \ $\beta \in (0, 1)$, \ then
 \[
   \cD_\ff\text{-}\hspace*{-1mm}\lim_{n\to\infty} \,
   \cD_\ff\text{-}\hspace*{-1mm}\lim_{N\to\infty} \,
    n^{-1+\frac{\beta}{2}} N^{-\frac{1}{2}} \,
   \hS^{(N,n)}
   = \sqrt{\frac{2\lambda\psi_1\Gamma(\beta) }{(2-\beta)(1-\beta)}}
     \left(\cB_{1-\frac{\beta}{2}}(t)-t \cB_{1-\frac{\beta}{2}}(1)\right)_{t\in\RR_+} ,
 \]
where the process \ $\cB_{1-\frac{\beta}{2}}$ \ is given by \eqref{fractional_BM}.

If \ $\beta \in (-1, 1)$, \ then
 \[
   \cD_\ff\text{-}\hspace*{-1mm}\lim_{N\to\infty} \,
   \cD_\ff\text{-}\hspace*{-1mm}\lim_{n\to\infty} \,
   N^{-\frac{1}{1+\beta}}n^{-\frac{1}{2}} \,
   \hS^{(N,n)}
   = \left(\cY_{1+\beta}(t)-t \cY_{1+\beta}(1)\right)_{t\in \RR_+} ,
 \]
 where the process \ $\cY_{1+\beta}$ \ is given in Theorem \ref{joint_aggregation_random_4}.

If \ $\beta \in (1, \infty)$, \ then
 \[
   \cD_\ff\text{-}\hspace*{-1mm}\lim_{n\to\infty} \,
   \cD_\ff\text{-}\hspace*{-1mm}\lim_{N\to\infty} \,
   (nN)^{-\frac{1}{2}} \hS^{(N,n)}
   = \cD_\ff\text{-}\hspace*{-1mm}\lim_{N\to\infty} \,
     \cD_\ff\text{-}\hspace*{-1mm}\lim_{n\to\infty} \,
     (nN)^{-\frac{1}{2}} \hS^{(N,n)}
   = \sigma (B_t-tB_1)_{t\in\RR_+},
 \]
 where \ $\sigma^2$ \ and the process \ $B$ \ are given in Theorem \ref{joint_aggregation_random_5}.
\end{Cor}

In Corollary \ref{Cor1}, the limit processes restricted on the time interval \ $[0,1]$ \ are bridges
 in the sense that they take the same value (namely, 0) at the time points 0 and 1, and especially, in case of
 \ $\beta\in(1,\infty)$, \ it is a Wiener bridge.
We note that no counterparts appear for the rest of the theorems because in those cases the limit processes are
 lines with random slopes, which result the constant zero process in this alternative case.
In case of \ $\beta\in(-1,0]$, \ by applying some smaller scaling factors,
 one could try to achieve a non-degenerate weak limit of \ $\hS^{(N,n)}$ \ by first taking
 the limit \ $N\to\infty$ \ and then that of \ $n\to\infty$.

\section{Proofs}
\label{Proofs}

Theorem \ref{joint_aggregation_random} is a counterpart of (2.7) of Theorem 2.1 in
  Pilipauskait\.e and Surgailis \cite{PilSur}.
We will present two proofs of Theorem \ref{joint_aggregation_random}, and we call
 the attention that both proofs are completely different from the proof of (2.7) in
 Theorem 2.1 in Pilipauskait\.e and Surgailis \cite{PilSur} (suspecting also that
 their result in question might be proved by our method as well).
Theorems \ref{joint_aggregation_random_2} and \ref{joint_aggregation_random_3} are counterparts
 of (2.8) and (2.9) of Theorem 2.1 in Pilipauskait\.e and Surgailis \cite{PilSur}.
The proofs of these theorems use the same technique, namely, expansions of characteristic functions,
 and we provide all the technical details.
Theorem \ref{joint_aggregation_random_4} is a counterpart of (2.10) of Theorem 2.1 in
 Pilipauskait\.e and Surgailis \cite{PilSur}.
We give two proofs of Theorem \ref{joint_aggregation_random_4}: the first one is based on
 expansions of characteristic functions (as the proof of (2.10) of Theorem 2.1 in Pilipauskait\.e and Surgailis \cite{PilSur}),
 the second one reduces to show that \ $\frac{\lambda(1+\alpha)}{(1-\alpha)^2}$ \ belongs to the domain of normal attraction
 of the \ $\frac{1+\beta}{2}$-stable law of \ $Y_{\frac{1+\beta}{2}}$.
Theorem \ref{joint_aggregation_random_5} is a counterpart of Theorem 2.3 in Pilipauskait\.e and Surgailis \cite{PilSur}.
The proof of Theorem \ref{joint_aggregation_random_5} is based on the multidimensional central limit theorem and checking
 convergence of covariances of some Gaussian processes.

The notations \ $\OO(1)$ \ and \ $|\OO(1)|$ \ stand for a possibly complex and
 respectively real sequence \ $(a_k)_{k\in\NN}$ \ that is bounded and can only depend
 on the parameters \ $\lambda$, \ $\psi_1$, \ $\beta$, \ and on some fixed \ $m \in \NN$
 \ and \ $\theta_1,\dots, \theta_m \in \RR$.
\ Further, we call the attention that several \ $\OO(1)$-s (respectively
 \ $|\OO(1)|$-s) \ in the same formula do not necessarily mean the same bounded
 sequence.

\noindent{\bf Proof of Proposition \ref{Pro_egyuttes_mom_gen}.}
First we prove \eqref{help1_alter} by induction.
Note that by \eqref{square0} the statement holds for \ $k = 0$.
\ We suppose that it holds for \ $0,\ldots,k$, \ and show that it is also true for
 \ $k+1$.
\ Using \eqref{square-1} it is easy to see that
 \begin{equation*}
  \begin{split}
  &F_{0,\dots,k,k+1}(z_0, \dots, z_k, z_{k+1})
   =\EE\left(z_0^{X_0} \cdots z_k^{X_k} z_{k+1}^{X_{k+1}}\right) \\
  &=\EE\left(z_0^{X_0} \cdots z_k^{X_k}
             \EE\left(z_{k+1}^{X_{k+1}} \mid X_0, \ldots, X_k\right)\right)
   =\EE\left(z_0^{X_0} \cdots z_k^{X_k}
             \EE\left(z_{k+1}^{X_{k+1}} \mid X_k\right)\right) \\
  &=\EE\left(z_0^{X_0} \cdots z_k^{X_k} \ee^{\lambda(z_{k+1}-1)}
             (1 - \alpha + \alpha z_{k+1})^{X_k}\right) .
  \end{split}
 \end{equation*}
On the one hand, for any \ $z_0, \dots, z_{k+1} \in \CC$, \ by the assumption of
 the induction,
 \begin{equation*}
  \begin{split}
   &F_{0,\dots,k,k+1}(z_0, \dots, z_k, z_{k+1})
    = \ee^{\lambda(z_{k+1}-1)}
      F_{0,\dots,k}(z_0, \dots, z_{k-1}, z_k (1 - \alpha + \alpha z_{k+1})) \\
   &=\exp\biggl\{\frac{\lambda}{1-\alpha}
                 \biggl[(1-\alpha)(z_{k+1}-1)
                        + \sum_{0\leq i\leq j \leq k-1}
                           \alpha^{j-i} (z_i - 1) z_{i+1} \cdots z_{j-1} (z_j - 1) \\
   &\phantom{=\exp\biggl\{\frac{\lambda}{1-\alpha}\biggl[}
                        +\mathrm{Sum}_1
                        + z_k (1 - \alpha + \alpha z_{k+1}) - 1 \biggr]\biggr\} ,
  \end{split}
 \end{equation*}
 with
 \[
   \mathrm{Sum}_1
   := \sum_{0\leq i\leq k-1}
       \alpha^{k-i} (z_i - 1) z_{i+1} \cdots
       z_{k-1} [z_k (1 - \alpha + \alpha z_{k+1}) - 1] .
 \]
On the other hand, the right hand side of \eqref{help1_alter} for \ $k+1$ \ has the
 form
 \begin{equation*}
  \begin{split}
  &\exp\Bigg\{\frac{\lambda}{1-\alpha}
              \Bigg[\sum_{0\leq i\leq j \leq k-1}
                     \alpha^{j-i} (z_i - 1) z_{i+1} \cdots z_{j-1} (z_j - 1)
                    + \mathrm{Sum}_2 + \mathrm{Sum}_3\Bigg]\Bigg\} ,
  \end{split}
 \end{equation*}
 where
 \begin{equation*}
  \begin{split}
   \mathrm{Sum}_2
   &=\sum_{0\leq i\leq k} \alpha^{k-i} (z_i - 1) z_{i+1} \cdots z_{k-1} (z_k - 1) \\
   &=(z_k - 1)
     + \sum_{0\leq i\leq k-1}
        \alpha^{k-i} (z_i - 1) z_{i+1} \cdots z_{k-1} (z_k - 1) ,
  \end{split}
 \end{equation*}
 and
 \begin{equation*}
  \begin{split}
   &\mathrm{Sum}_3
    = \sum_{0\leq i\leq k+1}
       \alpha^{k+1-i} (z_i - 1) z_{i+1} \cdots z_k (z_{k+1} - 1)\\
   &= (z_{k+1} - 1) + \alpha(z_k-1) (z_{k+1} - 1)
      + \sum_{0\leq i\leq k-1}
         \alpha^{k+1-i} (z_i - 1) z_{i+1} \cdots z_k (z_{k+1} - 1) .
  \end{split}
 \end{equation*}
Since
 \begin{equation*}
   \mathrm{Sum}_1
   = \sum_{0\leq i\leq k-1}
       \alpha^{k-i} (z_i - 1) z_{i+1} \cdots z_{k-1} (z_k - 1)
      + \sum_{0\leq i\leq k-1}
         \alpha^{k+1-i} (z_i - 1) z_{i+1} \cdots z_k (z_{k+1} - 1) ,
 \end{equation*}
 in order to show \eqref{help1_alter} for \ $k + 1$, \ it is enough to check that
 \[
   (1-\alpha)(z_{k+1}-1) + z_k(1-\alpha+\alpha z_{k+1})-1
   = (z_k-1) + (z_{k+1}-1) + \alpha(z_k-1)(z_{k+1}-1) ,
 \]
 which holds trivially.

Now we prove \eqref{help1}.
In formula \eqref{help1_alter}, for fixed indices \ $0 \leq i \leq j \leq k$ \ the
 term in the sum gives
 \begin{align*}
  &(z_i - 1) z_{i+1} \cdots z_{j-1} (z_j - 1) \\
  &\quad
   = (z_i \cdots z_j-1) - (z_i \cdots z_{j-1} - 1) - (z_{i+1} \cdots z_j - 1)
     + (z_{i+1} \cdots z_{j-1} - 1) ,
 \end{align*}
 meaning that the sum consists of similar terms as in \eqref{help1}.
We only have to show that the coefficients coincide in the formulas \eqref{help1}
 and \eqref{help1_alter}.
In \eqref{help1} the coefficient of \ $z_i \cdots z_j-1$ \ is
 \ $\lambda (1-\alpha)^{K_{i,j,k}} \alpha^{j-i}$.
\ In \eqref{help1_alter} this term may appear multiple times, depending on the
 indices \ $i$ \ and \ $j$.
\ If \ $i = 0$ \ and \ $j = k$, \ then it only appears once, with coefficient
 \ $\lambda \alpha^{j-i}/(1-\alpha)$, \ that is the same as in \eqref{help1}.
However, if \ $i = 0$ \ and \ $0 \leq j \leq k-1$ \ in \eqref{help1}, then the term
 also appears when the indices are \ $i$ \ and \ $j+1$ \ in \eqref{help1_alter},
 meaning that the coefficient is
 \[
   \lambda
   \left(\frac{\alpha^{j-i}}{1-\alpha} - \frac{\alpha^{j+1-i}}{1-\alpha}\right)
   = \lambda\alpha^{j-i} ,
 \]
 which is the same as in \eqref{help1}.
Similarly, if \ $1 \leq i \leq k$ \ and \ $j = k$ \ in \eqref{help1}, then the term
 also appears when the indices are \ $i - 1$ \ and \ $j$ \ in \eqref{help1_alter},
 meaning that the coefficient is
 \[
   \lambda
   \left(\frac{\alpha^{j-i}}{1-\alpha}
         - \frac{\alpha^{j - (i-1)}}{1-\alpha}\right)
   = \lambda \alpha^{j-i} ,
 \]
 which is the same as in \eqref{help1}.
If \ $1 \leq i \leq j \leq k-1$ \ in \eqref{help1}, then the term appears three more
 times, for the index pairs \ $(i-1, j), \ (i, j+1), \ (i-1, j+1)$ \ in
 \eqref{help1_alter}, resulting the coefficient
 \[
   \lambda
   \left(\frac{\alpha^{j-i}}{1-\alpha}
         - \frac{\alpha^{j-(i-1)}}{1-\alpha}
         - \frac{\alpha^{(j+1)-i}}{1-\alpha}
         + \frac{\alpha^{(j+1)-(i-1)}}{1-\alpha}\right)
   = \lambda \alpha^{j-i} \frac{1-2\alpha+\alpha^2}{1-\alpha}
   = \lambda \alpha^{j-i} (1 - \alpha) ,
 \]
 which is the same as in \eqref{help1}.
This completes the proof.
\proofend


\noindent{\bf Proof of Proposition \ref{simple_aggregation}.}
The distribution of \ $X_0$ \ is a Poisson distribution with parameter
 \ $(1 - \alpha)^{-1} \lambda$, \ thus
 \ $\cov(X_0, X_0) = \var(X_0) = (1 - \alpha)^{-1} \lambda$.
\ By \eqref{help1_alter}, we have
 \begin{align*}
  \begin{aligned}
   F_{0,k}(x_0, x_k)
   = \EE(x_0^{X_0} x_k^{X_k})
    = F_{0,\ldots,k}(x_0, 1, \ldots, 1, x_k)
   = \ee^{(1-\alpha)^{-1}\lambda[\alpha^k(x_0-1)(x_k-1) +(x_0-1)+(x_k-1)]}
  \end{aligned}
 \end{align*}
 for all \ $x_0, x_k \in \RR$ \ and \ $k \in \NN$.
\ Consequently,
 \begin{align*}
   \EE(X_0 X_k)
   = \frac{\partial^2 F_{0,k}(x_0, x_k)}
          {\partial x_0 \partial x_k}\bigg|_{(x_0,x_k)=(1,1)}
   = \frac{\lambda\alpha^k}{1-\alpha}
     + \frac{\lambda^2}{(1-\alpha)^2} , \qquad k \in \NN ,
 \end{align*}
 since
 \begin{align*}
  \frac{\partial^2 F_{0,k}(x_0, x_k)}{\partial x_0 \partial x_k}
  = F_{0,k}(x_0, x_k) \frac{\lambda^2}{(1-\alpha)^2}
     (\alpha^k (x_0 - 1) + 1) (\alpha^k (x_k - 1) + 1)
     + F_{0,k}(x_0, x_k) \frac{\lambda}{1-\alpha} \alpha^k.
 \end{align*}
Hence we obtain the formula for \ $\cov(X_0 , X_k)$.
\ The statement follows from the multidimensional central limit theorem.
Due to the continuous mapping theorem, it is sufficient to show the convergence
 \ $N^{-1/2} (S^{(N)}_0, S^{(N)}_1, \ldots, S^{(N)}_k)
    \distr (\cX_0, \cX_1, \ldots, \cX_k)$
 \ as \ $N \to \infty$ \ for all \ $k \in \ZZ_+$.
\ For all \ $k \in \ZZ_+$, \ the random vectors
 \ $\bigl(X^{(j)}_0 - \frac{\lambda}{1-\alpha}, X^{(j)}_1-\frac{\lambda}{1-\alpha},
          \ldots,X^{(j)}_k - \frac{\lambda}{1-\alpha}\bigr)$,
 \ $j \in \NN$, \ are independent, identically distributed having zero expectation
 vector and covariances
 \[
   \cov(X^{(j)}_{\ell_1}, X^{(j)}_{\ell_2})
   = \cov(X^{(j)}_0, X^{(j)}_{|\ell_2-\ell_1|})
   = \frac{\lambda \alpha^{|\ell_2-\ell_1|}}{1-\alpha} ,
   \qquad j \in \NN , \quad \ell_1, \ell_2 \in \{0, 1, \ldots, k\}  ,
 \]
 following from the strict stationarity of \ $X^{(j)}$ \ and from the form of
 \ $\cov(X_0 , X_k)$.
\proofend


\noindent{\bf Proof of Proposition \ref{simple_aggregation2}.}
It is known that
 \[
   M_k := X_k - \EE(X_k \mid \cF^X_{k-1})
        = X_k - (X_{k-1} \EE(\xi_{1,1}) + \EE(\vare_1))
        = X_k - \alpha X_{k-1} - \lambda ,
   \qquad k \in \NN ,
 \]
 are martingale differences with respect to the filtration
 \ $\cF_k^X := \sigma(X_0, \ldots, X_k)$, \ $k \in \ZZ_+$ \ with
 \begin{align}\label{help14}
   \EE(M_k^2 \mid \cF^X_{k-1}) = X_{k-1} \var(\xi_{1,1}) + \var(\vare_1)
   = \alpha (1 - \alpha) X_{k-1} + \lambda, \qquad k \in \NN .
 \end{align}
The functional martingale central limit theorem can be applied, see, e.g., Jacod and
 Shiryaev \cite[Theorem VIII.3.33]{JacShi}.
Indeed, by ergodicity, for each \ $t \in \RR_+$, \ we have
 \[
   \frac{1}{n} \sum_{k=1}^{\lfloor nt \rfloor} \EE(M_k^2 \mid \cF^X_{k-1})
   \as \Bigl( \alpha (1 - \alpha) \frac{\lambda}{1-\alpha} + \lambda\Bigr) t
   =  (1 + \alpha ) \lambda t
   \qquad \text{as \ $n \to \infty$.}
 \]
Moreover, the conditional Lyapunov condition holds, namely, again by
 ergodicity,
 \[
   \frac{1}{n^2} \sum_{k=1}^{\lfloor nt \rfloor} \EE(M_k^4 \mid \cF^X_{k-1})
   \as 0 \qquad \text{as \ $n \to \infty$,}
 \]
 since there exists a second order polynomial \ $P$ \ such that
 \ $E(M_k^4 \mid \cF_{k-1}) = P(X_{k-1})$, \ $k \in \NN$, \ see Ned\'enyi
 \cite[Formula (8)]{Ned}, or Barczy et al.\ \cite[Lemma A.2, part (ii)]{BarIspPap}
 together with the decomposition
 \ $M_k = \sum_{j=1}^{X_{k-1}} (\xi_{k,j} - \EE(\xi_{k,j}))
          + (\vare_k - \EE(\vare_k))$,
 \ $k \in \NN$.
\ Hence we obtain
 \[
   \biggl(\frac{1}{\sqrt{n}} \sum_{k=1}^{\lfloor nt \rfloor} M_k\biggr)_{t\in\RR_+}
   \distr \sqrt{\lambda(1+\alpha)} \, B \qquad
   \text{as \ $n \to \infty$.}
 \]
We have \ $X_k = \alpha X_{k-1} + M_k + \lambda$, \ $k \in \NN$, \ thus
 \ $\EE(X_k) = \alpha \EE(X_{k-1}) + \lambda$, \ $k \in \NN$, \ and hence
 \ $X_k - \EE(X_k) = \alpha (X_{k-1} - \EE(X_{k-1})) + M_k$, \ $k \in \NN$,
 \ yielding
 \[
   X_k - \EE(X_k) = \alpha^k (X_0 - \EE(X_0)) + \sum_{j=1}^k \alpha^{k-j} M_j ,
   \qquad k \in \NN .
 \]
Consequently, for each \ $n \in \NN$ \ and \ $t \in \RR_+$,
 \begin{align*}
  \frac{1}{\sqrt{n}} \sum_{k=1}^{\lfloor nt \rfloor} (X_k - \EE(X_k))
  &= \frac{1}{\sqrt{n}} (X_0 - \EE(X_0)) \sum_{k=1}^{\lfloor nt \rfloor} \alpha^k
     + \frac{1}{\sqrt{n}}
       \sum_{k=1}^{\lfloor nt \rfloor} \sum_{j=1}^k \alpha^{k-j} M_j \\
  &= (X_0 - \EE(X_0)) \frac{\alpha-\alpha^{\lfloor nt \rfloor+1}}{(1-\alpha)\sqrt{n}}
     + \frac{1}{\sqrt{n}}
        \sum_{j=1}^{\lfloor nt \rfloor}
         M_j \sum_{k=j}^{\lfloor nt \rfloor} \alpha^{k-j} \\
  &= (X_0 - \EE(X_0)) \frac{\alpha-\alpha^{\lfloor nt \rfloor+1}}{(1-\alpha)\sqrt{n}}
     + \frac{1}{\sqrt{n}}
        \sum_{j=1}^{\lfloor nt \rfloor}
         M_j \frac{1-\alpha^{\lfloor nt \rfloor-j+1}}{1-\alpha} ,
 \end{align*}
 implying the statement using Slutsky's lemma.
Indeed,
 \ $n^{-1/2} \sum_{j=1}^{\lfloor nt \rfloor} \alpha^{\lfloor nt \rfloor-j+1} M_j$
 \ converges in \ $L_1$ \ and hence in probability to 0 as \ $n \to \infty$, \ since, by \eqref{help14},
 \[
   \EE(|M_j|)
   \leq \sqrt{\EE(M_j^2)}
   = \sqrt{\alpha (1 - \alpha) \EE(X_{j-1}) + \lambda}
   = \sqrt{\lambda (1 + \alpha)} ,
 \]
 and hence,
 \begin{align*}
  \EE\biggl(\biggl|\frac{1}{\sqrt{n}}
                  \sum_{j=1}^{\lfloor nt \rfloor}
                   \alpha^{\lfloor nt \rfloor-j+1} M_j\biggr|\biggr)
  \leq \frac{\sqrt{\lambda (1 + \alpha)}}{\sqrt{n}}
        \sum_{j=1}^{\lfloor nt \rfloor}
         \alpha^{\lfloor nt \rfloor-j+1}
  = \frac{\sqrt{\lambda (1 + \alpha)}}{\sqrt{n}} \,
     \frac{\alpha(1-\alpha^{\lfloor nt \rfloor})}{1-\alpha}
   \to 0
 \end{align*}
 as \ $n \to \infty$.
\proofend


\noindent{\bf Proof of Theorem \ref{double_aggregation}.}
For all \ $N, m \in \NN$ \ and all \ $t_1, \ldots, t_m \in\RR_+$, \ by Proposition
 \ref{simple_aggregation2} and by the continuity theorem, we have
 \[
   \frac{1}{\sqrt{n}} (S^{(N,n)}_{t_1}, \ldots, S^{(N,n)}_{t_m})
   \distr \frac{\sqrt{\lambda(1+\alpha)}}{1-\alpha}
          \sum_{j=1}^N
          (B^{(j)}_{t_1}, \ldots, B^{(j)}_{t_m})
   \qquad \text{as \ $n \to \infty$,}
 \]
 where \ $B^{(j)} = (B^{(j)}_t)_{t\in\RR_+}$, \ $j \in \{1, \ldots, N\}$, \ are
 independent standard Wiener processes.
Since
 \[
   \frac{1}{\sqrt{N}} \sum_{j=1}^N (B^{(j)}_{t_1}, \ldots, B^{(j)}_{t_m})
   \distre (B_{t_1}, \ldots, B_{t_m}) , \qquad N \in \NN ,
 \]
 we obtain the first convergence.

For all \ $n \in \NN$ \ and for all \ $t_1, \ldots, t_m \in \RR_+$ \ with
 \ $t_1 < \ldots < t_m$, \ $m \in \NN$, \ by Proposition \ref{simple_aggregation}
 and by the continuous mapping theorem, we have
 \begin{align*}
  \frac{1}{N^{1/2}} \bigl(S^{(N,n)}_{t_1} , \ldots, S^{(N,n)}_{t_m}\bigr)
  &\distr \Biggl(\sum_{k=1}^{\lfloor nt_1\rfloor} \cX_k , \ldots,
                 \sum_{k=1}^{\lfloor nt_m\rfloor} \cX_k\Biggr) \\
  &\distre \cN_m\Biggl(\bzero,
                       \var\Biggl(\sum_{k=1}^{\lfloor nt_1\rfloor} \cX_k , \ldots,
                                  \sum_{k=1}^{\lfloor nt_m\rfloor}
                                   \cX_k\Biggr)\Biggr) ,
  \qquad N \to \infty ,
 \end{align*}
 where \ $(\cX_k)_{k\in\ZZ_+}$ \ is the stationary Gaussian process given in
 Proposition \ref{simple_aggregation} and
 \[
   \var\Biggl(\sum_{k=1}^{\lfloor nt_1\rfloor} \cX_k , \ldots,
              \sum_{k=1}^{\lfloor nt_m\rfloor} \cX_k\Biggr)
   = \Biggl(\frac{\lambda}{1-\alpha}
            \sum_{k=1}^{\lfloor nt_i\rfloor}
            \sum_{\ell=1}^{\lfloor nt_j\rfloor}
             \alpha^{|\ell-k|} \Biggr)_{i,j\in\{1,\ldots,m\}} .
 \]
By the continuity theorem, for all \ $\theta_1, \ldots, \theta_m \in \RR$,
 \ $m \in \NN$, \ we conclude
 \begin{align*}
  &\lim_{N\to\infty}
    \EE\biggl(\exp\biggl\{\ii
                          \sum_{j=1}^m
                           \theta_j n^{-1/2} N^{-1/2}
                           S^{(N,n)}_{t_j}\biggr\}\biggr) \\
  &= \exp\biggl\{-\frac{\lambda}{2n(1-\alpha)}
                  \sum_{i=1}^m \sum_{j=1}^m
                   \theta_i \theta_j
                   \sum_{k=1}^{\lfloor nt_i\rfloor}
                    \sum_{\ell=1}^{\lfloor nt_j\rfloor}
                     \alpha^{|\ell-k|}\biggr\} \\
  &\to \exp\biggl\{-\frac{(1+\alpha)\lambda}{2(1-\alpha)^2}
                    \sum_{i=1}^m \sum_{j=1}^m
                     \theta_i \theta_j (t_i \land t_j)\biggr\}
 \end{align*}
 as \ $n \to \infty$.
\ Indeed, for all \ $s, t \in \RR_+$ \ with \ $s \leq t$, \ we have
 \begin{equation}\label{sum1}
\begin{split}
 &\frac{1}{n}
  \sum_{k=1}^{\lfloor ns\rfloor}
   \sum_{\ell=1}^{\lfloor nt\rfloor}
    \alpha^{|\ell-k|}
  = \frac{1}{n}
     \sum_{k=1}^{\lfloor ns\rfloor}
      \sum_{\ell=1}^{k-1}
       \alpha^{k-\ell}
     + \frac{\lfloor ns\rfloor}{n}
     + \frac{1}{n}
       \sum_{k=1}^{\lfloor ns\rfloor}
        \sum_{\ell=k+1}^{\lfloor nt\rfloor}
         \alpha^{\ell-k} \\
  &= \frac{1}{n}
     \sum_{k=1}^{\lfloor ns\rfloor}
      \frac{\alpha-\alpha^k}{1-\alpha}
     + \frac{\lfloor ns\rfloor}{n}
     + \frac{1}{n}
       \sum_{k=1}^{\lfloor ns\rfloor}
        \frac{\alpha-\alpha^{\lfloor nt\rfloor-k+1}}{1-\alpha} \\
  &
    = \frac{1}{n(1-\alpha)}
      \left(\lfloor ns \rfloor \alpha
            - \alpha \frac{1- \alpha^{\lfloor ns\rfloor}}{1-\alpha}\right)
      + \frac{\lfloor ns\rfloor}{n} \\
  & \quad
      + \frac{1}{n(1-\alpha)}
        \left(\lfloor ns\rfloor \alpha
              - \alpha^{\lfloor nt\rfloor - \lfloor ns\rfloor + 1}
                \frac{1- \alpha^{\lfloor ns\rfloor}}{1-\alpha} \right) \\
  &= \frac{1+\alpha}{1-\alpha} \, \frac{\lfloor ns\rfloor}{n}
     - \frac{\alpha}{(1-\alpha)^2n}
       (1+\alpha^{\lfloor nt\rfloor-\lfloor ns\rfloor})
       (1-\alpha^{\lfloor ns\rfloor})
   \to \frac{1+\alpha}{1-\alpha} \, s
 \end{split}
\end{equation}
 as \ $n \to \infty$.
\ This implies the second convergence.
\proofend


\noindent{\bf Proof of Proposition \ref{simple_aggregation_random}.}
We have
 \[
   \EE\Bigl(X_k - \frac{\lambda}{1-\alpha}\Bigr)
   = \EE\Bigl[\EE\Bigl(X_k - \frac{\lambda}{1-\alpha} \, \Big| \, \alpha\Bigr)\Bigr]
   = 0 , \qquad k \in \ZZ_+ ,
 \]
 and hence, for all \ $k \in  \ZZ_+$,
 \begin{align}\label{help15}
 \begin{split}
  &\cov\Bigl(X_0 - \frac{\lambda}{1-\alpha}, X_k - \frac{\lambda}{1-\alpha}\Bigr)
   = \EE\Bigl[\Bigl(X_0 - \frac{\lambda}{1-\alpha}\Bigr)
              \Bigl(X_k - \frac{\lambda}{1-\alpha}\Bigr)\Bigr] \\
  &= \EE\Bigl\{\EE\Bigl[\Bigl(X_0 - \frac{\lambda}{1-\alpha}\Bigr)
                        \Bigl(X_k - \frac{\lambda}{1-\alpha}\Bigr)
                        \,\Big|\, \alpha\Bigr]\Bigr\}
   = \EE\Bigl(\frac{\lambda\alpha^k}{1-\alpha}\Bigr) ,
 \end{split}
 \end{align}
 where we applied \eqref{cov}.
Now the statement follows from the multidimensional central limit theorem in the same
 way as in the proof of Proposition \ref{simple_aggregation}.
\proofend


\noindent{\bf Proof of Proposition \ref{simple_aggregation_random3}.}
For each \ $n \in \NN$ \ and each \ $t \in \RR_+$, \ put
 \[
   \tT^{(n)}_t := \frac{1}{\sqrt{n}} \sum_{k=1}^{\lfloor nt \rfloor} \tS^{(1)}_k .
 \]
For each \ $m \in \NN$, \ each \ $t_1, \ldots, t_m \in \RR_+$, \ and each bounded
 continuous function \ $g : \RR^m \to \RR$, \ we have
 \begin{align*}
  &\EE(g(\tT^{(n)}_{t_1}, \ldots, \tT^{(n)}_{t_m}))
   = \int_0^1
      \EE(g(\tT^{(n)}_{t_1}, \ldots, \tT^{(n)}_{t_m}) \mid \alpha = a)
      \, \PP_\alpha(\dd a) \\
  &= \int_0^1
      \EE\biggl(g\biggl(\frac{1}{\sqrt{n}}
                        \sum_{k=1}^{\lfloor nt_1 \rfloor}
                         \Bigl(X_k - \frac{\lambda}{1-a}\Bigr) , \ldots ,
                        \frac{1}{\sqrt{n}}
                        \sum_{k=1}^{\lfloor nt_m \rfloor}
                         \Bigl(X_k - \frac{\lambda}{1-a}\Bigr)\biggr)
                \;\bigg|\; \alpha = a\biggr)
      \PP_\alpha(\dd a) .
 \end{align*}
Proposition \ref{simple_aggregation2}, the portmanteau theorem and the boundedness of \ $g$ \ justify the
 usage of the dominated convergence theorem, and we obtain
 \begin{align*}
  &\lim_{n\to\infty} \EE(g(\tT^{(n)}_{t_1}, \ldots, \tT^{(n)}_{t_m}))
   = \int_0^1
      \EE\biggl(g\biggl(\frac{\sqrt{\lambda(1+a)}}{1-a} B_{t_1}, \ldots,
                        \frac{\sqrt{\lambda(1+a)}}{1-a} B_{t_m}\biggr)\biggr)
      \PP_\alpha(\dd a) \\
  &= \int_0^1
      \EE\biggl(g\biggl(\frac{\sqrt{\lambda(1+\alpha)}}{1-\alpha} B_{t_1}, \ldots,
                        \frac{\sqrt{\lambda(1+\alpha)}}{1-\alpha} B_{t_m}\biggr)
                \bigg| \alpha = a\biggr)
      \PP_\alpha(\dd a) \\
  &= \EE\biggl(g\biggl(\frac{\sqrt{\lambda(1+\alpha)}}{1-\alpha} B_{t_1}, \ldots,
                       \frac{\sqrt{\lambda(1+\alpha)}}{1-\alpha}
                       B_{t_m}\biggr)\biggr) ,
 \end{align*}
 hence we obtain the statement by the portmanteau theorem.
\proofend


\noindent{\bf Proof of Proposition \ref{simple_aggregation_random2}.}
For all \ $k\in\ZZ_+$, \ by the strict stationarity of \ $(X_k)_{k\in\ZZ_+}$ \ and \eqref{help15},
 we have
 \begin{align}\label{help16}
 \begin{split}
  \cov(X_0 , X_k)
  &= \EE\left[\left(X_0 - \EE\left(\frac{\lambda}{1-\alpha}\right)\right)
              \left(X_k - \EE\left(\frac{\lambda}{1-\alpha}\right)\right)\right] \\
  &= \EE\left[\left(X_0 - \frac{\lambda}{1-\alpha}\right)
              \left(X_k - \frac{\lambda}{1-\alpha}\right)\right]
    + \EE\left[\left(\frac{\lambda}{1-\alpha}
                     - \EE\left(\frac{\lambda}{1-\alpha}\right)\right)^2\right] \\
  &= \lambda \EE\left(\frac{\alpha^k}{1-\alpha}\right)
     + \lambda^2 \var\left(\frac{1}{1-\alpha}\right) ,
 \end{split}
 \end{align}
 since
 \begin{align*}
  &\EE\left[\left(X_k - \frac{\lambda}{1-\alpha}\right)
            \left(\frac{\lambda}{1-\alpha}
                  - \EE\left(\frac{\lambda}{1-\alpha}\right)\right)\right] \\
  &= \EE\left\{\EE\left[\left(X_k - \frac{\lambda}{1-\alpha}\right)
                       \left(\frac{\lambda}{1-\alpha}
                             - \EE\left(\frac{\lambda}{1-\alpha}\right)\right)
                       \,\bigg|\, \alpha\right]\right\} \\
  &= \EE\left\{\left(\frac{\lambda}{1-\alpha}
                     - \EE\left(\frac{\lambda}{1-\alpha}\right)\right)
               \EE\left(X_k - \frac{\lambda}{1-\alpha}
                        \,\bigg|\, \alpha\right)\right\}
   = 0
 \end{align*}
 for all \ $k \in \ZZ_+$.

The statement follows from the multidimensional central limit theorem as in the
 proof of Proposition \ref{simple_aggregation}.
Indeed, for all \ $k\in\ZZ_+$, \ the random vectors
 \[
   \left(X^{(j)}_0 - \lambda \EE\left(\frac{1}{1-\alpha}\right),
         X^{(j)}_1 - \lambda \EE\left(\frac{1}{1-\alpha}\right),
         \ldots, X^{(j)}_k - \lambda \EE\left(\frac{1}{1-\alpha}\right)\right),
   \quad j \in \NN ,
 \]
 are independent, identically distributed having zero expectation vector and
 covariances
 \[
   \cov(X^{(j)}_{\ell_1}, X^{(j)}_{\ell_2})
   = \cov(X^{(j)}_0, X^{(j)}_{|\ell_2 -\ell_1|})
   = \lambda \EE\left(\frac{\alpha^{|\ell_2 - \ell_1|}}{1-\alpha}\right)
     + \lambda^2 \var\left(\frac{1}{1-\alpha}\right)
 \]
 for \ $j \in \NN$ \ and \ $\ell_1, \ell_2 \in \{0, 1, \ldots, k\}$, \ following
 from the strict stationarity of \ $X^{(j)}$ \ and from the form of
 \ $\cov(X_0, X_k)$ \ given in \eqref{help16}.
\proofend


\noindent{\bf Proof of Proposition \ref{simple_aggregation_random4}.}
We have a decomposition \ $S^{(1)}_k = \tS^{(1)}_k + R^{(1)}_k$, \ $k \in \ZZ_+$,
 \ with
 \[
   R^{(1)}_k := \EE(X^{(1)}_k \mid \alpha^{(1)}) - \EE(X^{(1)}_k)
   = \frac{\lambda}{1-\alpha^{(1)}}
     - \EE\biggl(\frac{\lambda}{1-\alpha^{(1)}}\biggr) , \qquad k \in \ZZ_+ .
 \]
We have
 \begin{align*}
  \Biggl(\frac{1}{n} \sum_{k=1}^{\lfloor nt \rfloor} R^{(1)}_k\Biggr)_{t\in\RR_+}
  &=\Biggl(\frac{\lfloor nt \rfloor}{n}
           \Biggl(\frac{\lambda}{1-\alpha^{(1)}}
                  -\EE\biggl(\frac{\lambda}
                                  {1-\alpha^{(1)}}\biggr)\Biggr)\Biggr)_{t\in\RR_+}\\
  &\distrf
  \biggl(\biggl(\frac{\lambda}{1-\alpha}
                - \EE\biggl(\frac{\lambda}{1-\alpha}\biggr)\biggr) t\biggr)_{t\in\RR_+}
 \end{align*}
 as \ $n \to \infty$.
Moreover, by Proposition \ref{simple_aggregation_random3},
 \ $\cD_\ff$-$\lim_{n\to\infty}
              \bigl(n^{-1/2}
                    \sum_{k=1}^{\lfloor nt \rfloor} \tS^{(1)}_k\bigr)_{t\in\RR_+}$
 \ exists, hence
 \[
   \Biggl(\frac{1}{n}
          \sum_{k=1}^{\lfloor nt \rfloor}
           \tS^{(1)}_k\Biggr)_{t\in\RR_+}
   \distrf 0 \qquad \text{as \ $n \to \infty$,}
 \]
 implying that for all \ $m \in \NN$ \ and all \ $t_1, \ldots, t_m \in \RR_+$, \ we
 have
 \[
   \Biggl(\frac{1}{n} \sum_{k=1}^{\lfloor nt_1 \rfloor} \tS^{(1)}_k , \ldots,
          \frac{1}{n} \sum_{k=1}^{\lfloor nt_m \rfloor} \tS^{(1)}_k\Biggr)
   \stoch 0 \qquad \text{as \ $n \to \infty$.}
 \]
By Slutsky's lemma we conclude the statement.
\proofend


\noindent{\bf First proof of Theorem \ref{joint_aggregation_random}.}
By Remark \ref{condexp}, condition \ $\beta \in (0, 1)$ \ implies
 \ $\EE\bigl(\frac{1}{1-\alpha}\bigr) < \infty$.
\ Hence, by Proposition \ref{simple_aggregation_random} and the continuous mapping
 theorem, it suffices to show that
 \begin{equation}\label{jar}
  \cD_\ff\text{-}\hspace*{-1mm}\lim_{n\to\infty} \,
  \Biggl(\frac{1}{n^{1-\frac{\beta}{2}}}
         \sum_{k=1}^{\lfloor nt \rfloor} \tcY_k\Biggr)_{t\in\RR_+}
  = \sqrt{\frac{2\lambda\psi_1\Gamma(\beta)}{(2-\beta)(1-\beta)}}
    \cB_{1-\frac{\beta}{2}} .
 \end{equation}
We are going to apply Theorem 4.3 in Beran et al.\ \cite{BerFenGhoKul} with \ $m = 1$ \ for
 the strictly stationary Gaussian process
 \ $\Bigl(\tcY_k / \sqrt{\var(\tcY_0)}\Bigr)_{k\in\ZZ_+}$, \ where, by \eqref{covariance},
 \[
   \var(\tcY_0) = \lambda \EE\left(\frac{1}{1-\alpha}\right) , \qquad
   \cov(\tcY_0, \tcY_k) = \lambda \EE\left(\frac{\alpha^k}{1-\alpha}\right) ,
   \qquad k \in \ZZ_+ ,
 \]
 hence
 \[
   \cov\left(\frac{\tcY_0}{\sqrt{\var(\tcY_0)}},
             \frac{\tcY_k}{\sqrt{\var(\tcY_0)}}\right)
   = \frac{\EE\left(\frac{\alpha^k}{1-\alpha}\right)}
          {\EE\left(\frac{1}{1-\alpha}\right)} ,
   \qquad k \in \ZZ_+ .
 \]
In order to check the conditions of Theorem 4.3 in Beran et al.\ \cite{BerFenGhoKul}, first we show that
 \begin{equation}\label{beta_distribution}
  k^\beta \EE\left(\frac{\alpha^k}{1-\alpha}\right)
  = k^\beta \int_0^1 a^k (1 - a)^{\beta-1} \psi(a) \, \dd a
  \to \psi_1 \Gamma(\beta) \qquad \text{as \ $k \to \infty$,}
 \end{equation}
meaning that the covariance function of the process \ $(\tcY_k)_{k \in \ZZ_+}$ \ is regularly varying with index \ $-\beta$.
\ First note that, by Stirling's formula,
 \begin{align*}
  \lim_{k \to \infty} k^\beta \int_0^1 a^k (1 - a)^{\beta-1} \psi_1 \, \dd a
   &= \lim_{k \to \infty}
       \psi_1 \frac{k^\beta \Gamma(k+1) }{\Gamma(k+\beta+1)} \, \Gamma(\beta) \\
   &= \psi_1 \Gamma(\beta)
      \lim_{k \to \infty}
       \sqrt{\frac{k}{k+\beta}} \left(\frac{k}{k+\beta}\right)^{k+\beta}\ee^\beta
    = \psi_1 \Gamma(\beta) .
 \end{align*}
Next, for arbitrary \ $\delta \in (0, \psi_1)$, \ there exists \ $\vare \in (0, 1)$
 \ such that \ $|\psi(a) - \psi_1| \leq \delta$ \ for all \ $a \in [1 - \vare, 1)$,
 \ and hence
 \[
   k^\beta \int_{1-\vare}^1 a^k (1 - a)^{\beta-1} |\psi(a) - \psi_1| \, \dd a
   \leq \delta \sup_{k \in \NN} k^\beta \int_0^1 a^k (1 - a)^{\beta-1} \, \dd a
 \]
 can be arbitrary small.
Further, observe
 \begin{align*}
  &k^\beta \int_0^{1-\vare} a^k (1 - a)^{\beta-1} \psi(a) \, \dd a
   \leq \frac{k^\beta (1-\vare)^k}{\vare}
        \int_0^{1-\vare} (1 - a)^\beta \psi(a) \, \dd a \\
  &\leq \frac{k^\beta (1-\vare)^k}{\vare}
         \int_0^1 (1 - a)^\beta \psi(a) \, \dd a
   = \frac{k^\beta (1-\vare)^k}{\vare}
   \to 0 \qquad \text{as \ $k \to \infty$.}
 \end{align*}
In a similar way, we have
 \begin{align*}
  &k^\beta \int_0^{1-\vare} a^k (1 - a)^{\beta-1} \psi_1 \, \dd a
   \leq \psi_1 k^\beta (1-\vare)^k
        \int_0^{1-\vare} (1 - a)^{\beta-1} \, \dd a \\
  &\leq \psi_1 k^\beta (1-\vare)^k
        \int_0^1 (1 - a)^{\beta-1} \, \dd a
   = \psi_1 \, \frac{k^\beta (1-\vare)^k}{\beta}
   \to 0 \qquad \text{as \ $k \to \infty$,}
 \end{align*}
 hence
 \[
   k^\beta \int_0^{1-\vare} a^k (1 - a)^{\beta-1} |\psi(a) - \psi_1| \, \dd a
   \to 0 \qquad \text{as \ $k \to \infty$,}
 \]
 implying \eqref{beta_distribution}.
Applying \eqref{beta_distribution}, we conclude
 \[
   k^\beta
   \cov\left(\frac{\tcY_0}{\sqrt{\var(\tcY_0)}},
             \frac{\tcY_k}{\sqrt{\var(\tcY_0)}}\right)
   = k^\beta
     \frac{\EE\left(\frac{\alpha^k}{1-\alpha}\right)}
          {\EE\left(\frac{1}{1-\alpha}\right)}
   \to \frac{\psi_1 \Gamma(\beta)}{\EE\left(\frac{1}{1-\alpha}\right)}
 \]
 as \ $k \to \infty$.
\ Consequently, by Theorem 4.3 in Beran et al.\ \cite{BerFenGhoKul},
 \[
   \left(\frac{1}{n^{1-\frac{\beta}{2}} L_1(n)^{1/2}}
         \sum_{k=1}^{\lfloor nt \rfloor}
          \frac{\tcY_k}
               {\sqrt{\lambda\EE\left(\frac{1}{1-\alpha}\right)}}\right)_{t\in\RR_+}
    \distr \cZ_{1,1-\frac{\beta}{2}}\distre \
   \cB_{1-\frac{\beta}{2}} , \qquad \text{as \ $n \to \infty$,}
 \]
where \ $\cZ_{1,1-\frac{\beta}{2}}$ \ is the Hermite-Rosenblatt process defined in Definition 3.24 of Beran et al.\ \cite{BerFenGhoKul}, and
 \[
   L_1(n) = \frac{\psi_1 \Gamma(\beta)}{\EE\left(\frac{1}{1-\alpha}\right)} \, C_1 ,
   \qquad n \in \NN , \qquad \text{with} \qquad
   C_1 = \frac{2}{(1 - \beta)(2 - \beta)}.
 \]
The fact that the Hermite-Rosenblatt process \ $\cZ_{1,1-\frac{\beta}{2}}$ \ coincides in law with \ $\cB_{1-\frac{\beta}{2}}$ \ is shown in Beran et al.\ \cite{BerFenGhoKul}, see Definition 3.23, the representation in formula (3.111), and page 195 of \cite{BerFenGhoKul} for details.
Hence we obtain the statement.
\proofend


\noindent{\bf Second proof of Theorem \ref{joint_aggregation_random}.}
As in the first proof of Theorem \ref{joint_aggregation_random}, it suffices to show
 \eqref{jar}.
As for every \ $n \in \NN$ \ the process
 \ $n^{-1+\beta/2} \sum_{k=1}^{\lfloor nt\rfloor} \tcY_k$, \ $t\in\RR_+$, \ is
 Gaussian, so is the limit process.
Also, it is clear that both processes have zero mean.
Therefore, it suffices to show that the covariance function of
 \ $n^{-1+\beta/2} \sum_{k=1}^{\lfloor nt\rfloor} \tcY_k$, \ $t\in\RR_+$,
 \ converges to that of the limit process in \eqref{jar}.
By \eqref{covariance}, the covariance function of
 \ $n^{-1+\beta/2} \sum_{k=1}^{\lfloor nt\rfloor} \tcY_k$, \ $t\in\RR_+$, \ for any
 \ $t_1, t_2\in \RR_+$ \ takes the form
 \[
   \cov\Biggl(n^{-1+\beta/2}
              \sum_{k=1}^{\lfloor nt_1\rfloor}
               \tcY_k,n^{-1+\beta/2}
               \sum_{\ell=1}^{\lfloor nt_2\rfloor}\tcY_\ell\Biggr)
   = n^{-2+\beta} \lambda
     \EE\Biggl(\sum_{k=1}^{\lfloor nt_1\rfloor}\sum_{\ell=1}^{\lfloor nt_2\rfloor}
                \frac{\alpha^{|k-\ell|}}{1-\alpha}\Biggr) .
 \]
By \eqref{sum1}, for time points \ $0 \leq t_1 \leq t_2$, \ we get
 \begin{equation}\label{help12}
  \begin{split}
   &\frac{1}{1-\alpha}
    \sum_{k=1}^{\lfloor nt_1\rfloor} \sum_{\ell=1}^{\lfloor nt_2\rfloor}
     \alpha^{|k-\ell|}
    = \frac{(1-\alpha^2)\lfloor nt_1\rfloor
            - \alpha \left(1-\alpha^{\lfloor nt_2\rfloor}
            - \alpha^{\lfloor nt_1\rfloor}
            + \alpha^{\lfloor nt_2\rfloor-\lfloor nt_1\rfloor}\right)}
           {(1-\alpha)^3} \\
   &= \frac{\alpha (\alpha^{\lfloor nt_2 \rfloor}-1)
            + \lfloor nt_2\rfloor (1-\alpha^2)/2}
           {(1-\alpha)^3}
      + \frac{\alpha (\alpha^{\lfloor nt_1 \rfloor}-1)
              + \lfloor nt_1\rfloor (1-\alpha^2)/2}
             {(1-\alpha)^3} \\
   &\quad
      - \frac{\alpha(\alpha^{\lfloor nt_2 \rfloor-\lfloor nt_1 \rfloor}-1)
      + (\lfloor nt_2\rfloor-\lfloor nt_1 \rfloor) (1-\alpha^2)/2}{(1-\alpha)^3} .
  \end{split}
 \end{equation}
We are going to show that for any \ $0 \leq t_1 \leq t_2$ \ we get
 \begin{equation*}
  \begin{split}
   &n^{-2+\beta}
    \EE\biggl(\frac{\alpha(\alpha^{\lfloor nt_2 \rfloor-\lfloor nt_1 \rfloor}-1)
                    + (\lfloor nt_2\rfloor-\lfloor nt_1 \rfloor) (1-\alpha^2)/2}
                   {(1-\alpha)^3}\biggr) \\
   &\to \psi_1
        \int_0^\infty
         \left(\ee^{-y(t_2-t_1)} - 1 + y (t_2 - t_1)\right) y^{\beta-3} \, \dd y
    = \frac{\psi_1\Gamma(\beta)}{(2-\beta)(1-\beta)} \, (t_2-t_1)^{2-\beta} ,
  \end{split}
 \end{equation*}
 as \ $n\to \infty$, \ where the equality follows with repeated partial integration.
This will imply that the limit of the sequence of the covariance functions in question is
 \begin{equation*}
  \frac{2\lambda \psi_1 \Gamma(\beta)}{(2-\beta)(1-\beta)} \,
  \frac{t_1^{2-\beta}+t_2^{2-\beta}-|t_2-t_1|^{2-\beta}}{2} ,
 \end{equation*}
 that is the covariance function of a fractional Brownian motion with parameter
 \ $1 - \beta/2$ \ multiplied by
 \ $\sqrt{2\lambda\psi_1\Gamma(\beta)/((2-\beta)(1-\beta))}$, \ as desired.

By substituting \ $a = 1 - y/n$ \ we get that
 \begin{equation*}
  \begin{split}
  &n^{-2+\beta}
   \EE\biggl(\frac{\alpha(\alpha^{\lfloor nt_2 \rfloor-\lfloor nt_1 \rfloor}-1)
                   +(\lfloor nt_2\rfloor-\lfloor nt_1 \rfloor) (1-\alpha^2)/2}
                  {(1-\alpha)^3}\biggr) \\
  &= n^{-2+\beta}
     \int_0^1
      \frac{a(a^{\lfloor nt_2 \rfloor-\lfloor nt_1 \rfloor}-1)
            +(\lfloor nt_2 \rfloor-\lfloor nt_1 \rfloor) (1-a^2)/2}
           {(1-a)^3}
      (1-a)^\beta \psi(a) \, \dd a \\
  &= n^{-2+\beta}
     \int_0^n
      \biggl[\left(1 - \frac{y}{n}\right)
             \left(\left(1 - \frac{y}{n}\right)^{\lfloor nt_2\rfloor-\lfloor nt_1\rfloor}-1\right) \\
  &\phantom{=n^{-2+\beta}\int_0^n \biggl[}
             + (\lfloor nt_2 \rfloor - \lfloor nt_1 \rfloor)
               \left(1-\left(1-\frac{y}{n}\right)^2\right)\frac{1}{2}\biggr]
      \left(\frac{y}{n}\right)^{\beta-3}
      \psi \left(1-\frac{y}{n}\right) \frac{\dd y}{n} \\
  &= \int_0^n D_n(y) \, \dd y
\end{split}
\end{equation*}
 with
 \begin{align*}
  D_n(y)
  := \biggl[\left(1-\frac{y}{n}\right)
             \left(\left(1 - \frac{y}{n}\right)^{\lfloor nt_2\rfloor-\lfloor nt_1\rfloor}-1\right)
             + \frac{\lfloor nt_2 \rfloor-\lfloor nt_1 \rfloor}{n} \,
               y\left(1-\frac{y}{2n}\right)\biggr]
      y^{\beta-3} \psi \left(1 - \frac{y}{n}\right)
 \end{align*}
 for \ $y \in [0, n]$.
\ First note that, for any \ $\vare \in (0, 1)$ \ and \ $n > 1/\vare$, \ we have
 \begin{align}\label{help9}
  \begin{split}
   \biggl|\int_{n\vare}^n D_n(y) \, \dd y\biggr|
   &\leq \int_{n\vare}^n
          \left(1\cdot2 + \left(t_2-t_1+\frac{1}{n}\right)y\right)
          y^{\beta-3} \psi\left(1 - \frac{y}{n}\right) \, \dd y \\
   &\leq \int_{n\vare}^n
          (1\cdot2 y + (t_2 - t_1 + 1) y) y^{\beta-3}
          \psi\left(1 - \frac{y}{n}\right) \, \dd y \\
   &= (t_2 - t_1 + 3)
      \int_0^{1-\vare} \big(n(1-a)\big)^{\beta-2} \psi(a) n \, \dd a \\
   &= (t_2 - t_1 + 3) n^{\beta-1}
      \int_0^{1-\vare}(1-a)^{\beta-2} \psi(a) \, \dd a \\
   &\leq (t_2 - t_1 + 3) n^{\beta-1} \vare^{-2}
         \int_0^1 (1 - a)^\beta \psi(a) \, \dd a
    \to 0, \qquad n \to \infty .
  \end{split}
 \end{align}
We are going to show that
 \begin{equation*}
  \begin{split}
   \int_0^{n\varepsilon} D_n(y) \, \dd y
   \to \psi_1
       \int_0^\infty
        \left(\ee^{-y(t_2-t_1)} - 1 + y (t_2 - t_1)\right) y^{\beta-3} \, \dd y ,
   \qquad n \to \infty .
  \end{split}
 \end{equation*}
The pointwise convergence is evident, and we can give a dominating integrable
 function proving the above convergence.

Note that
 \[
   \left(1 - \frac{y}{n}\right)^{\lfloor nt_2\rfloor-\lfloor nt_1\rfloor}
   = 1 - (\lfloor nt_2 \rfloor - \lfloor nt_1 \rfloor) \frac{y}{n}
     + \sum_{k=2}^{\lfloor nt_2\rfloor-\lfloor nt_1\rfloor}
        \binom{\lfloor nt_2\rfloor-\lfloor nt_1\rfloor}{k}
        \left(-\frac{y}{n}\right)^k,
 \]
 where for any \ $y \in [0, 1]$ \ we have
 \begin{align*}
  &\left|\sum_{k=2}^{\lfloor nt_2\rfloor-\lfloor nt_1\rfloor}
          \binom{\lfloor nt_2 \rfloor-\lfloor nt_1 \rfloor}{k}
          \left(-\frac{y}{n}\right)^k\right| \\
  &\leq y^2
        \sum_{k=2}^{\lfloor nt_2\rfloor-\lfloor nt_1\rfloor}
         \frac{y^{k-2}}{k!}
         \frac{(\lfloor nt_2\rfloor-\lfloor nt_1\rfloor)\cdots
               (\lfloor nt_2\rfloor-\lfloor nt_1\rfloor-k+1)}
              {n^k} \\
  &\leq y^2 \sum_{k=0}^\infty \frac{(t_2-t_1+1)^k}{k!}
   = y^2 \ee^{t_2-t_1+1} .
 \end{align*}
Choose \ $\vare \in (0, 1)$ \ such that for every \ $x \in (1 - \vare, 1)$ \ we have
 \ $\psi(x) \leq 2 \psi_1$.
\ Then for any \ $n > 1/\vare$, \ applying Bernoulli's inequality, we obtain
 \begin{align*}
  \int_{0}^1 |D_n(y)| \, \dd y
  &= \int_0^1
      \biggl|\left[\left(1 - \frac{y}{n}\right)^{\lfloor nt_2\rfloor
                                                 -\lfloor nt_1\rfloor}
                   - 1
                   + (\lfloor nt_2\rfloor-\lfloor nt_1\rfloor)\frac{y}{n}\right] \\
  &\phantom{=\int_0^1\biggl|}
     - \frac{y}{n}
       \left[\left(1 - \frac{y}{n}\right)^{\lfloor nt_2\rfloor-\lfloor nt_1\rfloor}
             - 1
             + \frac{y(\lfloor nt_2\rfloor-\lfloor nt_1 \rfloor)}{2n}\right]\biggr|
       y^{\beta-3}\psi \left(1-\frac{y}{n}\right) \, \dd y \\
  &\leq \int_0^1
         \left(y^2\ee^{t_2-t_1+1}
               + \frac{y}{n}
                 \left(y (t_2 - t_1 + 1) + y \frac{t_2-t_1+1}{2}\right)\right)
                 y^{\beta-3} \psi\left(1 - \frac{y}{n}\right) \, \dd y \\
  &\leq 2 \psi_1 (\ee^{t_2 - t_1 + 1} + 2 (t_2 - t_1 + 1))
        \int_0^1 y^{\beta-1} \, \dd y
   < \infty .
 \end{align*}
Similarly to \eqref{help9}, for any \ $n > 1/\vare$, \ one gets
 \begin{equation*}
  \begin{split}
   &\int_1^{n\vare} |D_n(y)| \, \dd y
    \leq \int_1^{n\vare}
          (2 + t_2 - t_1 + 1) y^{\beta-2} \psi\left(1 - \frac{y}{n}\right)
          \, \dd y \\
   &\leq (t_2 - t_1 + 3) 2 \psi_1 \int_1^{n\vare} y^{\beta-2} \, \dd y
    \leq (t_2 - t_1 + 3) 2 \psi_1 \int_1^\infty y^{\beta-2} \, \dd y
    < \infty .
  \end{split}
 \end{equation*}
So the function
 \[
   2 \psi_1 (\ee^{t_2-t_1+1} + 2 (t_2 - t_1 + 1)) y^{\beta-1} \bbone_{[0,1)}(y)
   + (t_2 - t_1 + 3) 2 \psi_1 y^{\beta-2} \bbone_{[1,\infty)}(y)
 \]
 can be chosen as a dominating integrable function.
\proofend


\noindent{\bf Proof of Theorem \ref{joint_aggregation_random_2}.}
To prove this limit theorem it is enough to show that for any \ $n \in \NN$,
 \[
   \cD_\ff\text{-}\hspace*{-1mm}\lim_{N\to\infty} \,
   N^{-\frac{1}{2(1+\beta)}} \, \tS^{(N,n)}
   = (\lfloor nt\rfloor V_{2(1+\beta)})_{t\in\RR_+} .
 \]
For this, by the continuous mapping theorem, it is enough to verify that for any
 \ $m \in \NN$,
 \begin{equation*}
  \begin{split}
   &N^{-\frac{1}{2(1+\beta)}}
    \sum_{j=1}^N
     \left(X_1^{(j)}-\frac{\lambda}{1-\alpha^{(j)}}, \dots,
           X_m^{(j)}-\frac{\lambda}{1-\alpha^{(j)}}\right)
   \distr V_{2(1+\beta)}(1, \dots, 1)
  \end{split}
 \end{equation*}
 as \ $N \to \infty$.
\ So, by the continuity theorem, we have to check that for any \ $m \in \NN$ \ and
 \ $\theta_1, \dots, \theta_m \in \RR$ \ the convergence
 \begin{equation*}
  \begin{split}
   &\EE\left(\exp\left\{\ii
                        \sum_{k=1}^m
                         \theta_k
                         \left(N^{-\frac{1}{2(1+\beta)}}
                               \sum_{j=1}^N
                                \left(X_k^{(j)}
                                      - \frac{\lambda}
                                             {1-\alpha^{(j)}}\right)\right)\right\}
       \right) \\
   &\qquad\qquad
    = \EE\left(\exp\left\{\ii
                          N^{-\frac{1}{2(1+\beta)}}
                          \sum_{j=1}^N \sum_{k=1}^m
                           \theta_k
                           \left(X_k^{(j)}
                                 - \frac{\lambda}{1-\alpha^{(j)}}\right)\right\}
         \right) \\
   &\qquad\qquad
    = \left[\EE\left(\exp\left\{\ii
                                N^{-\frac{1}{2(1+\beta)}}
                                \sum_{k=1}^m
                                 \theta_k
                                 \left(X_k
                                       - \frac{\lambda}{1-\alpha}\right)\right\}
               \right)\right]^N \\
   &\qquad\qquad
    \to \EE\left(\ee^{\ii \sum_{k=1}^m \theta_k V_{2(1+\beta)}}\right)
    = \ee^{-K_\beta |\sum_{k=1}^m \theta_k|^{2(1+\beta)}} \qquad
    \text{as \ $N \to \infty$}
  \end{split}
 \end{equation*}
 holds.
Note that it suffices to show
 \[
   \Theta_N
   := N
      \left[1 - \EE\left(\exp\left\{\ii N^{-\frac{1}{2(1+\beta)}}
                               \sum_{k=1}^m
                                \theta_k
                                \left(X_k
                                      - \frac{\lambda}
                                             {1-\alpha}\right)\right\}\right)
           \right]
   \to K_\beta \Biggl|\sum_{k=1}^m \theta_k\Biggr|^{2(1+\beta)}
 \]
 as \ $N \to \infty$, \ since it implies that
 \ $(1 - \Theta_N/N)^N
    \to \ee^{-K_\beta |\sum_{k=1}^m \theta_k|^{2(1+\beta)}}$
 \ as \ $N \to \infty$.
\ By applying \eqref{help1_alter} to the left hand side, we get
 \begin{equation*}
  \begin{split}
   \Theta_N
   &= N \EE\Bigl[1 - F_{0,\dots,m-1}\Bigl(\ee^{\ii N^{-\frac{1}{2(1+\beta)}}\theta_1},
                                      \dots,
                                      \ee^{\ii N^{-\frac{1}{2(1+\beta)}}\theta_m}
                                       \,\Big|\, \alpha\Bigr)
      \ee^{-\ii N^{-\frac{1}{2(1+\beta)}}\frac{\lambda}{1-\alpha}
            \sum_{k=1}^m\theta_k}
      \Bigr]\\
   &= N \EE\left[1 - \ee^{\frac{\lambda}{1-\alpha}A_N(\alpha)}\right]
    = N \int_0^1
         \left(1 - \ee^{\frac{\lambda}{1-a}A_N(a)}\right) \psi(a) (1-a)^\beta
         \, \dd a,
  \end{split}
 \end{equation*}
  where \ $F_{0,\ldots,m-1}(z_0,\ldots,z_{m-1}\mid\alpha):= \EE(z_0^{X_0}z_1^{X_1}\cdots z_{m-1}^{X_{m-1}}\mid\alpha)$,
  \ $z_0,\ldots,z_{m-1}\in\CC$, \ and
 \begin{align*}
  A_N(a)
  &:= - \frac{\ii(\theta_1+\cdots+\theta_m)}{N^{\frac{1}{2(1+\beta)}}} \\
  &\phantom{:=}
      + \sum_{1\leq\ell\leq j\leq m}
         a^{j-\ell}
         \bigl(\ee^{\ii N^{-\frac{1}{2(1+\beta)}}\theta_\ell} - 1\bigr)
         \ee^{\ii N^{-\frac{1}{2(1+\beta)}}(\theta_{\ell+1}+\dots+\theta_{j-1})}
         \bigl(\ee^{\ii N^{-\frac{1}{2(1+\beta)}}\theta_j} - 1\bigr)
 \end{align*}
 for \ $a \in [0, 1]$.
\ Let us show that for any \ $\vare \in (0, 1)$ \ we have
 \ $\sup_{a\in(0,1-\vare)} |N A_N(a)| \to 0$ \ as \ $N \to \infty$.
\ Using \eqref{exp:3}, for any \ $\vare \in (0, 1)$ \ we get
 \begin{equation*}
  \begin{split}
   &\sup_{a\in(0,1-\vare)} N |A_N(a)|
    = \sup_{a\in(0,1-\vare)}
       N \Bigg|\sum_{k=1}^m
                \Bigl(\ee^{\ii N^{-\frac{1}{2(1+\beta)}}\theta_k} - 1
                      - \ii N^{-\frac{1}{2(1+\beta)}}\theta_k\Bigr) \\
   &\qquad\qquad
      + \sum_{1\leq \ell <j \leq m} a^{j-\ell}
         \Bigl(\ee^{\ii N^{-\frac{1}{2(1+\beta)}}\theta_\ell} - 1\Bigr)
         \ee^{\ii N^{-\frac{1}{2(1+\beta)}}(\theta_{\ell+1}+\dots+\theta_{j-1})}
         \Bigl(\ee^{\ii N^{-\frac{1}{2(1+\beta)}}\theta_j} - 1\Bigr)\Bigg| \\
   &\leq N \left(\sum_{k=1}^m
                  N^{-\frac{1}{1+\beta}} \frac{\theta_k^2}{2}
                 + \sum_{1\leq\ell<j\leq m}
                    N^{-\frac{1}{1+\beta}} {|\theta_\ell| |\theta_j|}\right)
    = N^{\frac{\beta}{1+\beta}}
      \frac{\left(\sum_{k=1}^m |\theta_k| \right)^2}{2}
    \to 0
  \end{split}
 \end{equation*}
 as \ $N \to \infty$, \ since \ $\beta/(1+\beta) < 0$.
\ Therefore, by Lemma \ref{ordo}, substituting \ $a = 1 - z^{-1} N^{-\frac{1}{1+\beta}}$,
 \ the statement of the theorem will follow from
 \begin{gather}\label{limsup_an}
  \begin{aligned}
   &\limsup_{N\to\infty}
     N \int^1_{1-\vare}
        \Bigl|1 - \ee^{\frac{\lambda}{1-a}A_N(a)}\Bigr| (1-a)^\beta \, \dd a \\
   &= \limsup_{N\to\infty}
       \int^\infty_{\vare^{-1}N^{-\frac{1}{1+\beta}}}
        \Bigl|1 - \ee^{\lambda zN^{\frac{1}{1+\beta}}
                       A_N\bigl(1 - z^{-1}N^{-\frac{1}{1+\beta}}\bigr)}\Bigr|
        z^{-(2+\beta)} \, \dd z
    < \infty
  \end{aligned}
 \end{gather}
 for all \ $\vare \in (0, 1)$ \ and
 \begin{gather}\label{lim_an}
  \begin{aligned}
   &\lim_{\vare\downarrow0}
     \limsup_{N\to\infty}
      \left|N \int_{1-\vare}^1
               \Bigl(1 - \ee^{\frac{\lambda}{1-a}A_N(a)}\Bigr) (1-a)^\beta \, \dd a
            - I\right| \\
   &= \lim_{\vare\downarrow0}
       \limsup_{N\to\infty}
        \left|\int^\infty_{\vare^{-1}N^{-\frac{1}{1+\beta}}}
               \Bigl(1 - \ee^{\lambda zN^{\frac{1}{1+\beta}}
                              A_N\bigl(1 - z^{-1}N^{-\frac{1}{1+\beta}}\bigr)}\Bigr)
               z^{-(2+\beta)} \, \dd z
             - I\right|
    = 0
  \end{aligned}
 \end{gather}
 with
 \begin{align*}
  I &:= \int_0^\infty
         \biggl(1-\ee^{-\frac{\lambda z}{2}
                        \bigl(\sum_{k=1}^m \theta_k\bigr)^2}\biggr)
         z^{-(2+\beta)}
         \, \dd z \\
    &= \left(\frac{\lambda}{2} \left|\sum_{k=1}^m \theta_k\right|^2\right)^{1+\beta}
        \int_0^\infty (1-\ee^{-z}) z^{-(2+\beta)} \, \dd z
     = \psi_1^{-1} K_\beta \left|\sum_{k=1}^m \theta_k\right|^{2(1+\beta)} ,
 \end{align*}
 where the last equality is justified by Lemma 2.2.1 in Zolotarev \cite{Zol} (be
 careful for the misprint in \cite{Zol}: a negative sign is superfluous) or by
 Li \cite[formula (1.28)]{Li}.
Next we check \eqref{limsup_an} and \eqref{lim_an}.

By Taylor expansion,
 \begin{gather*}
  \ee^{\ii N^{-\frac{1}{2(1+\beta)}} \theta_\ell} - 1
  = \ii N^{-\frac{1}{2(1+\beta)}} \theta_\ell + N^{-\frac{1}{1+\beta}} \OO(1)
  = N^{-\frac{1}{2(1+\beta)}} \OO(1), \\
  \ee^{\ii N^{-\frac{1}{2(1+\beta)}} \theta_\ell} - 1
  - \ii N^{-\frac{1}{2(1+\beta)}}\theta_{\ell}
  = - N^{-\frac{1}{1+\beta}} \frac{\theta_{\ell}^2}{2}
    + N^{-\frac{3}{2(1+\beta)}} \OO(1)
 \end{gather*}
 for all \ $\ell \in\{1, \ldots, m\}$, \ resulting
 \begin{equation}\label{an}
  \lambda z N^{\frac{1}{1+\beta}} A_N\biggl(1 - \frac{1}{zN^{\frac{1}{1+\beta}}}\biggr)
  = - \frac{\lambda z\bigl(\sum_{k=1}^m\theta_k\bigr)^2}{2}
    + \frac{z\OO(1)}{N^{\frac{1}{2(1+\beta)}}}
    + \frac{\OO(1)}{N^{\frac{1}{1+\beta}}}
 \end {equation}
 for \ $z > N^{-\frac{1}{1+\beta}}$.
\ Indeed, for \ $z > N^{-\frac{1}{1+\beta}}$, \ we have
 \begin{equation*}
  \begin{split}
   &A_N\left(1-\frac{1}{zN^{\frac{1}{1+\beta}}}\right) \\
   &= \sum_{k=1}^m
       \bigl(\ee^{\ii N^{-\frac{1}{2(1+\beta)}}\theta_k} - 1
             - \ii N^{-\frac{1}{2(1+\beta)}}\theta_k\bigr) \\
   &\quad
      + \sum_{1\leq\ell<j\leq m}
         \left(1 - \frac{1}{zN^{\frac{1}{1+\beta}}}\right)^{j-\ell}
         \bigl(\ee^{\ii N^{-\frac{1}{2(1+\beta)}}\theta_\ell} - 1\bigr) \,
         \ee^{\ii N^{-\frac{1}{2(1+\beta)}}(\theta_{\ell+1}+\dots+\theta_{j-1})}
         \bigl(\ee^{\ii N^{-\frac{1}{2(1+\beta)}}\theta_j} - 1\bigr) \\
   &= \sum_{k=1}^m
       \biggl(- \frac{\theta_k^2}{2N^{\frac{1}{1+\beta}}}
              + \frac{\OO(1)}{N^{\frac{3}{2(1+\beta)}}}\biggr) \\
   &\quad
      + \sum_{1\leq\ell<j\leq m}
         \left(1 + \frac{\OO(1)}{zN^{\frac{1}{1+\beta}}}\right)
         \biggl(\frac{\ii\theta_\ell}{N^{\frac{1}{2(1+\beta)}}}
                + \frac{\OO(1)}{N^{\frac{1}{1+\beta}}}\biggr)
         \biggl(1 + \frac{\OO(1)}{N^{\frac{1}{2(1+\beta)}}}\biggr)
         \biggl(\frac{\ii\theta_j}{N^{\frac{1}{2(1+\beta)}}}
                + \frac{\OO(1)}{N^{\frac{1}{1+\beta}}}\biggr) \\
   &= - \frac{\sum_{k=1}^m \theta_k^2}{2N^{\frac{1}{1+\beta}}}
      + \frac{\OO(1)}{N^{\frac{3}{2(1+\beta)}}}
      - \frac{\sum_{1\leq\ell<j\leq m}\theta_\ell\theta_j}{N^{\frac{1}{1+\beta}}}
      + \frac{\OO(1)}{N^{\frac{3}{2(1+\beta)}}}
      + \frac{\OO(1)}{zN^{\frac{2}{1+\beta}}} \\
   &= - \frac{\bigl(\sum_{k=1}^m\theta_k\bigr)^2}{2N^{\frac{1}{1+\beta}}}
      + \frac{\OO(1)}{N^{\frac{3}{2(1+\beta)}}}
      + \frac{\OO(1)}{zN^{\frac{2}{1+\beta}}} ,
  \end{split}
 \end{equation*}
 since by Bernoulli's inequality
 \[
   \Bigg|\left(1 - \frac{1}{zN^{\frac{1}{1+\beta}}}\right)^{j-\ell} - 1\Bigg|
   \leq \frac{j-\ell}{zN^{\frac{1}{1+\beta}}}
   \leq \frac{m}{zN^{\frac{1}{1+\beta}}} ,
 \]
 yielding that
 \[
   \left(1 - \frac{1}{zN^{\frac{1}{1+\beta}}}\right)^{j-\ell}
   = 1 + \frac{\OO(1)}{zN^{\frac{1}{1+\beta}}} .
 \]
By \eqref{an}, for \ $z \in [1, \infty)$ \ and for large enough \ $N$ \ we have
 \begin{align*}
  \lambda z N^{\frac{1}{1+\beta}} \Re A_N\bigl(1 - z^{-1} N^{-\frac{1}{1+\beta}}\bigr)
  &= - \frac{\lambda z (\sum_{k=1}^m\theta_k)^2}{2} \,
       \left(1 - \frac{\Re\OO(1)}{N^{\frac{1}{2(1+\beta)}}}\right)
     + \frac{\Re\OO(1)}{N^{\frac{1}{1+\beta}}} \\
  &\leq - \frac{\lambda z (\sum_{k=1}^m\theta_k)^2}{4}
        + \frac{|\OO(1)|}{N^{\frac{1}{1+\beta}}}
   \leq 0 ,
 \end{align*}
 hence we obtain
 \begin{align}\label{help10A}
  \begin{split}
   &\int_1^\infty
     \left|1 - \ee^{\lambda z N^{\frac{1}{1+\beta}}
                    A_N(1 - z^{-1} N^{-\frac{1}{1+\beta}})}\right|
     z^{-(\beta+2)} \, \dd z \\
   &\leq \int_1^\infty
          \left(1 + \ee^{\lambda z N^{\frac{1}{1+\beta}}
                         \Re A_N(1 - z^{-1} N^{-\frac{1}{1+\beta}})}\right)
          z^{-(\beta+2)} \, \dd z
    \leq 2 \int_1^\infty z^{-(\beta+2)} \, \dd z
    < \infty .
  \end{split}
 \end{align}
Again by \eqref{an}, for \ $\vare \in (0, 1)$,
 \ $z \in \bigl(\vare^{-1} N^{-\frac{1}{1+\beta}}, 1\bigr]$ \ and for
 large enough \ $N$, \ we have
 \begin{align*}
  &\Bigl|\lambda z N^{\frac{1}{1+\beta}}
        A_N\Bigl(1 - z^{-1} N^{-\frac{1}{1+\beta}}\Bigr)\Bigr|
   \leq \frac{\lambda z (\sum_{k=1}^m\theta_k)^2}{2}
        + \frac{z|\OO(1)|}{N^{\frac{1}{2(1+\beta)}}}
        + \frac{|\OO(1)|}{N^{\frac{1}{1+\beta}}} \\
  &\leq z \biggl(\frac{\lambda(\sum_{k=1}^m\theta_k)^2}{2}
                 + \frac{|\OO(1)|}{N^{\frac{1}{2(1+\beta)}}}
                 + \vare |\OO(1)|\biggr)
   \leq  z |\OO(1)|
   \leq |\OO(1)| ,
 \end{align*}
 since \ $N^{-\frac{1}{1+\beta}} < z \vare$.
\ Hence, using \eqref{exp:1}, we obtain
 \begin{align*}
   &\int_{\vare^{-1} N^{-\frac{1}{1+\beta}}}^1
     \left|1 - \ee^{\lambda z N^{\frac{1}{1+\beta}}
                    A_N\bigl(1 - z^{-1} N^{-\frac{1}{1+\beta}}\bigr)}\right|
     z^{-(2+\beta)} \, \dd z \\
   &\leq \int_{\vare^{-1} N^{-\frac{1}{1+\beta}}}^1
          \left|\lambda z N^{\frac{1}{1+\beta}}
                A_N\bigl(1 - z^{-1} N^{-\frac{1}{1+\beta}}\bigr)\right|
          \ee^{\left|\lambda z N^{\frac{1}{1+\beta}}
                     A_N\bigl(1 - z^{-1} N^{-\frac{1}{1+\beta}}\bigr)\right|}
          z^{-(2+\beta)} \, \dd z \\
   &\leq |\OO(1)| \ee^{|\OO(1)|}
         \int_0^1 z^{-(1+\beta)} \, \dd z
     < \infty ,
 \end{align*}
 which, together with \eqref{help10A}, imply \eqref{limsup_an}.

Now we turn to prove \eqref{lim_an}.
By \eqref{exp:4}, we have
 \begin{equation*}
  \begin{split}
   &\left|\int_0^{\vare^{-1} N^{-\frac{1}{1+\beta}}}
           {\Bigl(1 - \ee^{-\frac{\lambda z}{2}
                          (\sum_{k=1}^m \theta_k)^2}\Bigr)}
           z^{-(2+\beta)} \, \dd z\right|
    \leq \int_0^{\vare^{-1} N^{-\frac{1}{1+\beta}}}
          {\frac{\lambda z(\sum_{k=1}^m \theta_k)^2}{2}}
          z^{-(2+\beta)} \, \dd z \\
   &= \frac{\lambda(\sum_{k=1}^m \theta_k)^2}{2}
            \int_0^{\varepsilon^{-1}
             N^{-\frac{1}{1+\beta}}} z^{-(1+\beta)}
             \, \dd z
          = \frac{\lambda(\sum_{k=1}^m \theta_k)^2}{2(-\beta)}
            \left(\frac{1}{\vare N^{\frac{1}{1+\beta}}}\right)^{-\beta}
    \to 0
  \end{split}
 \end{equation*}
 as \ $N \to \infty$, \ hence \eqref{lim_an} reduces to check that
 \ $\lim_{\vare\downarrow0} \limsup_{N\to\infty} I_{N,\vare} = 0$, \ where
 \begin{align*}
  I_{N,\vare}
  := \int^\infty_{\vare^{-1}N^{-\frac{1}{1+\beta}}}
      \Bigl[\ee^{\lambda z N^{\frac{1}{1+\beta}}
                 A_N(1 - z^{-1} N^{-\frac{1}{1+\beta}})}
            - \ee^{-\frac{\lambda z}{2}(\sum_{k=1}^m \theta_k)^2}\Bigr]
      z^{-(2+\beta)} \, \dd z .
 \end{align*}
Applying again \eqref{an}, we obtain
 \[
   |I_{N,\vare}|
   \leq \int^\infty_{\vare^{-1}N^{-\frac{1}{1+\beta}}}
         \ee^{-\frac{\lambda z}{2}(\sum_{k=1}^m \theta_k)^2}
          \Bigl|\ee^{z N^{-\frac{1}{2(1+\beta)}} \OO(1)
                     + N^{-\frac{1}{1+\beta}} \OO(1)}
                - 1\Bigr|
          z^{-(2+\beta)} \, \dd z .
 \]
Here, for \ $\vare \in (0, 1)$ \ and
 \ $z \in (\vare^{-1} N^{-\frac{1}{1+\beta}}, \infty)$, \ we have
 \[
   \bigl|z N^{-\frac{1}{2(1+\beta)}} \OO(1) + N^{-\frac{1}{1+\beta}} \OO(1)\bigr|
   \leq z \bigl(N^{-\frac{1}{2(1+\beta)}} + \vare\bigr) |\OO(1)| ,
 \]
 and hence, by \eqref{exp:1}, we get
 \begin{align*}
  \Bigl|\ee^{z N^{-\frac{1}{2(1+\beta)}}\OO(1)+N^{-\frac{1}{1+\beta}}\OO(1)}-1\Bigr|
  &\leq \bigl|z N^{-\frac{1}{2(1+\beta)}} \OO(1) + N^{-\frac{1}{1+\beta}} \OO(1)\bigr|\,
        \ee^{\bigl|zN^{-\frac{1}{2(1+\beta)}}\OO(1)+N^{-\frac{1}{1+\beta}}\OO(1)\bigr|}\\
  &\leq z \bigl(N^{-\frac{1}{2(1+\beta)}} + \vare\bigr) |\OO(1)| \,
     \ee^{z \bigl(N^{-\frac{1}{2(1+\beta)}} + \vare\bigr) |\OO(1)|} .
 \end{align*}
Consequently, for large enough \ $N$ \ and small enough \ $\vare\in(0,1)$,
 \begin{align*}
  |I_{N,\vare}|
  &\leq \bigl(N^{-\frac{1}{2(1+\beta)}} + \vare\bigr) |\OO(1)|
        \int^\infty_{\vare^{-1}N^{-\frac{1}{1+\beta}}}
         \ee^{-\frac{\lambda z}{2}(\sum_{k=1}^m \theta_k)^2
              +z \bigl(N^{-\frac{1}{2(1+\beta)}} + \vare\bigr) |\OO(1)|}
         z^{-(1+\beta)} \, \dd z \\
  &\leq \bigl(N^{-\frac{1}{2(1+\beta)}} + \vare\bigr) |\OO(1)|
        \int_0^\infty
         \ee^{-\frac{\lambda z}{4}(\sum_{k=1}^m \theta_k)^2}
         z^{-(1+\beta)} \, \dd z ,
 \end{align*}
 that gets arbitrarily close to zero as \ $N$ \ approaches infinity and
 \ $\vare$ \ tends to \ $0$, \ since the integral is finite due to the
 fact that
 \[
   \Gamma(-\beta) \left(\frac{\lambda}{4} \left(\sum_{k=1}^m\theta_k\right)^2\right)^\beta
    \ee^{-\lambda z {(\sum_{k=1}^m\theta_k)^2/4}} \, z^{-(1+\beta)},
    \qquad z > 0,
  \]
 is the density function of a Gamma distributed random variable with
 parameters \ $-\beta$ \ and \ $\lambda(\sum_{k=1}^m\theta_k)^2/4$.
\ This yields \eqref{lim_an} completing the proof.
\proofend


\noindent{\bf Proof of Theorem \ref{joint_aggregation_random_3}.}
Similarly as in the proof of Theorem \ref{joint_aggregation_random_2}, it suffices
 to show that for any \ $m \in \NN$ \ and \ $\theta_1, \dots, \theta_m \in \RR$ \ we
 have the convergence
 \[
   N\left[1-\EE\left(\exp\left\{\frac{\ii}{\sqrt{N\log N}}
                              \sum_{k=1}^m
                               \theta_k
                               \left(X_k
                                     -\frac{\lambda}{1-\alpha}\right)\right\}\right)
    \right]
   \to
   \frac{\lambda\psi_1}{2} \left(\sum_{k=1}^m \theta_k\right)^2
 \]
 as \ $N \to \infty$.
\ By applying \eqref{help1_alter}, the left hand side equals
 \begin{equation*}
  \begin{split}
   &N\EE\left[1-F_{0,\dots,m-1}\left(\ee^{\frac{\ii\theta_1}{\sqrt{N\log N}}},
                                     \dots,
                                     \ee^{\frac{\ii\theta_m}{\sqrt{N\log N}}}
                                     \,\Big|\, \alpha\right)
                                     \ee^{-\frac{\ii\lambda
                                                 (\theta_1+\cdots+\theta_m)}
                                                {(1-\alpha)\sqrt{N\log N}}}\right]\\
   &= N \EE\left[1 - \ee^{\frac{\lambda}{1-\alpha}B_N(\alpha)}\right]
    = N \int_0^1
         \left(1 - \ee^{\frac{\lambda}{1-a}B_N(a)}\right) \psi(a) \, \dd a
  \end{split}
 \end{equation*}
 with
 \begin{align*}
  B_N(a)
  &:= \sum_{k=1}^m
       \left(\ee^{\frac{\ii\theta_k}{\sqrt{N\log N}}} - 1
             - \frac{\ii\theta_k}{\sqrt{N\log N}}\right) \\
  &\phantom{:=}
      + \sum_{1\leq\ell<j\leq m}
         a^{j-\ell}
         \bigl(\ee^{\frac{\ii\theta_\ell}{\sqrt{N\log N}}} - 1\bigr)
         \ee^{\frac{\ii(\theta_{\ell+1}+\cdots+\theta_{j-1})}{\sqrt{N\log N}}}
         \bigl(\ee^{\frac{\ii\theta_j}{\sqrt{N\log N}}} - 1\bigr) , \qquad
  a \in [0, 1] .
 \end{align*}
Just like in the proof of Theorem \ref{joint_aggregation_random_2} it is easy to see
 that for any \ $\vare \in (0, 1)$ \ we have
 \[
   \sup_{a\in(0,1-\vare)} |NB_N(a)|
   \leq \frac{(\sum_{k=1}^m \theta_k)^2}{2\log N}
   \to 0
 \]
 as \ $N \to \infty$.
\ Therefore, by Lemma \ref{ordo}, substituting \ $a = 1 - z/N$, \ the statement of
 the theorem will follow from
 \begin{gather}\label{limsup_bn}
   \limsup_{N\to\infty}
     N \int_{1-\vare}^1
        \Bigl|1 - \ee^{\frac{\lambda}{1-a}B_N(a)}\Bigr| \, \dd a
   = \limsup_{N\to\infty}
       \int_0^{\vare N}
        \Bigl|1 - \ee^{\frac{\lambda N}{z} B_N(1-\frac{z}{N})}\Bigr| \, \dd z
    < \infty ,
 \end{gather}
 and
 \begin{gather}\label{lim_bn}
  \begin{aligned}
   \lim_{N\to\infty}
     N \int_{1-\vare}^1
        \Bigl(1 - \ee^{\frac{\lambda}{1-a}B_N(a)}\Bigr) \, \dd a
   = \lim_{N\to\infty}
       \int_0^{\vare N}
        \Bigl(1 - \ee^{\frac{\lambda N}{z} B_N(1-\frac{z}{N})}\Bigr) \, \dd z
    = \frac{\lambda}{2} \left(\sum_{k=1}^m \theta_k\right)^2
  \end{aligned}
 \end{gather}
 for all \ $\vare \in (0, 1)$.
\ Next we check \eqref{limsup_bn} and \eqref{lim_bn}.

Using Taylor expansions, similarly as in the proof of Theorem
 \ref{joint_aggregation_random_2}, we get
 \begin{equation}\label{bn}
  \frac{\lambda N}{z} B_N\Bigl(1-\frac{z}{N}\Bigr)
  = - \frac{\lambda(\sum_{k=1}^m\theta_k)^2}{2z\log N}
    + \frac{\OO(1)}{zN^{1/2}(\log N)^{3/2}}
    + \frac{\OO(1)}{N\log N} .
 \end {equation}
Indeed, for \ $z \in [0, N]$ \ we have
 \begin{align*}
  B_N\Bigl(1-\frac{z}{N}\Bigr)
  &= \sum_{k=1}^m
      \Big(\ee^{\frac{\ii\theta_k}{\sqrt{N\log N}}} - 1
           - \frac{\ii\theta_k}{\sqrt{N\log N}}\Big) \\
  &\quad
     + \sum_{1\leq\ell<j\leq m}
        \left(1 - \frac{z}{N}\right)^{j-\ell}
        \bigl(\ee^{\frac{\ii\theta_\ell}{\sqrt{N\log N}}} - 1\bigr)
        \, \ee^{\frac{\ii(\theta_{\ell+1}+\dots+\theta_{j-1})}{\sqrt{N\log N}}}
        \bigl(\ee^{\frac{\ii\theta_j}{\sqrt{N\log N}}} - 1\bigr) \\
  &= \sum_{k=1}^m
      \biggl(- \frac{\theta_k^2}{2N\log N}
            + \frac{\OO(1)}{(N\log N)^{3/2}}\biggr) \\
  &\quad
     +\sum_{1\leq\ell<j\leq m}
       \biggl(1 + \frac{z\OO(1)}{N}\biggr)
       \biggl(\frac{\ii\theta_\ell}{\sqrt{N\log N}}
              + \frac{\OO(1)}{N\log N}\biggr) \\
  &\phantom{=+\sum_{1\leq\ell<j\leq m}}
       \times
       \biggl(1 + \frac{\OO(1)}{\sqrt{N\log N}}\biggr)
       \biggl(\frac{\ii\theta_j}{\sqrt{N\log N}}
              + \frac{\OO(1)}{N\log N}\biggr) \\
  &= - \frac{\sum_{k=1}^m\theta_k^2}{2N\log N}
     + \frac{\OO(1)}{(N\log N)^{3/2}}
     - \frac{\sum_{1\leq\ell<j\leq m}\theta_\ell\theta_j}{N\log N}
     + \frac{\OO(1)}{(N\log N)^{3/2}}
     + \frac{z\OO(1)}{N^2\log N} \\
  &= - \frac{\left(\sum_{k=1}^m\theta_k\right)^2}{2N\log N}
     + \frac{\OO(1)}{(N\log N)^{3/2}}
     + \frac{z\OO(1)}{N^2\log N} ,
 \end{align*}
 since, by Bernoulli's inequality,
 \[
   \bigg|\left(1 - \frac{z}{N}\right)^{j-\ell} - 1\bigg|
   \leq (j - \ell) \frac{z}{N}
   \leq m \frac{z}{N} ,
 \]
 yielding that
 \[
   \left(1 - \frac{z}{N}\right)^{j-\ell} = 1 + \frac{z}{N} \OO(1) .
 \]
By \eqref{bn}, for \ $z \in \bigl(0, \frac{1}{\log N}\bigr)$ \ and for large enough
 \ $N$ \ we have
 \begin{align*}
  &\frac{\lambda N}{z} \Re B_N\Bigl(1-\frac{z}{N}\Bigr)
   = - \frac{\lambda(\sum_{k=1}^m\theta_k)^2}{2z\log N} \,
       \left(1 - \frac{\Re\OO(1)}{\sqrt{N\log N}}\right)
     + \frac{\Re\OO(1)}{N\log N} \\
  &\leq - \frac{\lambda(\sum_{k=1}^m\theta_k)^2}{4z\log N} + \frac{|\OO(1)|}{N\log N}
   \leq - \frac{\lambda(\sum_{k=1}^m\theta_k)^2}{4} + \frac{|\OO(1)|}{N\log N} ,
 \end{align*}
 hence we obtain
 \begin{align}\label{help10}
  \begin{split}
   &\int_0^{\frac{1}{\log N}}
     \left|1 - \ee^{\frac{\lambda N}{z}B_N(1-\frac{z}{N})}\right| \dd z
    \leq \int_0^{\frac{1}{\log N}}
          \left(1 + \ee^{\frac{\lambda N}{z}\Re B_N(1-\frac{z}{N})}\right) \dd z \\
   &\leq \frac{1}{\log N}
         \biggl(1 + \exp\biggl\{-\frac{\lambda\left(\sum_{k=1}^m\theta_k\right)^2}{4}
                                + \frac{|\OO(1)|}{N\log N}\biggr\}\biggr)
    \to 0 \qquad \text{as \ $N \to \infty$.}
  \end{split}
 \end{align}
Note that
 \begin{gather}\label{help11}
  \frac{1}{\log N}
  \int_{\frac{1}{\log N}}^{\vare N} \frac{1}{z} \, \dd z
  = \frac{\log \vare + \log N + \log\log N}{\log N}
  \to 1
  \qquad \text{as \ $N \to \infty$,} \\
  \label{help11+}
  \frac{1}{\log N}
  \int_{\frac{1}{\log N}}^{\vare N} \frac{1}{z^2} \, \dd z
  = \frac{\vare N \log N - 1}{\vare N \log N}
  \to 1
  \qquad \text{as \ $N \to \infty$.}
 \end{gather}
By \eqref{bn}, for all \ $z \in \bigl(\frac{1}{\log N}, \vare N\bigr)$, \ we have
 \begin{equation}\label{BN2}
  \left|\frac{\lambda N}{z} B_N\left(1 - \frac{z}{N}\right)\right|
  \leq \frac{\lambda(\sum_{k=1}^m\theta_k)^2}{2z\log N}
       + \frac{|\OO(1)|}{zN^{1/2}(\log N)^{3/2}}
       + \frac{|\OO(1)|}{N\log N}
  = |\OO(1)| .
 \end{equation}
Thus, by \eqref{exp:1} and \eqref{help11}, we get
 \begin{equation*}
  \begin{split}
   &\limsup_{N\to\infty}
     \int_{\frac{1}{\log N}}^{\vare N}
      \left|1 - \ee^{\frac{\lambda N}{z}B_N(1-\frac{z}{N})}\right| \dd z \\
   &\leq \limsup_{N\to\infty}
          \int_{\frac{1}{\log N}}^{\vare N}
           \left|\frac{\lambda N}{z}B_N(1-\frac{z}{N})\right|
           \ee^{\left|\frac{\lambda N}{z}B_N(1-\frac{z}{N})\right|}
           \, \dd z \\
   &\leq \limsup_{N\to\infty}
          \ee^{|\OO(1)|}
          \int_{\frac{1}{\log N}}^{\vare N}
           \left[\frac{\lambda(\sum_{k=1}^m \theta_k)^2}{2z\log N}
                 + \frac{|\OO(1)|}{zN^{1/2}(\log N)^{3/2}}
                 + \frac{|\OO(1)|}{N\log N}\right]
            \dd z
    < \infty ,
  \end{split}
 \end{equation*}
 which, together with \eqref{help10}, imply \eqref{limsup_bn}.

Now we turn to prove \eqref{lim_bn}.
By \eqref{help10}, the convergence \eqref{lim_bn} reduces to prove that
 \begin{align*}
  \left|\int_{\frac{1}{\log N}}^{\vare N}
         \left(1 - \ee^{\frac{\lambda N}{z}B_N(1-\frac{z}{N})}\right) \dd z
        - \frac{\lambda\left(\sum_{k=1}^m\theta_k\right)^2}{2}\right|
  \to 0 \qquad \text{as \ $N \to \infty$.}
 \end{align*}
Using \eqref{help11}, it is enough to check that
 \begin{align*}
  \left|\int_{\frac{1}{\log N}}^{\vare N}
         \left(\ee^{\frac{\lambda N}{z}B_N(1-\frac{z}{N})} - 1
               + \frac{\lambda\left(\sum_{k=1}^m\theta_k\right)^2}
                      {2z\log N}\right) \dd z\right|
  \to 0 \qquad \text{as \ $N \to \infty$.}
 \end{align*}
By applying \eqref{exp:2}, \eqref{bn} and \eqref{BN2}, for large enough \ $N$ \ we
 get
 \begin{align*}
  &\left|\int_{\frac{1}{\log N}}^{\vare N}
          \left[\left(\ee^{\frac{\lambda N}{z}B_N(1-\frac{z}{N})} - 1\right)
                +\frac{\lambda\left(\sum_{k=1}^m \theta_k\right)^2}
                      {2z\log N} \right] \dd z \right| \\
  &\leq \int_{\frac{1}{\log N}}^{\vare N}
         \Bigg[\frac{1}{2}
               \left|\frac{\lambda N}{z} B_N\Bigl(1 - \frac{z}{N}\Bigr)\right|^2
               \ee^{\left|\frac{\lambda N}{z}B_N(1-\frac{z}{N})\right|}
   + \left|\frac{\lambda N}{z} B_N\Bigl(1-\frac{z}{N}\Bigr)
           + \frac{\lambda\left(\sum_{k=1}^m\theta_k\right)^2}
                  {2z\log N} \right|\Bigg] \, \dd z \\
  &\leq \int_{\frac{1}{\log N}}^{\vare N}
     \Bigg[\frac{1}{2}\left(\frac{\lambda \left(\sum_{k=1}^m \theta_k\right)^2}
                                 {2z\log N}
                            + \frac{|\OO(1)|}{zN^{1/2}(\log N)^{3/2}}
                            + \frac{|\OO(1)|}{N\log N}\right)^2
           \ee^{|\OO(1)|}\\
 &\phantom{\leq \int_{\frac{1}{\log N}}^{\vare N} \Bigg[}
           + \frac{|\OO(1)|}{zN^{1/2}(\log N)^{3/2}}
           + \frac{|\OO(1)|}{N\log N}\Bigg] \dd z \\
 &\leq \int_{\frac{1}{\log N}}^{\vare N}
        \Bigg[\frac{3}{2} \left(\frac{|\OO(1)|}{z^2(\log N)^2}
                                + \frac{|\OO(1)|}{z^2N(\log N)^3}
                                +  \frac{|\OO(1)|}{N^2(\log N)^2}\right)
              +\frac{|\OO(1)|}{zN^{1/2}(\log N)^{3/2}}
              +\frac{|\OO(1)|}{N\log N}\Bigg] \dd z   ,
 \end{align*}
 which converges to \ $0$ \ as \ $N \to \infty$ \ using \eqref{help11} and
 \eqref{help11+}.
This yields \eqref{lim_bn} completing the proof.
\proofend


\noindent{\bf First proof of Theorem \ref{joint_aggregation_random_4}.}
By Proposition \ref{simple_aggregation_random3}, we have
 \begin{align*}
  \cD_\ff\text{-}\hspace*{-1mm}\lim_{n\to\infty} \,
  \Bigg(n^{-\frac{1}{2}} \,
        \sum_{k=1}^{\lfloor nt\rfloor}
         (X^{(1)}_k - \EE(X^{(1)}_k \mid \alpha^{(1)}))\Bigg)_{t\in\RR_+}
   = \frac{\sqrt{\lambda(1+\alpha)}}{1-\alpha} B ,
 \end{align*}
 where \ $(B_t)_{t\in\RR_+}$ \ is a standard Wiener process and \ $\alpha$ \ is a
 random variable having a density function of the form \eqref{alpha} with
 \ $\beta \in(-1, 1)$ \ and \ $\psi_1 \in (0, \infty)$, \ and being independent of
 \ $B$.
\ Let \ $\cW_t := \frac{\sqrt{\lambda(1+\alpha)}}{1-\alpha} B_t$, \ $t \in \RR_+$,
 \ and \ $(\cW^{(i)}_t)_{t\in\RR_+}$, \ $i \in \NN$, \ be its independent copies.
It remains to prove that
 \[
   \cD_\ff\text{-}\hspace*{-1mm}\lim_{N\to\infty}
    \biggl(N^{-\frac{1}{1+\beta}} \sum_{i=1}^N \cW^{(i)}_t\biggr)_{t\in\RR_+}
   = \cY_{1+\beta} .
 \]
Using the continuity theorem and the continuous mapping theorem, it is enough to prove that
 for all \ $m \in \NN$, \ $\theta_1, \ldots, \theta_m \in \RR$ \ and
 \ $0 =: t_0 < t_1 < t_2 < \cdots < t_m$,
 \begin{align*}
  &\EE\Biggl(\exp\Biggl\{\ii \sum_{j=1}^m
                              \theta_j
                              \Biggl(N^{-\frac{1}{1+\beta}}
                              \sum_{i=1}^N
                               (\cW^{(i)}_{t_j}
                                - \cW^{(i)}_{t_{j-1}})\Biggr)\Biggr\}\Biggr)
  = \Biggl[\EE\Biggl(\exp\Biggl\{\ii N^{-\frac{1}{1+\beta}}
                                  \sum_{j=1}^m
                                   \theta_j
                                   (\cW_{t_j}
                                    - \cW_{t_{j-1}})\Biggr\}\Biggr)\Biggr]^N \\
  &\to \EE\Biggl(\exp\Biggl\{\ii \sum_{j=1}^m
                                  \theta_j
                                  (\cY_{1+\beta}(t_j)
                                   - \cY_{1+\beta}(t_{j-1}))\Biggr\}\Biggr)
   = \EE\Biggl(\exp\Biggl\{\ii \sum_{j=1}^m
                              \theta_j
                              \sqrt{Y_{(1+\beta)/2}}
                              (B_{t_j} - B_{t_{j-1}})\Biggr\}\Biggr) \\
  &= \EE\Biggl(\exp\Biggl\{-\frac{1}{2} Y_{(1+\beta)/2}
                            \sum_{j=1}^m
                             \theta_j^2 (t_j - t_{j-1})\Biggr\}\Biggr)
   = \exp\Biggl\{-k_\beta
                  \biggl(\frac{1}{2}
                         \sum_{j=1}^m
                          \theta_j^2
                          (t_j - t_{j-1})\biggr)^{\frac{1+\beta}{2}}\Biggr\}
   = \ee^{-k_\beta \, \omega^{\frac{1+\beta}{2}}}
 \end{align*}
 as \ $N \to \infty$, \ where
 \ $\omega := \frac{1}{2} \sum_{j=1}^m \theta_j^2 (t_j - t_{j-1})$.
\ Note that, using the independence of \ $\alpha$ \ and \ $B$, \ it suffices to show
 \begin{align*}
  \Psi_N
  &:= N \Biggl[1 - \EE\Biggl(\exp\Biggl\{\ii N^{-\frac{1}{1+\beta}}
                                         \sum_{j=1}^m
                                          \theta_j
                                          (\cW_{t_j}
                                           - \cW_{t_{j-1}})\Biggr\}\Biggr)\Biggr] \\
  &= N \Biggl[1 - \EE\Biggl(\exp\Biggl\{-\frac{1}{2} N^{-\frac{2}{1+\beta}}
                                         \lambda(1+\alpha)(1-\alpha)^{-2}
                                         \sum_{j=1}^m
                                          \theta_j^2
                                          (t_j - t_{j-1})\Biggr\}\Biggr)\Biggr] \\
  &= N \int_0^1
        \Big(1 - \ee^{-\omega N^{-\frac{2}{1+\beta}}
                       \lambda(1+a)(1-a)^{-2}}\Big)
        \psi(a) (1-a)^\beta \, \dd a
   \to k_\beta \, \omega^{\frac{1+\beta}{2}}
 \end{align*}
 as \ $N \to \infty$, \ since it implies that
 \ $(1 - \Psi_N/N)^N \to \ee^{-k_\beta \, \omega^{\frac{1+\beta}{2}}}$ \ as
 \ $N \to \infty$.
\ For all \ $\vare \in (0, 1)$,
 \[
   \sup_{a\in(0,1-\vare)}
   \bigl|- N \omega
           N^{-\frac{2}{1+\beta}} (1+a)(1-a)^{-1}\bigr|
   = \omega N^{\frac{-1+\beta}{1+\beta}} (2-\vare) \vare^{-1}
   \to 0
 \]
 as \ $N \to \infty$.
\ Therefore, by Lemma \ref{ordo}, substituting \ $a = 1 - N^{-\frac{1}{1+\beta}}y$,
 \ the statement of the theorem will follow from
 \begin{equation}\label{limsup_Psin}
  \begin{aligned}
   &\limsup_{N\to\infty}
     N \int_{1-\vare}^1
        \Big|1 - \ee^{-\omega N^{-\frac{2}{1+\beta}}\lambda(1+a)(1-a)^{-2}}\Big|
       (1-a)^\beta \, \dd a \\
   &= \limsup_{N\to\infty}
       \int_0^{\vare N^{\frac{1}{1+\beta}}}
        \Big|1 - \ee^{-\omega\lambda(2-N^{-\frac{1}{1+\beta}}y)y^{-2}}\Big|
        y^\beta \, \dd y
   < \infty ,
  \end{aligned}
 \end{equation}
 and
 \begin{equation}\label{lim_Psin}
  \begin{aligned}
   &\lim_{N\to\infty}
     N \int_{1-\vare}^1
        \Big(1 - \ee^{-\omega N^{-\frac{2}{1+\beta}}\lambda(1+a)(1-a)^{-2}}\Big)
        (1-a)^\beta \, \dd a \\
   &= \lim_{N\to\infty}
       \int_0^{\vare N^{\frac{1}{1+\beta}}}
        \Big(1 - \ee^{-\omega\lambda(2-N^{-\frac{1}{1+\beta}}y)y^{-2}}\Big)
        y^\beta \, \dd y
    = \psi_1^{-1} k_\beta \, \omega^{\frac{1+\beta}{2}}
  \end{aligned}
 \end{equation}
 for all \ $\vare \in (0, 1)$.
\ Next we prove \eqref{limsup_Psin} and \eqref{lim_Psin}.

For all \ $N \in \NN$ \ and \ $\vare \in (0, 1)$, \ using \eqref{exp:4}, \ we
 have
 \begin{align*}
  \int_0^{\vare N^{\frac{1}{1+\beta}}}
   \Big|1 - \ee^{-\omega\lambda(2-N^{-\frac{1}{1+\beta}}y)y^{-2}}\Big|
   y^\beta \, \dd y
  &\leq \int_0^\infty
        \Big|1 - \ee^{-2\omega\lambda y^{-2}}\Big|
        y^\beta \, \dd y \\
  &\leq \int_0^1 y^\beta \, \dd y
        + 2 \omega \lambda \int_1^\infty y^{\beta-2} \, \dd y
   < \infty ,
 \end{align*}
 hence we obtain \eqref{limsup_Psin}.

Now we turn to prove \eqref{lim_Psin}.
For all \ $\vare \in (0, 1)$, \ we have
 \begin{align}\label{help13}
  \begin{split}
   \left|\int_{\vare N^{\frac{1}{1+\beta}}}^\infty
          \Big(1 - \ee^{-2\omega\lambda y^{-2}}\Big)
          y^\beta \, \dd y\right|
   \leq 2 \omega \lambda
        \int_{\vare N^{\frac{1}{1+\beta}}}^\infty y^{\beta-2} \, \dd y
   = \frac{2\omega\lambda}{1-\beta} (\vare N^{\frac{1}{1+\beta}})^{\beta-1}
   \to 0
   \end{split}
 \end{align}
 as \ $N \to \infty$.
\ Further, using \eqref{exp:1},
 \begin{align*}
  &\left|\int_0^{\vare N^{\frac{1}{1+\beta}}}
          \Big(1 - \ee^{-\omega\lambda(2-N^{-\frac{1}{1+\beta}}y)y^{-2}}\Big)
          y^\beta \, \dd y
          - \int_0^{\vare N^{\frac{1}{1+\beta}}}
             \Big(1 - \ee^{-2\omega\lambda y^{-2}}\Big)
             y^\beta \, \dd y\right| \\
  &\leq \int_0^{\vare N^{\frac{1}{1+\beta}}}
         \bigl|\ee^{-\omega\lambda(2-N^{-\frac{1}{1+\beta}}y)y^{-2}}
               - \ee^{-2\omega\lambda y^{-2}}\bigr|
         \, y^\beta \, \dd y \\
  &= \int_0^{\vare N^{\frac{1}{1+\beta}}}
      \ee^{-2\omega\lambda y^{-2}}
      \bigl|\ee^{\omega \lambda N^{-\frac{1}{1+\beta}a} y^{-1}} - 1\bigr|
      \, y^\beta \, \dd y \\
  &\leq \omega \lambda N^{-\frac{1}{1+\beta}}
        \int_0^{\vare N^{\frac{1}{1+\beta}}}
         \ee^{-2\omega\lambda y^{-2}}
         \ee^{\omega \lambda N^{-\frac{1}{1+\beta}}y^{-1}}
         y^{\beta-1} \, \dd y \\
  &\leq \omega \lambda N^{-\frac{1}{1+\beta}}
        \int_0^{\vare N^{\frac{1}{1+\beta}}}
         \ee^{-(2-\vare)\omega\lambda y^{-2}}
         y^{\beta-1} \, \dd y
   \leq \omega \lambda N^{-\frac{1}{1+\beta}}
        \int_0^{\vare N^{\frac{1}{1+\beta}}} y^{\beta-1} \, \dd y \\
  &= \omega \lambda N^{-\frac{1}{1+\beta}} \,
     \frac{(\vare N^{\frac{1}{1+\beta}})^\beta}{\beta}
   = \omega \lambda \, \frac{\vare^\beta N^{\frac{\beta-1}{1+\beta}}}{\beta}
   \to 0 \qquad \text{as \ $N \to \infty$,}
 \end{align*}
 hence, using \eqref{help13}, we conclude
 \begin{align*}
  &\lim_{N\to\infty}
   \int_0^{\vare N^{\frac{1}{1+\beta}}}
    \Big(1 - \ee^{-\omega\lambda(2-N^{-\frac{1}{1+\beta}}y)y^{-2}}\Big)
    y^\beta \, \dd y
   = \int_0^\infty \Big(1 - \ee^{-2\omega\lambda y^{-2}}\Big) y^\beta \, \dd y \\
  &= \frac{1}{2} \, (2\omega\lambda)^{\frac{1+\beta}{2}}
     \int_0^\infty (1 - \ee^{-u}) u^{-\frac{3+\beta}{2}} \, \dd u
   = \psi_1^{-1} k_\beta \, \omega^{\frac{1+\beta}{2}} ,
 \end{align*}
 where the last equality follows by Lemma 2.2.1 in Zolotarev \cite{Zol}, thus we
 obtain \eqref{lim_Psin}.
For the characteristic function of \ $Y_{(1+\beta)/2}$, \ see the second proof.
\proofend


\noindent{\bf Second proof of Theorem \ref{joint_aggregation_random_4}.}
By Proposition \ref{simple_aggregation_random3}, we have
 \begin{align*}
  \cD_\ff\text{-}\hspace*{-1mm}\lim_{n\to\infty} \,
   \biggl(n^{-1/2}
          \sum_{k=1}^{\lfloor nt \rfloor}
           (X^{(1)}_k - \EE(X^{(1)}_k \mid \alpha^{(1)}))\biggr)_{t\in\RR_+}
  = \frac{\sqrt{\lambda(1+\alpha)}}{1-\alpha} B ,
 \end{align*}
 where \ $(B_t)_{t\in\RR_+}$ \ is a standard Wiener process and \ $\alpha$ \ is a
 random variable having a density function of the form \eqref{alpha} with
 \ $\beta \in (-1, 1)$ \ and \ $\psi_1 \in (0, \infty)$, \ and being independent of
 \ $B$.
\ Hence it remains to prove that
 \[
   \cD_\ff\text{-}\hspace*{-1mm}\lim_{N\to\infty}
    (T^{(N)}_t)_{t\in\RR_+}
   = \bigl(\sqrt{Y_{(1+\beta)/2}} \, B_t\bigr)_{t\in\RR_+} ,
 \]
 where
 \[
   T^{(N)}_t
   := \frac{1}{N^{\frac{1}{1+\beta}}}
       \sum_{j=1}^N \frac{\sqrt{\lambda(1+\alpha^{(j)})}}{1-\alpha^{(j)}} B^{(j)}_t ,
   \qquad t \in \RR_+ , \qquad N \in \NN ,
 \]
 and \ $\alpha^{(j)}$, \ $j \in \NN$, \ and \ $B^{(j)}$, \ $j \in \NN$, \ are
 independent copies of \ $\alpha$ \ and \ $B$, \ respectively, being independent of
 each other as well.
By the continuous mapping theorem, it is enough to show that for all \ $m \in \NN$ \ and
 \ $0 =: t_0 \leq t_1 < t_2 < \cdots < t_m$,
 \begin{align*}
   \left(T^{(N)}_{t_1} - T^{(N)}_{t_0}, \ldots,
         T^{(N)}_{t_m} - T^{(N)}_{t_{m-1}}\right)
   \distr \left(\sqrt{Y_{(1+\beta)/2}} (B_{t_1} - B_{t_0}), \ldots,
                \sqrt{Y_{(1+\beta)/2}} (B_{t_m} - B_{t_{m-1}})\right)
 \end{align*}
 as \ $N \to \infty$.
\ By the portmanteau theorem, it is enough to check that for all \ $m \in \NN$,
 \ $0 = t_0 \leq t_1 < t_2 < \cdots < t_m$, \ and for all bounded and continuous
 functions \ $g : \RR^m \to \RR$,
 \begin{align*}
  &\EE\Bigl(g\bigl(T^{(N)}_{t_1} - T^{(N)}_{t_0}, \ldots,
                     T^{(N)}_{t_m} - T^{(N)}_{t_{m-1}}\bigr)\Bigr)\\
  &\qquad \to \EE\left(g\left(\sqrt{Y_{(1+\beta)/2}} (B_{t_1} - B_{t_0}), \ldots,
                          \sqrt{Y_{(1+\beta)/2}}
                          (B_{t_m} - B_{t_{m-1}})\right)\right)
 \end{align*}
 as \ $N \to \infty$.
\ Since
 \begin{align*}
  &\EE\Bigl(g\bigl(T^{(N)}_{t_1} - T^{(N)}_{t_0}, \ldots,
                     T^{(N)}_{t_m} - T^{(N)}_{t_{m-1}}\bigr)\Bigr)\\
  &=\EE\biggl[\EE\Bigl[g\bigl(T^{(N)}_{t_1} - T^{(N)}_{t_0}, \ldots,
                                T^{(N)}_{t_m} - T^{(N)}_{t_{m-1}}\bigr)
                                \,\Big|\, \alpha^{(j)}, \, j \in \NN\Bigr]\biggr] \\
  &=\EE\left[g\left(\sqrt{N^{-\frac{2}{1+\beta}}
                      \sum_{j=1}^N
                       \frac{\lambda(1+\alpha^{(j)})}{(1-\alpha^{(j)})^2}}
                       (\tB_{t_1} - \tB_{t_0}), \ldots,
                      \sqrt{N^{-\frac{2}{1+\beta}}
                      \sum_{j=1}^N
                       \frac{\lambda(1+\alpha^{(j)})}{(1-\alpha^{(j)})^2}}
                       (\tB_{t_m} - \tB_{t_{m-1}})\right)\right] \\
  &=\EE\left[h\left(N^{-\frac{2}{1+\beta}}
                      \sum_{j=1}^N
                       \frac{\lambda(1+\alpha^{(j)})}{(1-\alpha^{(j)})^2},
                       \tB_{t_1}, \ldots, \tB_{t_m}\right)\right] ,
 \end{align*}
 where \ $(\tB_t)_{t\in\RR_+}$ \ is a standard Wiener process independent of
 \ $\alpha^{(j)}$, \ $j\in\NN$, \  and \ $h : \RR^{m+1} \to \RR$ \ is an
 appropriate bounded and continuous function.
Hence it is enough to prove that
 \begin{align}\label{stable_convergence}
  \frac{1}{N^{\frac{2}{1+\beta}}}
  \sum_{j=1}^N \frac{\lambda(1+\alpha^{(j)})}{(1-\alpha^{(j)})^2}
  \distr Y_{(1+\beta)/2}
  \qquad \text{as \ $N \to \infty$,}
 \end{align}
 i.e., it suffices to show that \ $\frac{\lambda(1+\alpha)}{(1-\alpha)^2}$ \ belongs
 to the domain of normal attraction of the \ $\frac{1+\beta}{2}$-stable law of
 \ $Y_{(1+\beta)/2}$. \
Indeed, then, by the continuity theorem,
 \[
    \left(N^{-\frac{2}{1+\beta}}
                      \sum_{j=1}^N
                       \frac{\lambda(1+\alpha^{(j)})}{(1-\alpha^{(j)})^2},
                       \tB_{t_1}, \ldots, \tB_{t_m}\right)
     \distr (Y_{(1+\beta)/2}, \tB_{t_1}, \ldots, \tB_{t_m})
    \qquad \text{as \ $N\to\infty$,}
 \]
 where we additionally suppose that \ $(\widetilde B_t)_{t\in\RR_{+}}$ \ is independent
 of \ $Y_{(1+\beta)/2}$ \ as well.
Hence, using again the portmanteau theorem,
 \begin{align*}
  &\EE\left[ h\left(N^{-\frac{2}{1+\beta}}
                      \sum_{j=1}^N
                       \frac{\lambda(1+\alpha^{(j)})}{(1-\alpha^{(j)})^2},
                       \tB_{t_1}, \ldots, \tB_{t_m}\right) \right]\\
   &\qquad\to \EE\Big[ h\big(Y_{(1+\beta)/2}, \tB_{t_1}, \ldots, \tB_{t_m}\big) \Big]\\
   &\qquad \phantom{\to\,} = \EE\Big[g\Big( \sqrt{Y_{(1+\beta)/2}}(\tB_{t_1} - \tB_{t_0}),
                                     \ldots, \sqrt{Y_{(1+\beta)/2}}(\tB_{t_m} - \tB_{t_{m-1}})\Big) \Big]
 \end{align*}
 as \ $N\to\infty$, \ as desired.
Note that
 \begin{align*}
  \lim_{x\to-\infty} \vert x\vert^{\frac{1+\beta}{2}}
   \PP\biggl(\frac{\lambda(1+\alpha)}{(1-\alpha)^2} < x\biggr)
  = 0 ,
 \end{align*}
 and
 \begin{align*}
  \lim_{x\to\infty}
   x^{\frac{1+\beta}{2}}
   \PP\biggl(\frac{\lambda(1+\alpha)}{(1-\alpha)^2} > x\biggr)
  = \frac{\psi_1 (2\lambda)^{\frac{1+\beta}{2}}}{1+\beta} .
 \end{align*}
Indeed, the first convergence follows immediately due to the positivity of
 \ $\frac{\lambda(1+\alpha)}{(1-\alpha)^2}$, \ and using that
 \[
   \frac{\lambda(1+\alpha)}{(1-\alpha)^2}
   = 2 \lambda
     \left(\left(\frac{1}{1-\alpha} - \frac{1}{4}\right)^2 - \frac{1}{16}\right)
 \]
 and \ $\frac{1}{1-\alpha} \geq 1$, \ we have for all \ $x > 0$,
 \begin{align*}
  \frac{\lambda(1+\alpha)}{(1-\alpha)^2}  > x
  \qquad \Longleftrightarrow \qquad
  \alpha > 1 - \frac{1}{\frac{1}{4} + \sqrt{\frac{x}{2\lambda} + \frac{1}{16}}}
         =: 1 - \tth(\lambda,x) ,
 \end{align*}
 and hence
 \begin{align*}
  x^{\frac{1+\beta}{2}}
  &\PP\biggl(\frac{\lambda(1+\alpha)}{(1-\alpha)^2} > x\biggr)
   = x^{\frac{1+\beta}{2}} \int_{1-\tth(\lambda,x)}^1 (1-a)^\beta \psi(a) \, \dd a \\
  &= \int_0^{\sqrt{x}\tth(\lambda,x)}
      y^\beta \psi\Big(1-\frac{y}{\sqrt{x}}\Big) \, \dd y
   \to \psi_1\int_0^{\sqrt{2\lambda}} y^\beta\,\dd y
   =\frac{\psi_1(2\lambda)^{\frac{1+\beta}{2}}}{1+\beta}
  \qquad \text{as \ $x \to \infty$,}
 \end{align*}
 as desired.
Indeed, one can use the dominated convergence theorem, since there exists
 \ $\vare \in(0, 1)$ \ such that \ $|\psi(x) - \psi_1| < 2 \psi_1$ \ for all
 \ $x \in (1 - \vare, 1)$, \ and if \ $y \in (0, \sqrt{x}\tth(\lambda, x))$, \ then
 \ $1 - \frac{y}{\sqrt{x}} \in (1 - \tth(\lambda,x), 1)$ \ and hence, for large enough
 \ $x$, \ we get
 \[
   \psi\Big(1 - \frac{y}{\sqrt{x}}\Big)
   \leq 3 \psi_1 ,
   \qquad y \in (0, \sqrt{x} \tth(\lambda,x)) .
 \]
Since \ $\sqrt{x}\widetilde h(\lambda,x)\leq \sqrt{2\lambda}$, $x\in\RR_{++}$,
 \ this yields that \ $3 \psi_1 y^\beta \bbone_{[0, \sqrt{2\lambda}]}(y)$, \ $y \in \RR_+$,
 \ serves as an integrable dominating function for large enough \ $x$.
\ Consequently \eqref{stable_convergence} holds, see, e.g., Puplinskait\.{e} and
 Surgailis \cite[Remark 2.1]{PupSur2}.
Indeed, the characteristic function of the random variable \ $Y_{(1+\beta)/2}$
 \ takes the form
 \begin{align*}
  &\EE(\ee^{\ii \theta Y_{(1+\beta)/2}})\\
  & = \exp\left\{ -|\theta|^{\frac{1+\beta}{2}}
                 \frac{\Gamma(2-\frac{1+\beta}{2})}{1-\frac{1+\beta}{2}}
                 \frac{\psi_1 (2\lambda)^{\frac{1+\beta}{2}}}{1+\beta}
                 \left( \cos\left( \frac{\pi(1+\beta)}{4} \right)
                         - \ii \sign(\theta) \sin\left( \frac{\pi(1+\beta)}{4} \right) \right)
       \right\} \\
  & = \exp\left\{ -|\theta|^{\frac{1+\beta}{2}}
                  \Gamma\Big(1-\frac{1+\beta}{2}\Big)
                 \frac{\psi_1 (2\lambda)^{\frac{1+\beta}{2}}}{1+\beta}
                 \left( \cos\left(\sign(\theta) \frac{\pi(1+\beta)}{4} \right)
                         - \ii \sin\left( \sign(\theta) \frac{\pi(1+\beta)}{4} \right) \right)
       \right\} \\
  & = \exp\Big\{-k_\beta |\theta|^{\frac{1+\beta}{2}}
                 \ee^{-\ii \sign(\theta)\frac{\pi(1+\beta)}{4}}\Big\},
    \qquad \theta\in\RR.
 \end{align*}
This can be also seen using, for example, Theorem C.3 in Zolotarev \cite{Zol}
 (with the choices \ $\alpha = \frac{1+\beta}{2}$, \ $\beta = 1$ , \ $\gamma = 0$ \ and \ $\lambda=k_\beta$).
\proofend


\noindent{\bf Proof of Theorem \ref{joint_aggregation_random_5}.}
Since \ $\EE\bigl(\frac{1}{1-\alpha}\bigr) < \infty$, \ by Proposition
 \ref{simple_aggregation_random}, we have
 \[
   \frac{1}{\sqrt{N}} \tS^{(N)}
   \distrf \tcY \qquad \text{as \ $N \to \infty$,}
 \]
 where the strictly stationary Gaussian process \ $(\tcY_k)_{k\in\ZZ_+}$ \ is given
 in Proposition \ref{simple_aggregation_random}.
Consequently, by the continuous mapping theorem, for all \ $n \in \NN$, \ we get
 \[
   \cD_\ff\text{-}\hspace*{-1mm}\lim_{N\to\infty} \,
    (nN)^{-\frac{1}{2}} \tS^{(N,n)}
   = \biggl(n^{-1/2}
            \sum_{k=1}^{\lfloor nt\rfloor} \tcY_k\biggr)_{t\in\RR_+} ,
 \]
 and hence it remains to prove that
 \[
   \biggl(n^{-1/2}
          \sum_{k=1}^{\lfloor nt\rfloor} \tcY_k\biggr)_{t\in\RR_+}
   \distrf \sigma B
   \qquad \text{as \ $n \to \infty$.}
 \]
Since the processes
 \ $\bigl(n^{-1/2} \sum_{k=1}^{\lfloor nt\rfloor} \tcY_k\bigr)_{t\in\RR_+}$,
 \ $n \in \NN$, \ and \ $\sigma B$ \ are zero mean Gaussian processes, it
 suffices to show that the covariance function of
 \ $\bigl(n^{-1/2} \sum_{k=1}^{\lfloor nt\rfloor} \tcY_k\bigr)_{t\in\RR_+}$ \ converges pointwise
 to that of \ $\sigma B$ \ as \ $n \to \infty$.
\ For all \ $0 \leq t_1 \leq t_2$,
 \begin{align*}
  \cov&\biggl(n^{-1/2} \sum_{k=1}^{\lfloor nt_1\rfloor} \tcY_k,
              n^{-1/2} \sum_{k=1}^{\lfloor nt_2\rfloor} \tcY_k\biggr)
  = \frac{\lambda}{n}
    \EE\Biggl(\sum_{k=1}^{\lfloor nt_1\rfloor} \sum_{\ell=1}^{\lfloor nt_2\rfloor}
               \frac{\alpha^{|k-\ell|}}{1-\alpha}\Biggr) \\
  &\to \lambda \EE\biggl(\frac{1+\alpha}{(1-\alpha)^2}\biggr)\min(t_1, t_2)
   = \cov(\sigma B_{t_1}, \sigma B_{t_2})
 \qquad \text{as \ $n \to \infty$,}
 \end{align*}
 since one can use the decomposition \eqref{help12} together with
 \begin{align*}
  \frac{1}{n}
  \EE\biggl(\frac{\alpha(\alpha^{\lfloor nt_2\rfloor-\lfloor nt_1\rfloor}-1)
                  +(\lfloor nt_2\rfloor-\lfloor nt_1 \rfloor) (1-\alpha^2)/2}
                 {(1-\alpha)^3}\biggr)
  \to (t_2 - t_1) \EE\left(\frac{1+\alpha}{2(1-\alpha)^2}\right)
 \end{align*}
 as \ $n \to \infty$.
\ Indeed, by the dominated convergence theorem,
 \begin{align*}
  \frac{1}{n}
  \EE\biggl(\frac{\alpha(\alpha^{\lfloor nt_2 \rfloor-\lfloor nt_1 \rfloor}-1)}
                 {(1-\alpha)^3}\biggr)
  \to 0 \qquad \text{as \ $n \to \infty$,}
 \end{align*}
 where the pointwise convergence follows by
 \[
   \left|\frac{\alpha(\alpha^{\lfloor nt_2\rfloor-\lfloor nt_1\rfloor}-1)}
              {(1-\alpha)^3}\right|
   \leq \frac{1}{(1-\alpha)^3} ,
 \]
 and \ $(t_2 - t_1 + 1) \frac{\alpha}{(1-\alpha)^2}$ \ serves as an integrable
 dominating function, since, by Remark \ref{condexp},
 \ $\EE\left(\frac{\alpha}{(1-\alpha)^2}\right) < \infty$, \ and
 \begin{align*}
  \frac{1}{n}
  \left|\frac{\alpha(\alpha^{\lfloor nt_2\rfloor-\lfloor nt_1\rfloor}-1)}
             {(1-\alpha)^3}\right|
  &= \frac{\alpha(1+\alpha+\alpha^2+\cdots+\alpha^{\lfloor nt_2 \rfloor-\lfloor nt_1 \rfloor - 1})}{n(1-\alpha)^2}\\
  &\leq \frac{\alpha(\lfloor nt_2\rfloor-\lfloor nt_1\rfloor)}{n(1-\alpha)^2}
   \leq (t_2 - t_1 + 1) \frac{\alpha}{(1-\alpha)^2} .
 \end{align*}
For the second convergence, first note that, by Proposition
 \ref{simple_aggregation_random3}, we have
 \begin{align*}
  \cD_\ff\text{-}\hspace*{-1mm}\lim_{n\to\infty} \,
   \Bigg(\frac{1}{\sqrt{n}} \,
         \sum_{k=1}^{\lfloor nt\rfloor}
          (X^{(1)}_k - \EE(X^{(1)}_k \mid \alpha^{(1)}))\Bigg)_{t\in\RR_+}
  = \frac{\sqrt{\lambda(1+\alpha)}}{1-\alpha} B,
 \end{align*}
 where \ $(B_t)_{t\in\RR_+}$ \ is a standard Wiener process and \ $\alpha$ \ is a
 random variable having a density function of the form \eqref{alpha} with
 \ $\beta \in(1, \infty)$ \ and \ $\psi_1 \in (0, \infty)$, \ and being independent of
 \ $B$.
\ Hence it remains to prove that
 \[
   \cD_\ff\text{-}\hspace*{-1mm}\lim_{N\to\infty}
   \frac{1}{\sqrt{N}}
   \sum_{j=1}^N \frac{\sqrt{\lambda(1+\alpha^{(j)})}}{1-\alpha^{(j)}} B^{(j)}
   = \sigma B,
 \]
 where \ $\alpha^{(j)}$, \ $j \in \NN$, \ and \ $B^{(j)}$, \ $j\in\NN$, \ are
 independent copies of \ $\alpha$ \ and \ $B$, \ respectively, being independent of
 each other as well.
Similarly to the second proof of Theorem \ref{joint_aggregation_random_4}, it is
 enough to show that
 \begin{align*}
  \frac{1}{N} \sum_{j=1}^N\frac{\lambda(1+\alpha^{(j)})}{(1-\alpha^{(j)})^2}
  \distr \sigma^2
  \qquad \text{as \ $N \to \infty$.}
 \end{align*}
This readily follows by the strong law of large numbers, since
 \ $\EE\bigl(\frac{\lambda(1+\alpha)}{(1-\alpha)^2}\bigr) < \infty$ \ due to
 Remark \ref{condexp}.
\proofend


\noindent{\bf Proof of Theorem \ref{joint_aggregation_random_4.5}.}
We have a decomposition
 \begin{align}\label{help_decomposition}
  S^{(N,n)}_t = R^{(N,n)}_t + \tS^{(N,n)}_t ,  \qquad t \in \RR_+ ,
 \end{align}
 with
 \begin{align*}
  \begin{split}
   R^{(N,n)}_t
   := \sum_{j=1}^N \sum_{k=1}^{\lfloor nt \rfloor}
        (\EE(X^{(j)}_k\mid\alpha^{(j)}) - \EE(X^{(j)}_k))
    = \lfloor nt \rfloor
      \sum_{j=1}^N
       \left(\frac{\lambda}{1-\alpha^{(j)}}
             - \EE\left(\frac{\lambda}{1-\alpha^{(j)}}\right)\right)
  \end{split}
 \end{align*}
 for \ $t \in \RR_+$.
\ By Theorem \ref{joint_aggregation_random}, for each \ $n \in \NN$,
 \ $\cD_\ff\text{-}\hspace*{-1mm}\lim_{N\to\infty} \, N^{-\frac{1}{2}} \tS^{(N,n)}$
 \ exists, hence
 \begin{equation}\label{remainder1}
  \cD_\ff\text{-}\hspace*{-1mm}\lim_{N\to\infty} \,
  N^{-\frac{1}{1+\beta}} \tS^{(N,n)}
  = \cD_\ff\text{-}\hspace*{-1mm}\lim_{N\to\infty} \,
     N^{\frac{\beta - 1}{2(1+\beta)}} N^{-\frac{1}{2}} \tS^{(N,n)}
  = 0 .
 \end{equation}
The distribution of the random variable
 \ $\lambda (1 - \alpha)^{-1} - \EE(\lambda (1 - \alpha)^{-1})$ \ belongs to the
 domain of attraction of an \ $(1 + \beta)$-stable distribution.
Indeed, we have
 \begin{align*}
  &\lim_{x\to\infty}
    x^{1+\beta}
    \PP\left(\frac{\lambda}{1-\alpha} - \EE\left(\frac{\lambda}{1-\alpha}\right)
             > x\right) \\
  &= \lim_{x\to\infty}
      x^{1+\beta}
      \PP\left(\alpha
               > 1- \frac{1}{\lambda^{-1} x + \EE((1 - \alpha)^{-1})}\right) \\
  &= \lim_{x\to\infty}
      \frac{1}{x^{-(1+\beta)}}
      \int_{1-(\lambda^{-1}x+\EE((1-\alpha)^{-1}))^{-1}}^1
       \psi(a) (1 - a)^\beta \, \dd a \\
  &= \lim_{x\to\infty}
      \frac{-\psi(1-(\lambda^{-1}x+\EE((1-\alpha)^{-1}))^{-1})
             (\lambda^{-1}x+\EE((1-\alpha)^{-1}))^{-\beta-2} \lambda^{-1}}
           {-(1+\beta)x^{-(1+\beta)-1}}
   = \frac{\psi_1 \lambda^{1+\beta}}{1+\beta}
 \end{align*}
 by L'H\^opital's rule.
Further, using that \ $\PP(\lambda(1-\alpha)^{-1} > 0) = 1$,
 \begin{align*}
  &\lim_{x\to-\infty}
    |x|^{1+\beta}
    \PP\left(\frac{\lambda}{1-\alpha} - \EE\left(\frac{\lambda}{1-\alpha}\right)
             \leq x\right)
  = \lim_{x\to-\infty} |x|^{1+\beta}\cdot 0 = 0 .
 \end{align*}
Consequently, for each \ $n \in \NN$,
 \[
   \cD_\ff\text{-}\hspace*{-1mm}\lim_{N\to\infty} \,
    N^{-\frac{1}{1+\beta}} R^{(N,n)}
   = \bigl(\lfloor nt \rfloor Z_{1+\beta}\bigr)_{t\in\RR_+} ,
 \]
 see, e.g., Puplinskait\.{e} and Surgailis \cite[Remark 2.1]{PupSur2}.
Indeed, the characteristic function of the random variable \ $Z_{1+\beta}$ \ takes
 the form
 \begin{align*}
  &\EE(\ee^{\ii \theta Z_{1+\beta}})\\
  & = \exp\left\{ -|\theta|^{1+\beta}
                 \frac{\Gamma(2-(1+\beta))}{1-(1+\beta)}
                 \frac{\psi_1\lambda^{1+\beta}}{1+\beta}
                 \left( \cos\left( \frac{\pi(1+\beta)}{2} \right)
                         - \ii \sign(\theta) \sin\left( \frac{\pi(1+\beta)}{2} \right) \right)
       \right\} \\
  & = \exp\Big\{-|\theta|^{1+\beta} \frac{\Gamma(1-\beta)}{-\beta}
                  \frac{\psi_1\lambda^{1+\beta}}{1+\beta}
                  \ee^{-\ii \sign(\theta)\frac{\pi(1+\beta)}{2}}\Big\}\\
  & = \exp\Big\{ -|\theta|^{1+\beta} \omega_\beta(\theta)\Big\},
    \qquad \theta\in\RR.
 \end{align*}
Together with \eqref{remainder1}, we obtain the first convergence.

By Theorem \ref{joint_aggregation_random_4}, for each \ $N \in \NN$,
 \ $\cD_\ff\text{-}\hspace*{-1mm}\lim_{n\to\infty} \, n^{-\frac{1}{2}} \tS^{(N,n)}$ \ exists and hence
 \begin{align*}
  \cD_\ff\text{-}\hspace*{-1mm}\lim_{n\to\infty} \, n^{-1} \tS^{(N,n)}
  = \cD_\ff\text{-}\hspace*{-1mm}\lim_{n\to\infty} \,
     n^{-\frac{1}{2}} n^{-\frac{1}{2}} \tS^{(N,n)}
  = 0 ,
 \end{align*}
 and
 \begin{align*}
  \cD_\ff\text{-}\hspace*{-1mm}\lim_{n\to\infty} \, n^{-1} R^{(N,n)}
  =\left(t
         \sum_{j=1}^N
          \left(\frac{\lambda}{1-\alpha^{(j)}}
                -\EE\left(\frac{\lambda}{1-\alpha^{(j)}}\right)\right)\right)_{t\in\RR_+} .
 \end{align*}
Based on the above considerations, using the decomposition
 \eqref{help_decomposition} as well, we obtain the second convergence.
\proofend


\noindent{\bf Proof of Theorem \ref{joint_aggregation_random_5.5}.}
First note that, since \ $\beta>1$, \ by Remark \ref{condexp},
 \ $\var((1-\alpha)^{-1})<\infty$.
\ Hence, by the central limit theorem, for each \ $n \in \NN$, \
 \[
   \cD_\ff\text{-}\hspace*{-1mm}\lim_{N\to\infty} N^{-\frac{1}{2}}R^{(N,n)}
   =\left(\lfloor nt\rfloor W_{\lambda^2\var((1-\alpha)^{-1})} \right)_{t\in\RR_+} .
 \]
Consequently,
 \[
   \cD_\ff\text{-}\hspace*{-1mm}\lim_{n\to\infty} \,
   \cD_\ff\text{-}\hspace*{-1mm}\lim_{N\to\infty} \,
    n^{-1} N^{-\frac{1}{2}} R^{(N,n)}
   =  (W_{\lambda^2\var((1-\alpha)^{-1})} t)_{t\in\RR_+} .
 \]
By Theorem \ref{joint_aggregation_random_5},
 \ $\cD_\ff\text{-}\hspace*{-1mm}\lim_{n\to\infty} \,
    \cD_\ff\text{-}\hspace*{-1mm}\lim_{N\to\infty} \,
     (nN)^{-\frac{1}{2}} \tS^{(N,n)}$
 \ exists, hence
 \begin{align*}
  \cD_\ff\text{-}\hspace*{-1mm}\lim_{n\to\infty} \,
  \cD_\ff\text{-}\hspace*{-1mm}\lim_{N\to\infty} \,
   n^{-1} N^{-\frac{1}{2}} \tS^{(N,n)}
  = \cD_\ff\text{-}\hspace*{-1mm}\lim_{n\to\infty} \,
    \cD_\ff\text{-}\hspace*{-1mm}\lim_{N\to\infty} \,
     n^{-\frac{1}{2}} (nN)^{-\frac{1}{2}} \tS^{(N,n)}
  = 0 .
 \end{align*}
Using the decomposition \eqref{help_decomposition}, we have the first convergence.

Similarly, for each \ $N \in \NN$,
 \[
   \cD_\ff\text{-}\hspace*{-1mm}\lim_{n\to\infty} n^{-1}R^{(N,n)}
   = \left( \sum_{j=1}^N \left(\frac{\lambda}{1-\alpha^{(j)}}
                      -\EE\left(\frac{\lambda}{1-\alpha^{(j)}}\right)\right) t \right)_{t\in\RR_+} ,
 \]
 and, by the central limit theorem,
 \[
   \cD_\ff\text{-}\hspace*{-1mm}\lim_{N\to\infty} \,
   \cD_\ff\text{-}\hspace*{-1mm}\lim_{n\to\infty} \,
    n^{-1} N^{-\frac{1}{2}} R^{(N,n)}
   =  (W_{\lambda^2\var((1-\alpha)^{-1})} t)_{t\in\RR_+} .
 \]
By Theorem \ref{joint_aggregation_random_5}, we also have
 \[
   \cD_\ff\text{-}\hspace*{-1mm}\lim_{N\to\infty} \,
   \cD_\ff\text{-}\hspace*{-1mm}\lim_{n\to\infty} \,
   n^{-1} N^{-\frac{1}{2}} \tS^{(N,n)}
   = 0 ,
 \]
 which yields the second convergence using the decomposition
 \eqref{help_decomposition} as well.
\proofend

\vspace*{5mm}

\appendix

\vspace*{5mm}

\noindent{\bf\Large Appendices}

\section{Non-Markov property of the randomized INAR(1) model}
\label{Markov}

The aim of this appendix is to show that the randomized INAR(1) process \ $(X_k)_{k\in\ZZ_+}$ \
 defined in Section \ref{RINAR} does not have the Markov property provided that \ $\alpha$ \ is non-degenerate.
We show that if \ $\alpha$ \ is non-degenerate, then
 \[
   \PP(X_2 = 0 \mid X_1 = 1, X_0 = 0) \ne \PP(X_2 = 0 \mid X_1 = 1),
 \]
 implying our statement.
By the strict stationarity of \ $(X_k)_{k\in\ZZ_+}$, \ the conditional independence
 of \ $\xi_{1,1}$, \ $\vare_1$ \ and \ $X_0$ \ given \ $\alpha$, \ and \eqref{help_stac_RINAR1_1}--\eqref{help_stac_RINAR1_3},
 we have
 \begin{align*}
  \PP(X_2 = 0 \mid X_1 = 1)
  & = \PP(X_1 = 0 \mid X_0 = 1)
    = \frac{\PP(X_1 = 0, X_0 = 1)}{\PP(X_0 = 1)}
    = \frac{\PP(\xi_{1,1} = 0, \vare_1 = 0, X_0 = 1)}{\PP(X_0 = 1)} \\
  & = \frac{\int_0^1
           \PP(\xi_{1,1} = 0, \vare_1 = 0, X_0 = 1 \mid \alpha = a) \, \PP_\alpha(\dd a)}
         {\int_0^1 \PP(X_0 = 1 \mid \alpha = a) \, \PP_\alpha(\dd a)}\\
  & = \frac{\int_0^1
           \PP(\xi_{1,1} = 0 \mid \alpha = a) \PP(\vare_1 = 0 \mid \alpha = a)
           \PP(X_0 = 1 \mid \alpha = a) \, \PP_\alpha(\dd a)}
         {\int_0^1 \PP(X_0 = 1 \mid \alpha = a) \, \PP_\alpha(\dd a)} \\
  & = \frac{\int_0^1 (1-a) \ee^{-\lambda} \frac{\lambda}{1-a} \, \ee^{-\frac{\lambda}{1-a}} \, \PP_\alpha(\dd a)}
         {\int_0^1 \frac{\lambda}{1-a} \, \ee^{-\frac{\lambda}{1-a}} \, \PP_\alpha(\dd a)}
    = \frac{\int_0^1 \ee^{ -\lambda -\frac{\lambda}{1-a} } \, \PP_\alpha(\dd a)}
         {\int_0^1 \frac{1}{1-a} \, \ee^{-\frac{\lambda}{1-a}} \, \PP_\alpha(\dd a)}.
 \end{align*}
Similarly, we have
 \begin{align*}
  \PP(X_2 = 0 \mid X_1 = 1, X_0=0)
  & = \frac{\PP(X_2 = 0, X_1 = 1, X_0 = 0)}{\PP(X_1 = 1, X_0=0)}
    = \frac{\PP(\xi_{2,1} = 0, \vare_2 = 0, \vare_1=1,X_0 = 0)}{\PP(\vare_1=1,X_0 = 0)} \\
  & = \frac{\int_0^1 (1-a) \ee^{-\lambda} \lambda \ee^{-\lambda} \, \ee^{-\frac{\lambda}{1-a}} \, \PP_\alpha(\dd a)}
         {\int_0^1  \lambda \ee^{-\lambda} \, \ee^{-\frac{\lambda}{1-a}} \, \PP_\alpha(\dd a)}
    = \frac{\int_0^1 (1-a) \ee^{ -\lambda -\frac{\lambda}{1-a} } \, \PP_\alpha(\dd a)}
         {\int_0^1 \ee^{-\frac{\lambda}{1-a}} \, \PP_\alpha(\dd a)}.
 \end{align*}
By Cauchy--Schwarz's inequality, we have
 \[
   \left( \int_0^1  \ee^{-\frac{\lambda}{1-a}} \, \PP_\alpha(\dd a) \right)^2
      \leq \int_0^1  (1-a) \ee^{-\frac{\lambda}{1-a}} \, \PP_\alpha(\dd a)
           \int_0^1  \frac{1}{1-a} \ee^{-\frac{\lambda}{1-a}} \, \PP_\alpha(\dd a),
 \]
 and equality holds if and only if there exists some positive constant \ $C>0$ \ such that
 \ $(1-a) \ee^{-\frac{\lambda}{1-a}} = C \frac{1}{1-a} \ee^{-\frac{\lambda}{1-a}}$ \
 $\PP_\alpha$-almost every \ $a\in(0,1)$, \ which is equivalent to the fact that there exists
 \ $C\in(0,1)$ \ such that \ $\PP_\alpha$ \ is the Dirac measure concentrated on the point
 \ $1-\sqrt{C}$.
\ Consequently, \ $\PP(X_2 = 0 \mid X_1 = 1, X_0=0)\geq   \PP(X_2 = 0 \mid X_1 = 1)$ \ and equality holds
 if and only if \ $\PP_\alpha$ \ is a Dirac measure concentrated on some point in \ $(0,1)$,
 \ i.e., \ $\alpha$ \ is degenerate.
Hence if \ $\alpha$ \ is non-degenerate, then the randomized INAR(1) process \ $(X_k)_{k\in\ZZ_+}$ \ does not have
 the Markov property.
If \ $\alpha$ \ is degenerate, then \ $(X_k)_{k\in\ZZ_+}$ \ is a usual INAR(1) model being a Markov chain.

\section{Approximations of the exponential function and some of its integrals}
\label{App1}


In this appendix we collect some useful approximations of the exponential function and
 some of its integrals.

We will frequently use the following the well-known inequalities:
 \begin{gather}\label{exp:4}
  1 - \ee^{-x} \leq x, \qquad x \in \RR, \\  \label{exp:3}
  |\ee^{\ii u} - 1| \leq |u|, \qquad |\ee^{\ii u} - 1 - \ii u| \leq u^2/2 ,
  \qquad u \in \RR .
 \end{gather}

The next lemma is about how the inequalities in \eqref{exp:3} change
 if we replace \ $u\in\RR$ \ by an arbitrary complex number.

\begin{Lem}\label{exp}
We have
 \begin{gather}\label{exp:1}
  |\ee^z - 1| \leq |z| \ee^{|z|} , \qquad z \in \CC , \\
  \label{exp:2}
  |\ee^z - 1 - z| \leq \frac{|z|^2}{2}\ee^{|z|} , \qquad z \in \CC .
 \end{gather}
\end{Lem}

\noindent{\bf Proof.}
For any \ $z \in \CC$ \ we have
 \begin{equation*}
  \begin{split}
   |\ee^z - 1|
   &= \left|z + \frac{z^2}{2!} + \frac{z^3}{3!} + \dots\right|
    \leq |z| \left(1 + \frac{|z|}{2!} + \frac{|z|^2}{3!} + \dots\right) \\
   &\leq |z| \left(1 + \frac{|z|}{1!} + \frac{|z|^2}{2!} + \dots\right)
    = |z| \ee^{|z|} ,
  \end{split}
 \end{equation*}
 \begin{equation*}
  \begin{split}
   |\ee^z - 1 - z| &
   = \left|\frac{z^2}{2!} + \frac{z^3}{3!} + \dots\right|
   \leq \frac{|z|^2}{2}
        \left(1 + \frac{|z|}{3} + \frac{|z|^2}{3\cdot4} + \dots\right) \\
   &\leq \frac{|z|^2}{2} \left(1 + \frac{|z|}{1!} + \frac{|z|^2}{2!} + \dots\right)
    = \frac{|z|^2}{2} \ee^{|z|} ,
  \end{split}
 \end{equation*}
 since \ $3\cdot4\cdots(n+2) \geq n!$ \ for any \ $n \in \NN$.
\proofend


\begin{Lem}\label{ordo}
Suppose that \ $(0, 1) \ni x \mapsto \psi(x) (1 - x)^\beta$ \ is a probability density,
 where \ $\psi$ \ is a function on \ $(0, 1)$ \ having a limit
 \ $\lim_{x\uparrow 1} \psi(x) = \psi_1 \in (0, \infty)$ \ (and necessarily \ $\beta\in(-1,\infty)$).
\ For all \ $a \in (0, 1)$, \ let \ $(z_N(a))_{N\in\NN}$ \ be a sequence of complex
 numbers such that
 \begin{gather}\label{help_ordo1}
  \lim_{N\to\infty} \sup_{a\in(0,1-\vare)} |Nz_N(a)| = 0 \qquad
  \text{for all \ $\vare \in (0, 1)$,} \\ \nonumber
  \limsup_{N\to\infty}
   N \int_{1-\vare_0}^1
      \left|1 - \ee^{\frac{\lambda}{1-a}z_N(a)}\right| (1-a)^\beta \, \dd a
   < \infty \qquad \text{for some \ $\vare_0 \in (0, 1)$,} \\ \nonumber
  \lim_{\vare\downarrow0}
   \limsup_{N\to\infty}
    \left|N \int_{1-\vare}^1
             \left(1 - \ee^{\frac{\lambda}{1-a}z_N(a)}\right) (1-a)^\beta \, \dd a
          - I\right|
  = 0
 \end{gather}
 with some \ $I \in \CC$.
\ Then
 \[
   \lim_{N\to\infty}
    N \int_0^1
       \left(1 - \ee^{\frac{\lambda}{1-a}z_N(a)}\right) \psi(a) (1-a)^\beta
       \, \dd a
   = \psi_1 I .
 \]
\end{Lem}

\noindent{\bf Proof.}
Using dominated convergence theorem, first we check that
  \begin{align}\label{help17}
   \lim_{N\to\infty}
    N \int_0^{1-\vare}
       \left(1 - \ee^{\frac{\lambda}{1-a}z_N(a)}\right) \psi(a) (1 - a)^\beta
       \, \dd a
   = 0
   \qquad \text{for all \ $\vare\in(0,1)$.}
 \end{align}
By applying \eqref{exp:1} and using \eqref{help_ordo1}, for any \ $\vare \in (0, 1)$
 \ and \ $a \in (0, 1 - \vare)$, \ we get
 \begin{equation}\label{ordo1}
  \Bigl|N \left(1 - \ee^{\frac{\lambda}{1-a}z_N(a)}\right)\Bigr|
  \leq N\Bigl|\frac{\lambda}{1-a}z_N(a)\Bigr|
       \ee^{\left|\frac{\lambda}{1-a}z_N(a)\right|}
  \to 0
 \end{equation}
 as \ $N \to \infty$.
\ Further, if \ $\vare \in (0, 1)$ \ and \ $a \in (0, 1 - \vare)$, \ then
 \[
   \left|N\left(1 - \ee^{\frac{\lambda}{1-a}z_N(a)}\right)\right|
   \leq \frac{\lambda}{\vare}
        \sup_{N\in\NN} \sup_{a\in (0,1-\vare)} |N z_N(a)| \,
             \ee^{\frac{\lambda}{\vare}
                  \sup_{N\in\NN}\sup_{a\in (0,1-\vare)} |z_N(a)|} =: C_\vare ,
 \]
 where \ $C_\vare \in \RR_+$.
\ Since \ $\int_0^1 \psi(a) (1-a)^{\beta} \, \dd a = 1$, \ we have
 \[
   \left| N\int_0^{1-\vare}
            \left(1 - \ee^{\frac{\lambda}{1-a}z_N(a)}\right) \psi(a)(1-a)^\beta
            \, \dd a\right|
   \leq \int_0^{1-\vare} C_\vare \psi(a) (1-a)^\beta \, \dd a
   < \infty .
 \]
Therefore, \ $(0, 1 - \vare) \ni a \mapsto C_\vare \psi(a)(1-a)^\beta$ \ serves as a
 dominating integrable function.
Thus the pointwise convergence in \eqref{ordo1} results \eqref{help17}.
Moreover, for all \ $\vare \in (0, 1)$, \ we have
 \begin{align*}
  &\left|N \int_0^1
           \left(1 - \ee^{\frac{\lambda}{1-a}z_N(a)}\right) \psi(a) (1-a)^\beta
           \, \dd a
        - \psi_1 I\right| \\
  &\qquad\qquad
   \leq \left|N \int_0^{1-\vare}
                 \left(1 - \ee^{\frac{\lambda}{1-a}z_N(a)}\right) \psi(a) (1-a)^\beta
                 \, \dd a\right| \\
  &\qquad\qquad\quad
        + \left|N \int_{1-\vare}^1
                   \left(1 - \ee^{\frac{\lambda}{1-a}z_N(a)}\right)
                   (\psi(a) - \psi_1) (1-a)^\beta
                   \, \dd a\right| \\
  &\qquad\qquad\quad
        + \psi_1 \left|N \int_{1-\vare}^1
                          \left(1 - \ee^{\frac{\lambda}{1-a}z_N(a)}\right)
                          (1-a)^\beta
                          \, \dd a
                       - I\right| ,
 \end{align*}
 where
 \begin{align*}
  &\left|N \int_{1-\vare}^1
            \left(1 - \ee^{\frac{\lambda}{1-a}z_N(a)}\right) (\psi(a) - \psi_1)
            (1-a)^\beta
            \, \dd a\right| \\
  &\leq N \sup_{a\in[1-\vare,1)} |\psi(a) - \psi_1|
        \int_{1-\vare}^1
         \left|1 - \ee^{\frac{\lambda}{1-a}z_N(a)}\right| (1-a)^\beta
         \, \dd a ,
 \end{align*}
 with \ $\sup_{a\in[1-\vare,1)} |\psi(a) - \psi_1| \to 0$ \ as
 \ $\vare \downarrow 0$, \ by the assumption.
First taking \ $\limsup_{N\to\infty}$ \ and then \ $\vare \downarrow 0$, \ using \eqref{help17},
 we obtain the statement.
\proofend


\section{A representation of fractional Brownian motion due to Pilipauskait\.{e} and Surgailis \cite{PilSur}}
\label{App2}

We recall an integral representation of the fractional Brownian motion with
 Hurst parameter in \ $\bigl(\frac{1}{2}, 1\bigr)$ \ due to Pilipauskait\.{e} and Surgailis
 \cite{PilSur} in order to connect our results with the ones in Pilipauskait\.{e}
 and Surgailis \cite{PilSur} and in Puplinskait\.{e} and Surgailis \cite{PupSur1},
 \cite{PupSur2}.

For all \ $\beta \in (0, 1)$ \ let us consider the stochastic process given by
 \begin{equation}\label{fractionaldef}
  \tcB_{1-\frac{\beta}{2}}(t)
  := \int_{\RR_+\times\RR} (f(x,t-s) - f(x,-s)) \, Z(\dd x, \dd s) ,
  \qquad t\in\RR_+,
 \end{equation}
 where
 \begin{equation}\label{fractionaldef2}
  f(x, t) := \begin{cases}
              (1 - \ee^{-xt})/x & \text{if \ $x\in\RR_{++}$ \ and \ $t\in\RR_{++},$} \\
              0 & \text{otherwise,}
             \end{cases}
 \end{equation}
 with respect to a Gaussian random measure \ $Z(\dd x, \dd s)$ \ on
 \ $\RR_+ \times \RR$ \ with zero mean, variance
 \ $\nu(\dd x, \dd s) := (2-\beta)(1-\beta)/\Gamma(\beta) x^\beta \, \dd x \, \dd s$
 \ and characteristic function
 \ $\EE(\ee^{\ii \theta Z(A)}) = \ee^{-\theta^2\nu(A)/2}$ \ for each Borel set
 \ $A \subset \RR_+ \times \RR$ \ with \ $\nu(A) < \infty$ \ and \ $\theta\in\RR$.
\ Note that, by Pilipauskait\.{e} and Surgailis \cite[page 1014]{PilSur},
 \ $(\tcB_{1-\frac{\beta}{2}}(t))_{t\in\RR_+}$ \ is a fractional Brownian motion
 multiplied by some constant.
In what follows we check that this constant is in fact one.
It suffices to show that the variance of the process defined in
 \eqref{fractionaldef} at time \ $1$ \ is \ $1$.
\ By \eqref{fractionaldef} and \eqref{fractionaldef2} (see also formula (2.4) in
 Pilipauskait\.{e} and Surgailis \cite{PilSur}) one can easily see that the
 variance of \ $\tcB_{1-\frac{\beta}{2}}(1)$ \ takes the form
 \begin{align*}
  \EE(\tcB_{1-\frac{\beta}{2}}(1)^2)
  = \frac{(2-\beta)(1-\beta)}{\Gamma(\beta)}
  \int_0^\infty \int_{-\infty}^\infty
   \left(f(x,1-s)-f(x,-s) \right)^2 x^\beta \, \dd s \, \dd x ,
 \end{align*}
 where, for \ $x\in\RR_{++}$ \ and \ $t \in \RR_+$,
 \begin{align*}
  &\int_{-\infty}^\infty \big( f(x,t-s)-f(x,-s) \big)^2 x^{\beta}\,\dd s\\
  &\quad
   = \int_{-\infty}^0
      \left(\frac{1-\ee^{-x(t-s)}}{x} - \frac{1-\ee^{-x(-s)}}{x} \right)^2
      x^{\beta} \, \dd s
     + \int_0^t \left(\frac{1-\ee^{-x(t-s)}}{x}\right)^2 x^{\beta} \, \dd s \\
  &\quad
   = \int_{-\infty}^0 \ee^{2xs}(1 - \ee^{-xt})^2 x^{\beta - 2} \, \dd s
     + \int_0^t (1-\ee^{-x(t-s)})^2 x^{\beta - 2} \, \dd s \\
  &\quad
   = \frac{1}{2}(1 - \ee^{-xt})^2 x^{\beta-3}
     + \int_0^t (1-\ee^{-x(t-s)})^2 x^{\beta-2} \, \dd s .
 \end{align*}
Hence, with repeated partial integration, we have
 \begin{align*}
  \EE(\tcB_{1-\frac{\beta}{2}}(1)^2)
  &= \frac{(2-\beta)(1-\beta)}{\Gamma(\beta)}
     \int_0^\infty
      \left[\frac{1}{2} (1-\ee^{-x})^2 x^{\beta-3}
            + \int_0^1 (1-\ee^{-x(1-s)})^2 x^{\beta-2} \, \dd s\right] \dd x \\
  &= \frac{(2-\beta)(1-\beta)}{\Gamma(\beta)}
     \int_0^\infty
      \left(\ee^{-x} - 1 + x\right) x^{\beta-3} \, \dd x \\
  & = \frac{(2-\beta)(1-\beta)}{\Gamma(\beta)}
           \frac{1}{2-\beta}
           \int_0^\infty
          \left(-\ee^{-x} + 1 \right) x^{\beta-2} \, \dd x \\
  &= \frac{(2-\beta)(1-\beta)}{\Gamma(\beta)} \cdot
     \frac{\Gamma(\beta)}{(2-\beta)(1-\beta)}
   = 1 ,
 \end{align*}
 as desired.

\section*{Acknowledgements}
We are grateful to the referee for several valuable comments and suggestions, especially for initiating
the centralization by the empirical mean as well.

\bibliographystyle{plain}
\bibliography{aggr6}

\def\polhk#1{\setbox0=\hbox{#1}{\ooalign{\hidewidth
  \lower1.5ex\hbox{`}\hidewidth\crcr\unhbox0}}}
\begin{thebibliography}{10}

\bibitem{AloAlz}
M.A. Al-Osh and A.A. Alzaid.
\newblock First-order integer-valued autoregressive ({INAR}({$1$})) process.
\newblock {\em J. Time Ser.\ Anal.}, 8(3):261--275, 1987.

\bibitem{BarIspPap}
M.~Barczy, M.~Isp{\'a}ny, and G.~Pap.
\newblock Asymptotic behavior of conditional least squares estimators for
  unstable integer-valued autoregressive models of order 2.
\newblock {\em Scand.\ J. Stat.}, 41(4):866--892, 2014.

\bibitem{BerFenGhoKul}
J.~Beran, Y.~Feng, S.~Ghosh, and R.~Kulik.
\newblock {\em Long-{M}emory {P}rocesses. {P}robabilistic {P}roperties and
  {S}tatistical {M}ethods}.
\newblock Springer, Heidelberg, 2013.

\bibitem{BerSchGho}
J.~Beran, M.~Sch{\"u}tzner, and S.~Ghosh.
\newblock From short to long memory: aggregation and estimation.
\newblock {\em Comput.\ Statist.\ Data Anal.}, 54(11):2432--2442, 2010.

\bibitem{CelLeiPhi}
D.~Celov, R.~Leipus, and A.~Philippe.
\newblock Time series aggregation, disaggregation, and long memory.
\newblock {\em Liet.\ Mat.\ Rink.}, 47(4):466--481, 2007.

\bibitem{DomKaj}
C.~Dombry and I.~Kaj.
\newblock The on-off network traffic model under intermediate scaling.
\newblock {\em Queueing Syst.}, 69(1):29--44, 2011.

\bibitem{FosWil}
J.H. Foster and J.A. Williamson.
\newblock Limit theorems for the {G}alton--{W}atson process with time-dependent
  immigration.
\newblock {\em Z. Wahrscheinlichkeitstheorie und Verw.\ Gebiete}, 20:227--235,
  1971.

\bibitem{GaiKaj}
R.~Gaigalas and I.~Kaj.
\newblock Convergence of scaled renewal processes and a packet arrival model.
\newblock {\em Bernoulli}, 9(4):671--703, 2003.

\bibitem{GonGou}
E.~Gon{\c{c}}alves and Ch. Gouri{\'e}roux.
\newblock Agr\'egation de processus autor\'egressifs d'ordre {$1$}.
\newblock {\em Ann.\ \'Econom.\ Statist.}, (12):127--149, 1988.

\bibitem{Gra}
C.W.J. Granger.
\newblock Long memory relationships and the aggregation of dynamic models.
\newblock {\em J. Econometrics}, 14(2):227--238, 1980.

\bibitem{IglTer}
E.~Igl{\'o}i and G.~Terdik.
\newblock Long-range dependence through gamma-mixed {O}rnstein--{U}hlenbeck
  process.
\newblock {\em Electron.\ J. Probab.}, 4:no.\ 16, 33 pp.\ (electronic), 1999.

\bibitem{JacShi}
J.~Jacod and A.N. Shiryaev.
\newblock {\em Limit {T}heorems for {S}tochastic {P}rocesses}, volume 288 of
  {\em Grundlehren der Mathematischen Wissenschaften [Fundamental Principles of
  Mathematical Sciences]}.
\newblock Springer-Verlag, Berlin, second edition, 2003.

\bibitem{Jir}
M.~Jirak.
\newblock Limit theorems for aggregated linear processes.
\newblock {\em Adv.\ in Appl.\ Probab.}, 45(2):520--544, 2013.

\bibitem{Johnson}
Norman~L. Johnson, Samuel Kotz, and N.~Balakrishnan.
\newblock {\em Discrete multivariate distributions}.
\newblock Wiley Series in Probability and Statistics: Applied Probability and
  Statistics. John Wiley \& Sons, Inc., New York, 1997.
\newblock A Wiley-Interscience Publication.

\bibitem{Key}
E.S. Key.
\newblock Limiting distributions and regeneration times for multitype branching
  processes with immigration in a random environment.
\newblock {\em Ann.\ Probab.}, 15(1):344--353, 1987.

\bibitem{LeoSavZhi}
N.N. Leonenko, V.~Savani, and A.A. Zhigljavsky.
\newblock Autoregressive negative binomial processes.
\newblock {\em Ann.\ I.S.U.P.}, 51(1-2):25--47, 2007.

\bibitem{Li}
Z.~Li.
\newblock {\em Measure-valued branching {M}arkov processes}.
\newblock Probability and its {A}pplications ({N}ew {Y}ork). Springer,
  Heidelberg, 2011.

\bibitem{McK}
E.~McKenzie.
\newblock Some simple models for discrete variate time series.
\newblock {\em JAWRA Journal of the American Water Resources Association},
  21(4):645--650, 1985.

\bibitem{MikResRooSte}
T.~Mikosch, S.~Resnick, H.~Rootz{\'e}n, and A.~Stegeman.
\newblock Is network traffic approximated by stable {L}\'evy motion or
  fractional {B}rownian motion?
\newblock {\em Ann.\ Appl.\ Probab.}, 12(1):23--68, 2002.

\bibitem{Ned}
F.~Ned{\'e}nyi.
\newblock Conditional least squares estimators for multitype {G}alton--{W}atson
  processes.
\newblock {\em Acta Sci.\ Math.\ (Szeged)}, 81(1-2):325--348, 2015.

\bibitem{NedPap}
F.~Ned{\'e}nyi and G.~Pap.
\newblock Iterated scaling limits for aggregation of random coefficient {AR}(1)
  and {INAR}(1) processes.
\newblock {\em Statist. Probab. Lett.}, 118:16--23, 2016.

\bibitem{OppVia2004}
G.~Oppenheim and M.-C. Viano.
\newblock Aggregation of random parameters {O}rnstein-{U}hlenbeck or {AR}
  processes: some convergence results.
\newblock {\em J. Time Ser.\ Anal.}, 25(3):335--350, 2004.

\bibitem{PilSur}
V.~Pilipauskait{\.e} and D.~Surgailis.
\newblock Joint temporal and contemporaneous aggregation of random-coefficient
  {AR}(1) processes.
\newblock {\em Stochastic Process.\ Appl.}, 124(2):1011--1035, 2014.

\bibitem{PilSur2}
V.~Pilipauskait{\.e} and D.~Surgailis.
\newblock Joint aggregation of random-coefficient {${\rm AR}(1)$} processes
  with common innovations.
\newblock {\em Statist.\ Probab.\ Lett.}, 101:73--82, 2015.

\bibitem{PipTaqLev}
V.~Pipiras, M.S. Taqqu, and J.B. Levy.
\newblock Slow, fast and arbitrary growth conditions for renewal-reward
  processes when both the renewals and the rewards are heavy-tailed.
\newblock {\em Bernoulli}, 10(1):121--163, 2004.

\bibitem{PupSur1}
D.~Puplinskait{\.e} and D.~Surgailis.
\newblock Aggregation of random-coefficient {AR}(1) process with infinite
  variance and common innovations.
\newblock {\em Lith.\ Math.\ J.}, 49(4):446--463, 2009.

\bibitem{PupSur2}
D.~Puplinskait{\.e} and D.~Surgailis.
\newblock Aggregation of a random-coefficient {${\rm AR}(1)$} process with
  infinite variance and idiosyncratic innovations.
\newblock {\em Adv.\ in Appl.\ Probab.}, 42(2):509--527, 2010.

\bibitem{Rob}
P.M. Robinson.
\newblock Statistical inference for a random coefficient autoregressive model.
\newblock {\em Scand.\ J. Statist.}, 5(3):163--168, 1978.

\bibitem{Sen2}
E.~Seneta.
\newblock Functional equations and the {G}alton-{W}atson process.
\newblock {\em Adv.\ in Appl.\ Probab.}, 1:1--42, 1969.

\bibitem{Shi}
A.N. Shiryaev.
\newblock {\em Probability}, volume~95 of {\em Graduate Texts in Mathematics}.
\newblock Springer-Verlag, New York, second edition, 1996.
\newblock Translated from the first (1980) Russian edition by R.P. Boas.

\bibitem{SteHar}
F.W. Steutel and K.~van Harn.
\newblock Discrete analogues of self-decomposability and stability.
\newblock {\em Ann.\ Probab.}, 7(5):893--899, 1979.

\bibitem{Str}
D.W. Stroock.
\newblock {\em Probability {T}heory. An {A}nalytic {V}iew}.
\newblock Cambridge University Press, Cambridge, second edition, 2011.

\bibitem{TaqWilShe}
M.S. Taqqu, W.~Willinger, and R.~Sherman.
\newblock Proof of a fundamental result in self-similar traffic modeling.
\newblock {\em ACM SIGCOMM Computer Communication Review}, 27(2):5--23, 1997.

\bibitem{Wei}
Ch.H. Wei{\ss}.
\newblock Thinning operations for modeling time series of counts---a survey.
\newblock {\em AStA Adv.\ Stat.\ Anal.}, 92(3):319--341, 2008.

\bibitem{WilTaqSheWil}
W.~Willinger, M.S. Taqqu, R.~Sherman, and D.V. Wilson.
\newblock Self-similarity through high-variability: statistical analysis of
  ethernet lan traffic at the source level.
\newblock {\em IEEE/ACM Transactions on Networking}, 5(1):71--86, 1997.

\bibitem{Zaf2004}
P.~Zaffaroni.
\newblock Contemporaneous aggregation of linear dynamic models in large
  economies.
\newblock {\em J. Econometrics}, 120(1):75--102, 2004.

\bibitem{ZheBasDat2007}
H.~Zheng, I.V. Basawa, and S.~Datta.
\newblock First-order random coefficient integer-valued autoregressive
  processes.
\newblock {\em J. Statist.\ Plann.\ Inference}, 137(1):212--229, 2007.

\bibitem{Zol}
V.M. Zolotarev.
\newblock {\em One-{D}imensional {S}table {D}istributions}, volume~65 of {\em
  Translations of Mathematical Monographs}.
\newblock American Mathematical Society, Providence, RI, 1986.
\newblock Translated from the Russian by H.H. McFaden, Translation edited by
  Ben Silver.

\end{thebibliography}

\end{document}